%% file: main.tex
\documentclass{elsarticle}

\makeatletter
\def\ps@pprintTitle{
  \let\@oddfoot\@empty
}
\makeatother

\input{preamble.pap.tex}

\input{preamble.macros.tex}

\input{preamble.tikz.tex}

\usepackage{rotating}

\usepackage{setspace}
\usepackage{enumitem}
\usepackage{ulem}
\usepackage{tipa}

\newbool{fastcompile}
\setbool{fastcompile}{false}

\begin{document}
\title{Symmetric, optimization-based, cross-element compatible nodal distributions for high-order finite elements}

\author[rvt2]{Julian M. Kaufmann\fnref{fn1}}
\ead{jkaufma2@nd.edu}

\author[rvt2]{Matthew J. Zahr\fnref{fn2}\corref{cor1}}
\ead{mzahr@nd.edu}

\address[rvt2]{Department of Aerospace and Mechanical Engineering, University
	of Notre Dame, Notre Dame, IN 46556, United States}

\cortext[cor1]{Corresponding author}

\fntext[fn1]{Undergraduate Student, Department of Aerospace and Mechanical Engineering, University of Notre Dame}
\fntext[fn2]{Assistant Professor, Department of Aerospace and Mechanical
             Engineering, University of Notre Dame}


\begin{abstract}
We present a general framework to construct symmetric, well-conditioned, cross-element compatible nodal distributions that can be used for high-order and high-dimensional finite elements. Starting from the inherent symmetries of an element geometry, we construct node groups in a systematic and efficient manner utilizing the natural coordinates of each element, while ensuring nodes stay within the elements. Proper constraints on the symmetry group lead to nodal distributions that ensure cross-element compatibility (i.e., nodes of adjacent elements are co-located) on both homogeneous and mixed meshes. The final nodal distribution is defined as a minimizer of an optimization problem over symmetry group parameters with linear constraints that ensure nodes remain with an element and enforce other properties (e.g., cross-element compatibility). We demonstrate the merit of this framework by comparing the proposed optimization-based nodal distributions with other popular distributions available in the literature, and its robustness by generating optimized nodal distributions for otherwise difficult elements (such as simplex and pyramid elements). All nodal distributions are tabulated in the $\mathtt{optnodes}$ package \cite{optnodes}.
\end{abstract}
    
\maketitle

\section{Introduction}
\label{sec:intro}
\input{_intro.tex}

\section{Abstract formulation}
\label{sec:absform}
\input{_absform.tex}

\section{Reference domains}
\label{sec:refdom}
\input{_refdom.tex}

\section{Optimization-based nodal distributions}
\label{sec:optim}
\input{_optim.tex}

\section{Numerical experiments}
\label{sec:num_exp}
\input{_numexp.tex}

\section{Conclusion}
\label{sec:conclude}
\input{_conclude.tex}

\appendix
\section{Tabulated comparisons}
\label{sec:appendix}
\input{_appendix.tex}

\section*{Acknowledgments}
This work is supported by AFOSR award numbers FA9550-20-1-0236, FA9550-22-1-0002, FA9550-22-1-0004, ONR award number N00014-22-1-2299, and NSF award number CBET-2338843. The content of this publication does not necessarily reflect the position or policy of any of these supporters, and no official endorsement should be inferred.

\bibliographystyle{plain}
\bibliography{biblio}

\end{document}

%% file: preamble.pap.tex
\usepackage[margin=1in]{geometry}

\usepackage{bbm}
\usepackage{graphicx}
\usepackage{pstool}
\usepackage{wrapfig}
\usepackage{caption}
\usepackage{subcaption}
\usepackage[section]{placeins} 
\usepackage{fancyhdr} 
\usepackage{lastpage} 
\usepackage{extramarks} 

\usepackage{xcolor} 
\usepackage{enumerate} 
\usepackage{paralist} 
\usepackage{amsmath, amsthm, amssymb, mathtools}
\usepackage{mathabx, pifont, stmaryrd} 
\usepackage[explicit]{titlesec} 
\usepackage{etoolbox} 
\usepackage{bibentry} 
\makeatletter\let\saved@bibitem\@bibitem\makeatother 
\usepackage[colorlinks, bookmarksopen, bookmarksnumbered,
            citecolor=red,urlcolor=red]{hyperref} 
\makeatletter\let\@bibitem\saved@bibitem\makeatother 
\usepackage{comment} 

\usepackage{algorithm}
\usepackage{algorithmic}

\usepackage{optidef}
\usepackage{mathtools}

\newtheorem{remark}{Remark}

\theoremstyle{definition}

%% file: preamble.macros.tex
\usepackage{bm} 
\usepackage{amsfonts} 

\newcommand{\spanV}{\ensuremath{\mathrm{span\hspace{0.05cm}}}}

\newcommand{\optunc}[2]{\underset{#1}{\text{minimize}} ~~ #2}

\newcommand{\Ical}{\ensuremath{\mathcal{I}}}

\newcommand{\Pcal}{\ensuremath{\mathcal{P}}}

\newcommand{\Rcal}{\ensuremath{\mathcal{R}}}

\newcommand{\Cbb}{\ensuremath{\mathbb{C}}}

\newcommand{\Nbb}{\ensuremath{\mathbb{N} }}
\newcommand{\Pbb}{\ensuremath{\mathbb{P} }}
\newcommand{\Qbb}{\ensuremath{\mathbb{Q} }}
\newcommand{\Rbb}{\ensuremath{\mathbb{R} }}

\newcommand\Abm{{\ensuremath{\bm{A}}}}
\newcommand\Bbm{{\ensuremath{\bm{B}}}}

\newcommand\Ibm{{\ensuremath{\bm{I}}}}

\newcommand\Nbm{{\ensuremath{\bm{N}}}}

\newcommand\Sbm{{\ensuremath{\bm{S}}}}

\newcommand\bbm{{\ensuremath{\bm{b}}}}

\newcommand\rbm{{\ensuremath{\bm{r}}}}

\newcommand\vbm{{\ensuremath{\bm{v}}}}
\newcommand\wbm{{\ensuremath{\bm{w}}}}
\newcommand\xbm{{\ensuremath{\bm{x}}}}

\newcommand\lambdabold{{\ensuremath{\boldsymbol{\lambda}}}}

\newcommand\nubold{{\ensuremath{\boldsymbol{\nu}}}}

\newcommand\sigmabold{{\ensuremath{\boldsymbol{\sigma}}}}

\newcommand\xibold{{\ensuremath{\boldsymbol{\xi}}}}



%% file: preamble.tikz.tex
\usepackage{tikz}
\usepackage{pgfplots}
\usepackage{pgfplotstable, booktabs}
\pgfplotsset{compat=1.9}

\usetikzlibrary{pgfplots.groupplots}
\usepgfplotslibrary{fillbetween}
\usetikzlibrary{calc,fit,matrix,arrows,automata,positioning,shapes}
\usetikzlibrary{arrows.meta}

\pgfplotsset{select coords between index/.style 2 args={
    x filter/.code={
        \ifnum\coordindex<#1\fi
        \ifnum\coordindex>#2\fi
    }
}}

\tikzset{
 invisible/.style={opacity=0},
 visible on/.style={alt={#1{}{invisible}}},
 alt/.code args={<#1>#2#3}{%
   \alt<#1>{\pgfkeysalso{#2}}{\pgfkeysalso{#3}}
 },
}


%% file: _intro.tex
It is well-known that Lagrange polynomial bases, commonly used in finite element methods, can quickly become ill-conditioned as the polynomial degree grows if they are based on uniform nodal distributions. On the other hand, the Gauss-Legendre-Lobatto (GLL) nodes, which are well-defined for elements with tensor-product structure, lead to near-optimal interpolation properties \cite{CHEN1995405}. Substantial research has been devoted to the discovery or generation of optimal nodal distributions for other element geometries commonly used as finite elements \cite{fekete2,doi:10.1137/S003614299630587X,chan2015comparison}. The task is complicated by other properties we wish to impose on the nodal distribution including: (1) symmetry, (2) sufficient nodes on each face of the element to uniquely define a polynomial restricted to the face, and (3) nodes of adjacent elements are co-located (cross-element compatibility). The first property ensures the approximation is invariant to rotations of elements within its symmetry group. The second and third properties help enforce $C^0$ continuity of the approximation for $H^1$-conforming finite elements.

Near-optimal nodal distributions have been well-studied for simplex elements. Explicitly defined nodal distributions are directly constructed from a formula or algorithm \cite{Isaac,blyth2006lobatto,luo2006lobatto,warburton2006explicit,doi:10.1137/S003614299630587X,doi:10.1137/S1064827598343723}, whereas implicitly defined distributions involve the solution of an optimization problem \cite{Fekete_def,CHEN1995405,taylor2000algorithm,rapetti2012generation}. Explicit distributions are simple to implement and efficient to define on-the-fly for any polynomial degree. On the other hand, implicit distributions are expensive to compute and must be tabulated, but usually have lower Lebesgue constants than explicit distributions. The Fekete nodes \cite{Fekete_def,taylor2000algorithm} maximize the determinant of the Vandermonde matrix, coincide with the asymptomatically optimal GLL nodes on the d-dimensional cube, reduce to the GLL nodes in one dimension and along edges of triangles, and exhibit logarithmic growth of the Lebesgue constant with polynomial degree. Nodal distributions based on minimizing the Lebesgue constant \cite{roth2005nodal,rapetti2012generation} have been shown to the yield the slowest growth of the Lebesgue constant; however, they are not necessarily symmetric. The explicit distributions of \cite{Isaac} were shown to be competitive with implicit distributions, particularly for tetrahedron where generating high-order implicit distributions is difficult and expensive. Several of these concepts have been extended to pyramid elements \cite{gassner2009polymorphic,bergot2010higher, bos2011geometric, chan2015comparison}, which commonly arise in hexahedra-dominated meshes; however, pyramids have received much less attention than simplices.

We propose a numerical framework to create symmetric, optimization-based nodal distributions that satisfy arbitrary linear constraints on the symmetry group parameters. Constraints will be used to ensure sufficient nodes lie on the element faces to uniquely define a polynomial of degree $p$ (global polynomial approximation degree) on the face. Linear constraints will optionally be used to ensure the distribution on a given face matches a given symmetric nodal distribution of the lower dimensional element. This is accomplished by defining nodal distributions from collections of symmetry orbits and optimizing the coefficients of the symmetry orbits. The Lebesgue constant is used as the optimization objective; however, the framework could be used with other objective functions to generate distributions for purposes other than interpolation. Our framework generalizes the approach of \cite{CHEN1995405} that only considered triangles and the work in \cite{rapetti2012generation} that does not consider optimization within symmetry groups. A similar approach of optimization within symmetry orbits has been used to generate quadrature nodes \cite{WANDZURAT20031829,WITHERDEN20151232}.  The novelty of our work lies in the extension of this idea to interpolation nodes, which requires a different objective function and a framework for handling constraints within a symmetry orbit (such as cross-element compatibility). Nodal distributions for the line (up to polynomial degree 30), the triangle and quadrilateral (up to polynomial degree 23), and the tetrahedron, hexahedron, triangular prism, and pyramid (up to polynomial degrees 9) are generated to show the merit of this framework.



%% file: _absform.tex
In this section, we detail an abstract, systematic procedure to construct symmetric nodal distributions that respect a general set of constraints that will be used to enforce cross-element compatibility. Specifically, we introduce the notion of an intrinsic symmetry orbit for an abstract element geometry that maps symmetry parameters to Cartesian coordinates of the points in the symmetry orbit (Section~\ref{sec:absform:natcoord}-\ref{sec:absform:symorb}). We equip these symmetry orbits with linear constraints (Section~\ref{sec:absform:csymorb}) to ensure nodal distributions can be constructed that satisfy requirements imposed by the finite element method. Finally, finite element nodal distributions are defined as appropriate collections of constrained symmetry orbits (Section~\ref{sec:absform:symorbcoll}). The abstract formulation is detailed for seven standard finite element geometries in Section~\ref{sec:refdom}.


\subsection{Natural coordinates}
\label{sec:absform:natcoord}
Consider an element geometry $\Omega \subset \Rbb^d$.
Nodal distributions over $\Omega$ are conveniently expressed in Cartesian coordinates, denoted $\xbm \in \Omega$, whereas symmetry orbits, which leverage the inherent structure of the element, may be naturally described in a different coordinate system. We introduce the natural coordinates of an element, denoted $\lambdabold \in \Rbb^{d'}$, that will be used to define symmetry orbits (Section~\ref{sec:absform:symorb}), and assume the natural and Cartesian coordinates are related by a linear mapping, $\Rcal : \Cbb_\lambda \rightarrow \Omega$,
\begin{equation} \label{eqn:nat_def}
\Rcal: \lambdabold \mapsto \Nbm \lambdabold + \nubold,
\end{equation}
where $\Nbm \in \Rbb^{d'\times d}$ and $\nubold \in \Rbb^d$ define the linear mapping. The mapping domain, $\Cbb_\lambda \subset\Rbb^{d'}$, is the intersection of $2c$ linear constraints to ensure the range of $\Rcal$ is $\Omega$
\begin{equation} \label{eqn:nat_bnd}
    \Cbb_\lambda \coloneqq \{\lambdabold\in\Rbb^{d'} \mid \vbm_\mathrm{l} \leq \Bbm_\lambda \lambdabold \leq \vbm_\mathrm{u}\},
\end{equation}
where $\vbm_\mathrm{l},\vbm_\mathrm{u}\in\Rbb^c$ and $\Bbm_\lambda \in \Rbb^{c\times d'}$. Section~\ref{sec:refdom} details the natural coordinates for each element geometry considered.


\subsection{Intrinsic symmetry orbits}
\label{sec:absform:symorb}
We will construct finite element nodal distributions from symmetry orbits to obtain distributions with symmetries that match those of the element. An element $\Omega$ has $M$ symmetry orbits, where each symmetry orbit maps a set of parameters into a collection of natural coordinates. Specifically, symmetry orbit $k$ of $\Omega$, denoted $\bar\Pcal_k : \Cbb_{\xi,k} \rightarrow \Cbb_\lambda^{m_k}$, maps a collection of $l_k$ symmetry parameters, denoted $\xibold_k \in \Rbb^{l_k}$, to a collection of $m_k$ natural coordinates, denoted $\bar{\lambdabold} = \{ \lambdabold_{k,1}, \dots, \lambdabold_{k,m_k}\} \in \Cbb_\lambda^{m_k}$, corresponding to all points produced by the symmetry orbit. The symmetry orbits are written as
\begin{equation} \label{eqn:symorb_def_all}
    \bar\Pcal_k : \xibold_k \mapsto \{ \Pcal_{k,1}(\xibold_k), \dots, \Pcal_{k,m_k}(\xibold_k) \},
\end{equation}
where $\Pcal_{k,i} : \Cbb_{\xi,k} \rightarrow \Cbb_\lambda$, the linear mapping defining the $i$th point of the $k$th symmetry orbit ($\lambdabold_{k,i}$), is defined as
\begin{equation}\label{eqn:symorb_def}
  \Pcal_{k,i} : \xibold_k \mapsto \Sbm_{k,i} \xibold_k + \sigmabold_{k,i}
\end{equation}
with $\Sbm_{k,i}\in\Rbb^{d' \times l_k}$ and $\sigmabold_{k,i}\in\Rbb^{d'}$.
The mapping domain, $\Cbb_{\xi,k} \subset \Rbb^{l_k}$, is the intersection of $2c$ linear constraints to ensure the range of $\bar\Pcal_k$ is $\Cbb_\lambda$
\begin{equation} \label{eqn:symorb_bnd}
    \Cbb_{\xi,k} \coloneqq \{\xibold\in\Rbb^{l_k} \mid \wbm_{\mathrm{l},k} \leq \Bbm_{\xi,k} \xibold \leq \wbm_{\mathrm{u},k}\},
\end{equation}
where $\wbm_{\mathrm{l},k},\wbm_{\mathrm{u},k}\in\Rbb^c$ and $\Bbm_{\xi,k} \in \Rbb^{c\times l_k}$. The bounds on $\xibold_k$ can be written explicitly from the bounds on the natural coordinates (\ref{eqn:nat_bnd}) and the linear mapping from symmetry-to-natural coordinates (\ref{eqn:symorb_def}) as
\begin{equation} \label{eqn:symorb_bnd2}
    \Bbm_{\xi,k} = \Bbm_\lambda \Sbm_{k,1}, \qquad
    \wbm_{\mathrm{l},k} = \vbm_\mathrm{l} - \Bbm_\lambda\sigmabold_{k,1}, \qquad
    \wbm_{\mathrm{u},k} = \vbm_\mathrm{u} - \Bbm_\lambda\sigmabold_{k,1}.
\end{equation}
The linear mapping of the first point only ($\Sbm_{k,1}$, $\sigmabold_{k,1}$) is used to defined the bounds on $\xibold_k$. By symmetry, all other points produced by the symmetry orbit will respect the bound on the natural coordinates.
Section~\ref{sec:refdom} details the symmetry orbits for each element geometry considered. Additionally, Sections~\ref{sec:refdom:tri}~and~\ref{sec:refdom:pyrmd} give a concrete examples of this structure.


\subsection{Constrained symmetry orbits}
\label{sec:absform:csymorb}
Nodal distributions for $H^1$-conforming finite element must satisfy a number of conditions not accounted for by the symmetry orbits. For example, a node should lie at each element vertex and an appropriate number of nodes should lie on each element edge and face to uniquely define a polynomial of appropriate degree when restricted to those entities (to ensure continuity). To ensure these and other necessary constraints will be satisfied, we incorporate $2a_k$ additional linear constraints on the symmetry variables of the form
\begin{equation} \label{eqn:symcoord_addbnd}
    \Cbb_{\xi,k}^{\Abm_k,\bbm_{\mathrm{l},k},\bbm_{\mathrm{u},k}} \coloneqq \{ \xibold \in \Rbb^{m_k} \mid \bbm_{\mathrm{l},k} \leq \Abm_k \xibold_k \leq \bbm_{\mathrm{u},k} \},
\end{equation}
where $\Abm_k \in \Rbb^{a_k\times m_k}$ and $\bbm_{\mathrm{l},k},\bbm_{\mathrm{u},k} \in \Rbb^{a_k}$ define the linear constraints. The additional constraints are enforced by using the symmetry orbits $\bar\Pcal_k$ with the domain $\Cbb_{\xi,k}\cap\Cbb_{\xi,k}^{\Abm_k,\bbm_{\mathrm{l},k},\bbm_{\mathrm{u},k}}$ in place of the domain introduced in Section~\ref{sec:absform:symorb} ($\Cbb_{\xi,k}$).


\subsection{Symmetry orbit collections}
\label{sec:absform:symorbcoll}
Finally, we define a parametrized nodal distribution $\Nbb_{\bar\xi} \subset \Rbb^d$ as the union of the nodes produced by a collection of $n$ symmetry orbits. Let $\Ical = (k_1,\dots,k_n)$ with $k_j \in \{1,\dots,M\}$ define the symmetry orbits used and $\{(\Abm_{k_j}^{(j)},\bbm_{\mathrm{l},k_j}^{(j)},\bbm_{\mathrm{u},k_j}^{(j)})\}_{j=1}^n$ denote the constraints imposed on each symmetry orbit in the collection. The parameter vector, $\bar\xibold \in \bar\Cbb_\xi$, is defined as
\begin{equation} \label{eqn:parvec}
    \bar\xibold = (\xibold_{k_1}^{(1)}, \dots, \xibold_{k_n}^{(n)}),
\end{equation}
where $\xibold_{k_j}^{(j)} \in \Cbb_{\xi,k_j} \cap \Cbb_{\xi,k_j}^{\Abm_{k_j}^{(j)},\bbm_{\mathrm{l},k_j}^{(j)},\bbm_{\mathrm{u},k_j}^{(j)}}$ are the symmetry parameters used for the $j$th symmetry orbit in the collection (corresponds to symmetry orbit $k_j$) and
\begin{equation}
    \bar\Cbb_\xi \coloneqq \left(\Cbb_{\xi,k_1} \cap \Cbb_{\xi,k_1}^{\Abm_{k_1}^{(1)},\bbm_{\mathrm{l},k_1}^{(1)},\bbm_{\mathrm{u},k_1}^{(1)}}\right) \times \cdots \times \left(\Cbb_{\xi,k_n} \cap \Cbb_{\xi,k_n}^{\Abm_{k_n}^{(n)},\bbm_{\mathrm{l},k_n}^{(n)},\bbm_{\mathrm{u},k_n}^{(n)}}\right).
\end{equation}
The set $\bar\Cbb_\xi$ is the Cartesian product of the intersection of linear constraints, which makes it an intersection of linear constraints in the larger space.
From these definitions, the nodal distribution is written as
\begin{equation}
    \Nbb_{\bar\xibold} = \bigcup_{j=1}^n \bar\Pcal_{k_j}(\xibold_{k_j}^{(j)}),
\end{equation}
where $\bar\xibold \in \bar\Cbb_\xi$. The cardinality of $\Nbb_{\bar\xibold}$ is $\bar{m} = \sum_{k\in\Ical} m_k$. We provide concrete examples of this framework in Section~\ref{sec:refdom} and detail systematic procedures to build constrained symmetry orbit collections that lead to admissible finite elements in Section~\ref{sec:optim}.

\begin{remark}
    \color{red}{All nodal distributions are generated in the reference element for each geometry, which is straight-sided and lies in the bi-unit block (Section~\ref{sec:refdom}). These nodal distribution can be mapped to (potentially curved) physical elements in a finite element mesh using the isoparametric concept to generate the corresponding nodal distribution over arbitrary curved elements of the same geometry.}
\end{remark}

\begin{remark}
    The proposed procedure to define nodal distributions guarantees the nodes will be symmetric with respect to the natural symmetries of the element geometry. However, it does not preserve other structure. For example, it does not guarantee that the nodal distribution will have tensor product structure, even for tensor product elements (e.g., quadrilaterals, hexahedron). However, this structure could be imposed through additional linear constraints on the symmetry orbits (Section~\ref{sec:absform:csymorb}). Tensor product structure is not imposed on any of the nodal distributions constructed in this work.
\end{remark}

%% file: _refdom.tex
In this section, we specialize the abstract formulation of Section~\ref{sec:absform} to the seven element geometries considered in this work. We detail the element geometry, natural coordinates, symmetry orbits, and corresponding bounds from \cite{WITHERDEN20151232}. The mappings, $\Rcal$ and $\bar\Pcal_k$, are explicitly defined, from which the linear mapping terms ($\Nbm$, $\nubold$, $\Sbm_{k,i}$, $\sigmabold_{k,i}$) can be determined. An example of the entire abstract framework in Section~\ref{sec:absform} is provided for the triangular element.

\subsection{Line}
Despite existence of highly effective nodal distribution for line elements (i.e., the GLL nodes), we use the line as a simple demonstration of the proposed optimization-based nodal distributions. We take the domain to be $\Omega = (-1, 1)$ and select the Cartesian coordinates as the natural coordinates, i.e., $d = d' = 1$ with $\Rcal(\lambdabold) = \lambdabold$, which requires $\Bbm_\lambda = 1$, $\vbm_\mathrm{l} = -1$, and $\vbm_\mathrm{u} = 1$.
The symmetry orbits for the line element are
\begin{equation}
    \begin{aligned}
        \bar\Pcal_1 &= \{0\}, \qquad &m_1 &= 1 \\
        \bar\Pcal_2(\alpha) &= \{\alpha,~-\alpha\}, \qquad &m_2 &= 2.
    \end{aligned}
\end{equation}
The local function space associated with the line element is the standard polynomial space of maximum degree $p$, $\Qbb^p_1 \coloneqq \spanV\{1, x, x^2, \dots, x^p\}$. The total number of nodes for polynomial degree $p$ is given by $p+1$.

\subsection{Triangle}
\label{sec:refdom:tri}
\begin{figure}
    \centering
    \input{_py/tridemo_symorb.tikz}
    \caption{Symmetry orbits for the triangular element. Legend: Points generated by $\bar\Pcal_1$ (\ref{line:tridemo1:orb1}), $\bar\Pcal_2(0.2)$ (\ref{line:tridemo1:orb2}), and $\bar\Pcal_3(0.3,0.2)$ (\ref{line:tridemo1:orb3}).}
    \label{fig:tri_symorbit}
\end{figure}
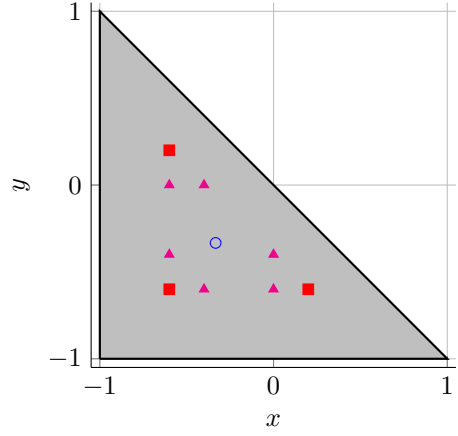

The domain for the triangular element is the bi-unit triangle,
$\Omega = \{(x,y) \mid x, y \geq -1, x+y \leq 0\}$. The natural coordinates are taken to be the barycentric coordintes with natural-to-Cartesian mapping given by $\Rcal(\lambdabold) = (-\lambda_1+\lambda_2-\lambda_3, -\lambda_1-\lambda_2+\lambda_3)$ and bounds
\begin{equation}
    \Bbm_\lambda = \begin{bmatrix} 1 & 0 & 0 \\ 0 & 1 & 0 \\ 0 & 0 & 1 \\ 1 & 1 & 1\end{bmatrix}, \qquad
    \vbm_\mathrm{l} = \begin{bmatrix} 0 \\ 0 \\ 0 \\ 1 \end{bmatrix}, \qquad
    \vbm_\mathrm{u} = \begin{bmatrix} 1 \\ 1 \\ 1 \\ 1 \end{bmatrix}.
\end{equation}
The symmetry orbits are given in terms of the barycentric coordinates as
\begin{equation} \label{eqn:tri_symorbs}
    \begin{aligned}
        \bar\Pcal_1 &= \{(1/3,~1/3,~1/3)\}, \qquad &m_1 &= 1 \\
        \bar\Pcal_2(\alpha) &= \text{Perm}(\alpha,\alpha,1-2\alpha), \qquad &m_2 &= 3 \\
        \bar\Pcal_3(\alpha,\beta) &= \text{Perm}(\alpha,\beta,1-\alpha-\beta), \qquad &m_3 &= 6, \\
    \end{aligned}
\end{equation}
where $\text{Perm}(a,b,c)$ returns a set of 3-vectors containing the unique permutations of the inputs (Figure~\ref{fig:tri_symorbit}). The local function space associated with the triangular element is the standard polynomial space of degree $p$, $\Pbb_2^p \coloneqq \spanV\{x^\alpha y^\beta \mid 0 \leq \alpha+\beta\leq p\}$. The total number of nodes for polynomial degree $p$ is given by $\frac{(p+1)(p+2)}{2}$.

We close this subsection with a concrete example of the complete framework of Section~\ref{sec:absform}. To this end, we select a collection consisting of five symmetry orbits, $\Ical = (1, 2, 2, 3, 3)$ and write the parameter vector $\bar\xibold$ according to (\ref{eqn:parvec}) with
\begin{equation}
    \xibold_1^{(1)} = (), \quad
    \xibold_2^{(2)} = \alpha_1, \quad
    \xibold_2^{(3)} = \alpha_2, \quad
    \xibold_3^{(4)} = (\alpha_3,\alpha_4), \quad
    \xibold_3^{(5)} = (\alpha_5,\alpha_6).
\end{equation}
Then, the explicit symmetry orbit formulation for $\xibold_3^{(4)}$ is given by
\begin{equation}
        \bar\Pcal_3 : \xibold_3^{(4)} \mapsto \{ \Pcal_{3,1}(\xibold_3^{(4)}), \Pcal_{3,2}(\xibold_3^{(4)}), \Pcal_{3,3}(\xibold_3^{(4)}),
        \Pcal_{3,4}(\xibold_3^{(4)}), \Pcal_{3,5}(\xibold_3^{(4)}), \Pcal_{3,6}(\xibold_3^{(4)})\},
\end{equation}
where
\begin{gather}
      \Pcal_{3,1}(\xibold_3^{(4)}) = \begin{bmatrix} \alpha_3 \\ \alpha_4 \\ 1-\alpha_3-\alpha_4 \end{bmatrix}, 
      \quad \Pcal_{3,2}(\xibold_3^{(4)}) = \begin{bmatrix}  \alpha_4 \\ \alpha_3 \\ 1-\alpha_3-\alpha_4 \end{bmatrix}, 
      \quad \Pcal_{3,3}(\xibold_3^{(4)}) = \begin{bmatrix}  \alpha_4 \\ 1-\alpha_3-\alpha_4 \\ \alpha_3  \end{bmatrix}, \\
      \quad \Pcal_{3,4}(\xibold_3^{(4)}) = \begin{bmatrix}  \alpha_3 \\ 1-\alpha_3-\alpha_4 \\ \alpha_4  \end{bmatrix},  
      \quad \Pcal_{3,5}(\xibold_3^{(4)}) = \begin{bmatrix}  1-\alpha_3-\alpha_4 \\ \alpha_4 \\ \alpha_3  \end{bmatrix},  
      \quad \Pcal_{3,6}(\xibold_3^{(4)}) = \begin{bmatrix} 1-\alpha_3-\alpha_4 \\ \alpha_3 \\  \alpha_4  \end{bmatrix}.
\end{gather}
Furthermore, the first point in each symmetry orbits in (\ref{eqn:tri_symorbs}) is expressed using the formalism of Section~\ref{sec:absform} as
\begin{equation}
    \Sbm_{1,1} = \emptyset, \quad \sigmabold_{1,1} = \begin{bmatrix} 1/3 \\ 1/3 \\ 1/3 \end{bmatrix}, \quad
    \Sbm_{2,1} = \begin{bmatrix} 1 \\ 1 \\ -2 \end{bmatrix}, \quad \sigmabold_{2,1} = \begin{bmatrix} 0 \\ 0 \\ 1 \end{bmatrix}, \quad
    \Sbm_{3,1} = \begin{bmatrix} 1 & 0 \\ 0 & 1 \\ -1 & -1\end{bmatrix}, \quad \sigmabold_{3,1} = \begin{bmatrix} 0 \\ 0 \\ 1 \end{bmatrix}.
\end{equation}

We take the first four orbits to be standard, i.e., the constraint matrix and vector are taken as $\Abm_{\Ical_j}^{(j)} = \bbm_{\mathrm{l},\Ical_j}^{(j)} = \bbm_{\mathrm{u},\Ical_j}^{(j)} = \emptyset$ for $j=1,...,4$, and impose the condition that $\alpha_6 = 0$ on the last orbit with
\begin{equation}
  \Abm_3^{(5)} = \begin{bmatrix} 0 & 1\end{bmatrix}, \qquad \bbm_{3,\mathrm{l}}^{(5)} = \bbm_{3,\mathrm{u}}^{(5)} = 0.
\end{equation}

\begin{figure}
    \centering
    \input{_py/tridemo_nddist.tikz}
    \caption{Nodal distribution for symmetry orbit collection described in Section~\ref{sec:refdom:tri} with $\alpha_1 = 0.25$, $\alpha_2 = 0.5$, $\alpha_3 = 0.1$, $\alpha_4 = 0.6$, $\alpha_5 = 0.7$, and $\alpha_6 = 0$. Legend: Points generated by $\bar\Pcal_1$ (\ref{line:tridemo2:orb1}), $\bar\Pcal_2(\alpha_1)$ (\ref{line:tridemo2:orb2a}), $\bar\Pcal_2(\alpha_2)$ (\ref{line:tridemo2:orb2b}), $\bar\Pcal_3(\alpha_3,\alpha_4)$ (\ref{line:tridemo2:orb3a}),
    and $\bar\Pcal_3(\alpha_5,\alpha_6)$ (\ref{line:tridemo2:orb3b}).}
    \label{fig:tri_nddist0}
\end{figure}
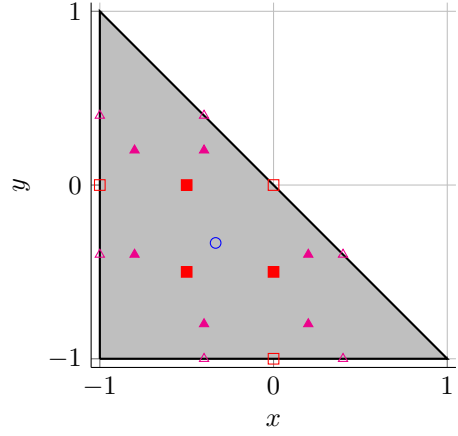

According to (\ref{eqn:symorb_bnd})-(\ref{eqn:symorb_bnd2}), this leads to the following constraints on $\alpha_1,\dots,\alpha_6$
\begin{equation}
    \begin{bmatrix} 0 \\ 0 \\ -1 \\ 0 \end{bmatrix} \leq
    \begin{bmatrix} 1 \\ 1 \\ -2 \\ 0 \end{bmatrix}\alpha_i \leq
    \begin{bmatrix} 1 \\ 1 \\ 0 \\ 0 \end{bmatrix}, \qquad
    \begin{bmatrix} 0 \\ 0 \\ -1 \\ 0 \end{bmatrix} \leq
    \begin{bmatrix} 1 & 0 \\ 0 & 1 \\ -1 & -1 \\ 0 & 0 \end{bmatrix}\begin{bmatrix} \alpha_j \\ \alpha_k\end{bmatrix} \leq
    \begin{bmatrix} 1 \\ 1 \\ 0 \\ 0 \end{bmatrix},
\end{equation}
where $i = 1,2$, $j=3,5$, and $k=4,6$. After removing redundant and trivial constraints and incorporating the $\alpha_6 = 0$ constraint from above, we have
\begin{equation}
    0 \leq \alpha_1, \alpha_2 \leq 1/2, \quad
    0 \leq \alpha_3, \alpha_4, \alpha_5 \leq 1, \quad \alpha_3+\alpha_4 \leq 1, \quad \alpha_6 = 0.
\end{equation}
The nodal distribution for the $\alpha_1 = 0.25$, $\alpha_2 = 0.5$, $\alpha_3 = 0.1$, $\alpha_4 = 0.6$, $\alpha_5 = 0.7$, and $\alpha_6 = 0$ is shown in Figure~\ref{fig:tri_nddist0}.

\subsection{Quadrilateral}
The domain for the quadrilateral element is the bi-unit quadrilateral,
$\Omega = \{(x,y) \mid -1 \leq x,y \leq 1\}$. The natural coordinates are, as in the line, Cartesian coordinates with $\Rcal(\lambdabold) = (\lambda_1, \lambda_2)$ and bounds
\begin{equation}
    \Bbm_\lambda = \begin{bmatrix} 1 & 0 \\ 0 & 1 \end{bmatrix}, \qquad
    \vbm_\mathrm{l} = \begin{bmatrix} -1 \\ -1 \end{bmatrix}, \qquad
    \vbm_\mathrm{u} = \begin{bmatrix} 1 \\ 1  \end{bmatrix}.
\end{equation}
The symmetry orbits are given as
\begin{equation} \label{eqn:quad_symorb}
    \begin{aligned}
        \bar\Pcal_1 &= \{(0,0)\}, \qquad &m_1 &= 1 \\
        \bar\Pcal_2(\alpha) &= \chi(\alpha,0), \qquad &m_2 &= 4 \\
        \bar\Pcal_3(\alpha) &= \chi(\alpha,\alpha), \qquad &m_3 &= 4, \\
        \bar\Pcal_4(\alpha,\beta) &= \chi(\alpha,\beta), \qquad &m_4 &= 8, \\
    \end{aligned}
\end{equation}
where $\chi(a,b)$ returns a set of 2-vectors containing the unique permutations of the inputs and their additive inverse.
The local function space associated with the quadrilateral element is the standard polynomial space of maximal degree $p$, $\Qbb^p_2 \coloneqq \spanV\{x^\alpha y^\beta \mid 0 \leq \alpha, \beta\leq p\}$. The total number of nodes for polynomial degree $p$ is given by $(p+1)^2$.

\subsection{Tetrahedron}
The domain for the tetrahedral element is the bi-unit tetrahedron,
$\Omega = \{(x,y,z) \mid x, y, z \geq -1, x+y+z \leq -1\}$. The natural coordinates are taken to be the barycentric coordintes with natural-to-Cartesian mapping given by $\Rcal(\lambdabold) = (-\lambda_1+\lambda_2-\lambda_3-\lambda_4, -\lambda_1-\lambda_2+\lambda_3-\lambda_4, -\lambda_1-\lambda_2-\lambda_3+\lambda_4)$ and bounds
\begin{equation}
    \Bbm_\lambda = \begin{bmatrix} 1 & 0 & 0 & 0 \\ 0 & 1 & 0 & 0\\ 0 & 0 & 1 & 0\\ 0 & 0 & 0 & 1\\ 1 & 1 & 1 & 1\end{bmatrix}, \qquad
    \vbm_\mathrm{l} = \begin{bmatrix} 0 \\ 0 \\ 0 \\ 0 \\1 \end{bmatrix}, \qquad
    \vbm_\mathrm{u} = \begin{bmatrix} 1\\1 \\ 1 \\ 1 \\ 1 \end{bmatrix}.
\end{equation}
The symmetry orbits are given in terms of the barycentric coordinates as
\begin{equation} \label{eqn:tri_symorb}
    \begin{aligned}
        \bar\Pcal_1 &= \{(1/4,~1/4,~1/4,~1/4)\}, \qquad &m_1 &= 1 \\
        \bar\Pcal_2(\alpha) &= \text{Perm}(\alpha,\alpha,\alpha, 1-3\alpha), \qquad &m_2 &= 4 \\
        \bar\Pcal_3(\alpha) &= \text{Perm}(\alpha,\alpha,1/2-\alpha,1/2-\alpha), \qquad &m_3 &= 6, \\
        \bar\Pcal_4(\alpha,\beta) &= \text{Perm}(\alpha,\alpha,\beta,1-2\alpha-\beta), \qquad &m_4 &= 12, \\
        \bar\Pcal_5(\alpha,\beta,\gamma) &= \text{Perm}(\alpha,\beta,\gamma,1-\alpha-\beta-\gamma), \qquad &m_5 &= 24, \\
    \end{aligned}
\end{equation}
where $\text{Perm}(a,b,c,d)$ returns a set of 4-vectors containing the unique permutations of the inputs. The local function space associated with the triangular element is the standard polynomial space of degree $p$, $\Pbb^p_3 \coloneqq \spanV\{x^\alpha y^\beta z^\gamma \mid 0 \leq \alpha+\beta+\gamma\leq p\}$. The total number of nodes for polynomial degree $p$ is given by $\frac{(p+1)(p+2)(p+3)}{6}$.

\subsection{Hexahedron}
The domain for the hexahedral element is the bi-unit hexahedron,
$\Omega = \{(x,y,z) \mid -1 \leq x,y,z \leq 1\}$. The natural coordinates are, again, Cartesian coordinates with $\Rcal(\lambdabold) = (\lambda_1, \lambda_2, \lambda_3)$ and bounds
\begin{equation}
    \Bbm_\lambda = \begin{bmatrix} 1 & 0 & 0 \\ 0 & 1 & 0 \\ 0 & 0 & 1 
    \end{bmatrix}, \qquad
    \vbm_\mathrm{l} = \begin{bmatrix} -1 \\ -1 \\ -1 \end{bmatrix}, \qquad
    \vbm_\mathrm{u} = \begin{bmatrix} 1 \\ 1  \\ 1 \end{bmatrix}.
\end{equation}
The symmetry orbits are given as
\begin{equation} \label{eqn:quad_symorb}
    \begin{aligned}
        \bar\Pcal_1 &= \{(0,0,0)\}, \qquad &m_1 &= 1 \\
        \bar\Pcal_2(\alpha) &= \chi(\alpha,0,0), \qquad &m_2 &= 6 \\
        \bar\Pcal_3(\alpha) &= \chi(\alpha,\alpha,\alpha), \qquad &m_3 &= 8, \\
        \bar\Pcal_4(\alpha) &= \chi(\alpha,\alpha,0), \qquad &m_4 &= 12, \\
        \bar\Pcal_5(\alpha,\beta) &= \chi(\alpha,\beta,0), \qquad &m_5 &= 24, \\
        \bar\Pcal_6(\alpha,\beta) &= \chi(\alpha,\alpha,\beta), \qquad &m_6 &= 24, \\
        \bar\Pcal_7(\alpha,\beta,\gamma) &= \chi(\alpha,\beta,\gamma), \qquad &m_7 &= 48, \\
    \end{aligned}
\end{equation}
where $\chi(a,b,c)$ returns a set of 3-vectors containing the unique permutations of the inputs and their additive inverse.
The local function space associated with the quadrilateral element is the standard polynomial space of maximal degree $p$, $\Qbb^p_3 \coloneqq \spanV\{x^\alpha y^\beta z^\gamma \mid 0 \leq \alpha, \beta, \gamma \leq p\}$. The total number of nodes for polynomial degree $p$ is given by $(p+1)^3$.

\subsection{Triangular prism}
The domain for the triangular prism is the bi-unit triangular prism,
$\Omega = \{(x,y,z) \mid x, y, z \geq -1, x+y \leq 0, z \leq 1\}$. The natural coordinates are taken to be the barycentric coordinates in the triangle dimension and Cartesian in the extruded dimension with natural-to-Cartesian mapping given by $\Rcal(\lambdabold) = (-\lambda_1+\lambda_2-\lambda_3, -\lambda_1-\lambda_2+\lambda_3, \lambda_4)$ and bounds
\begin{equation}
    \Bbm_\lambda = \begin{bmatrix} 1 & 0 & 0 & 0 \\ 0 & 1 & 0 & 0\\ 0 & 0 & 1 & 0\\ 1 & 1 & 1 & 0\\ 0 & 0 & 0 & 1\end{bmatrix}, \qquad
    \vbm_\mathrm{l} = \begin{bmatrix} 0 \\ 0 \\ 0 \\ 1 \\-1 \end{bmatrix}, \qquad
    \vbm_\mathrm{u} = \begin{bmatrix} 1 \\1 \\ 1 \\ 1 \\ 1 \end{bmatrix}.
\end{equation}
The symmetry orbits are given in terms of these extruded barycentric coordinates as
\begin{equation} \label{eqn:tri_symorb}
    \begin{aligned}
        \bar\Pcal_1 &= \{(1/3,~1/3,~1/3, 0)\}, \qquad &m_1 &= 1 \\
        \bar\Pcal_2(\gamma) &= (1/3,~1/3,~1/3,\pm\gamma), \qquad &m_2 &= 2 \\
        \bar\Pcal_3(\alpha) &= \text{Perm}_3(\alpha,\alpha,1-2\alpha,0), \qquad &m_3 &= 3, \\
        \bar\Pcal_4(\alpha,\gamma) &= \text{Perm}(\alpha,\alpha,1-2\alpha,\pm\gamma), \qquad &m_4 &= 6, \\
        \bar\Pcal_5(\alpha,\beta) &= \text{Perm}_3(\alpha,\beta,1-\alpha-\beta, 0), \qquad &m_5 &= 6, \\
        \bar\Pcal_6(\alpha,\beta,\gamma) &= \text{Perm}_3(\alpha,\beta,1-\alpha-\beta, \pm\gamma), \qquad &m_6 &= 12, \\
    \end{aligned}
\end{equation}
where $\text{Perm}_3(a,b,c,d)$ returns a set of 4-vectors containing the unique permutations of the first three inputs. The local function space associated with the triangular prism element is the extruded polynomial space of degree $p$, $\Pbb^p_2 \Qbb^p_1 \coloneqq \spanV\{x^\alpha y^\beta z^\gamma \mid 0 \leq \alpha+\beta,\gamma\leq p\}$. The total number of nodes for polynomial degree $p$ is given by $\frac{(p+1)^2(p+2)}{2}$.

\subsection{Square Pyramid}
\label{sec:refdom:pyrmd}
The domain for the square pyramid is the bi-unit square pyramid,
$\Omega = \{(x,y,z) \mid -1 \leq z \leq 1, (-1+z)/2 \leq x,y \leq (1-z)/2\}$. The natural coordinates are taken to be the Cartesian coordinates with $\Rcal(\lambdabold) = (\lambda_1, \lambda_2, \lambda_3)$ and bounds

\begin{equation}
    \Bbm_\lambda = \begin{bmatrix} 2 & 0 & 1 \\ 2 & 0 & -1 \\ 0 & 2 & 1\\ 0 & 2 & -1 \\ 0 & 0 & 1 \end{bmatrix}, \qquad
    \vbm_\mathrm{l} = \begin{bmatrix} -\infty \\ -1 \\ -\infty \\ -1 \\ -1 \end{bmatrix}, \qquad
    \vbm_\mathrm{u} = \begin{bmatrix} 1 \\ \infty \\ 1 \\ \infty \\ 1 \end{bmatrix}.
\end{equation}

The symmetry orbits are given in terms of Cartesian coordinates as

\begin{equation} \label{eqn:quad_pyr_symorb}
    \begin{aligned}
        \bar\Pcal_1(\gamma) &= \{(0,0,\gamma)\}, \qquad &m_1 &= 1 \\
        \bar\Pcal_2(\alpha,\gamma) &= \chi(\alpha,0,\gamma), \qquad &m_2 &= 4 \\
        \bar\Pcal_3(\alpha,\gamma) &= \chi(\alpha,\alpha,\gamma), \qquad &m_3 &= 4, \\
        \bar\Pcal_4(\alpha,\beta,\gamma) &= \chi(\alpha,\beta,\gamma), \qquad &m_4 &= 8, \\
    \end{aligned}
\end{equation}
where $\chi(a,b,c)$ returns a set of 3-vectors containing the unique permutations of the first two inputs and their additive inverse. The local function space associated with the pyramidal element is a rational space (polynomial trace space) defined as: $\spanV\{P_{ijk}(x,y,z) \mid 0 \leq i,j \leq p, 0 \leq k \leq p - \max(i,j)\}$. Here, $P_{ijk} : \Rbb^3 \rightarrow \Rbb$ is defined as
\begin{equation}
P_{ijk} : (x,y,z) \mapsto P_i^{0,0} \left( \frac{x}{1-\frac{z+1}{2}} \right) P_j^{0,0} \left( \frac{y}{1-\frac{z+1}{2}} \right) \left(1-\frac{z+1}{2}\right)^c P_k^{2(c+1),0} (z),
\end{equation}
where 
\begin{equation}
    c = \max(i,j), \quad 0 \leq i,j \leq p, \quad 0 \leq k \leq p-c
\end{equation}
and $P_n^{a,b}(x)$ is the Jacobi polynomial of order $n$ with respect to the weight $(1-x)^a(1-x)^b$.
The total number of nodes for degree $p$ is given by $\frac{\left(p+1\right)\left(p+2\right)\left(2p+3\right)}{6}$.

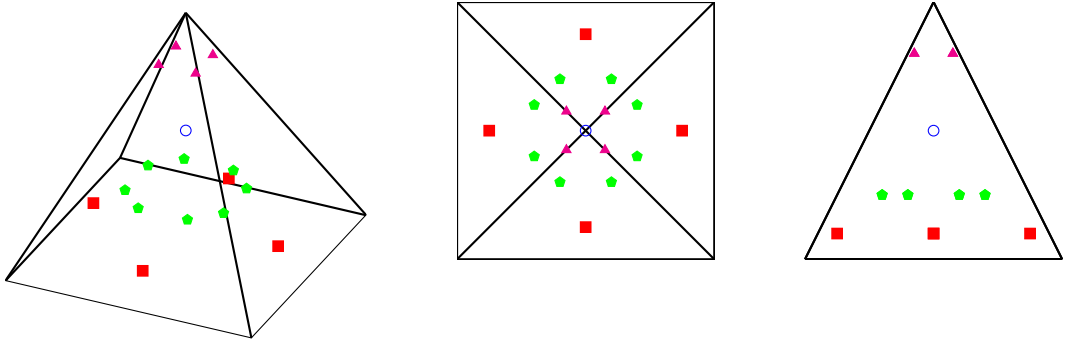
\begin{figure}
    \centering
    \input{_py/pyrmddemo_symorb.tikz}
    \caption{Symmetry orbits for the square pyramid element including isometric (\textit{left}), top (\textit{middle}), and side(\textit{right}) views. Legend: Points generated by $\bar\Pcal_1(0.0)$ (\ref{line:pyrmddemo1:orb1}), $\bar\Pcal_2(0.75, -0.8)$ (\ref{line:pyrmddemo1:orb2}), $\bar\Pcal_3(0.15, 0.5)$ (\ref{line:pyrmddemo1:orb3}), and $\bar\Pcal_4(0.2, 0.4, -0.5)$ (\ref{line:pyrmddemo1:orb4}).}
    \label{fig:pyrmd_symorbit}
\end{figure}

We close this subsection with another concrete example that illustrates the symmetry orbits of the square pyramid. Because the natural coordinates are the normal Cartesian coordinates, the symmetry orbits directly give the desired nodal coordinates. To illustrate this, we select a collection consisting of four symmetry orbits, $\Ical = (1, 2, 3, 4)$ and write the parameter vector $\bar\xibold$ according to (\ref{eqn:parvec}) with
\begin{equation}
    \xibold_1^{(1)} = (0.0), \quad
    \xibold_2^{(2)} = (0.75, -0.8), \quad
    \xibold_3^{(3)} = (0.15, 0.5), \quad
    \xibold_4^{(4)} = (0.2, 0.4, -0.5).
\end{equation}
The parameter $\xibold_3^{(3)}$ generates the nodal coordinates
\begin{gather}
      \Pcal_{3,1}(\xibold_3^{(3)}) = \begin{bmatrix} 0.15 \\ 0.15 \\ 0.5 \end{bmatrix}, 
      \quad \Pcal_{3,2}(\xibold_3^{(3)}) = \begin{bmatrix}  -0.15 \\ 0.15 \\ 0.5 \end{bmatrix}, 
      \quad \Pcal_{3,3}(\xibold_3^{(3)}) = \begin{bmatrix}  0.15 \\ -0.15 \\ 0.5  \end{bmatrix}, 
      \quad \Pcal_{3,4}(\xibold_3^{(3)}) = \begin{bmatrix}  -0.15 \\ -0.15 \\ 0.5 \end{bmatrix}.
\end{gather}
These and all the other nodal coordinates corresponding to $\bar\xibold$ are plotted in Figure~\ref{fig:pyrmd_symorbit} to show the structure of the symmetry orbits. 

%% file: _py/tridemo_symorb.tikz
\begin{tikzpicture}
\begin{axis}[
axis equal image,
axis x line*=bottom,
axis y line*=left,
grid=both,
width=0.45\textwidth,
xmin=-1.05,
xmax=1.05,
ymin=-1.05,
ymax=1.05,
xtick={-1, 0, 1},
ytick={-1, 0, 1},
xlabel={$x$},
ylabel={$y$}]
\addplot [fill=lightgray, opacity=0.6, forget plot]
coordinates {
(-1.00000000e+00, -1.00000000e+00)
( 1.00000000e+00, -1.00000000e+00)
(-1.00000000e+00,  1.00000000e+00)
(-1.00000000e+00, -1.00000000e+00)};

\addplot [solid, black, thick, forget plot]
coordinates {
(-1.00000000e+00, -1.00000000e+00)
( 1.00000000e+00, -1.00000000e+00)
(-1.00000000e+00,  1.00000000e+00)
(-1.00000000e+00, -1.00000000e+00)};

\addplot [blue, only marks, mark=o, mark size=2, mark options={solid}]
coordinates {
(-3.33333333e-01, -3.33333333e-01)};\label{line:tridemo1:orb1}

\addplot [red, only marks, mark=square*, mark size=2, mark options={solid}]
coordinates {
(-6.00000000e-01,  2.00000000e-01)
(-6.00000000e-01, -6.00000000e-01)
( 2.00000000e-01, -6.00000000e-01)};\label{line:tridemo1:orb2}

\addplot [magenta, only marks, mark=triangle*, mark size=2, mark options={solid}]
coordinates {
(-6.00000000e-01, -5.55111512e-17)
(-4.00000000e-01, -5.55111512e-17)
(-6.00000000e-01, -4.00000000e-01)
(-5.55111512e-17, -4.00000000e-01)
(-5.55111512e-17, -6.00000000e-01)
(-4.00000000e-01, -6.00000000e-01)};\label{line:tridemo1:orb3}

\end{axis}
\end{tikzpicture}

%% file: _py/tridemo_nddist.tikz
\begin{tikzpicture}
\begin{axis}[
axis equal image,
axis x line*=bottom,
axis y line*=left,
grid=both,
width=0.45\textwidth,
xmin=-1.05,
xmax=1.05,
ymin=-1.05,
ymax=1.05,
xtick={-1, 0, 1},
ytick={-1, 0, 1},
xlabel={$x$},
ylabel={$y$}]
\addplot [fill=lightgray, opacity=0.6, forget plot]
coordinates {
(-1.00000000e+00, -1.00000000e+00)
( 1.00000000e+00, -1.00000000e+00)
(-1.00000000e+00,  1.00000000e+00)
(-1.00000000e+00, -1.00000000e+00)};

\addplot [solid, black, thick, forget plot]
coordinates {
(-1.00000000e+00, -1.00000000e+00)
( 1.00000000e+00, -1.00000000e+00)
(-1.00000000e+00,  1.00000000e+00)
(-1.00000000e+00, -1.00000000e+00)};

\addplot [blue, only marks, mark=o, mark size=2, mark options={solid}]
coordinates {
(-3.33333333e-01, -3.33333333e-01)};\label{line:tridemo2:orb1}

\addplot [red, only marks, mark=square*, mark size=2, mark options={solid}]
coordinates {
(-5.00000000e-01,  0.00000000e+00)
(-5.00000000e-01, -5.00000000e-01)
( 0.00000000e+00, -5.00000000e-01)};\label{line:tridemo2:orb2a}

\addplot [red, only marks, mark=square, mark size=2, mark options={solid}]
coordinates {
( 0.00000000e+00, -1.00000000e+00)
( 0.00000000e+00,  0.00000000e+00)
(-1.00000000e+00,  0.00000000e+00)};\label{line:tridemo2:orb2b}

\addplot [magenta, only marks, mark=triangle*, mark size=2, mark options={solid}]
coordinates {
( 2.00000000e-01, -4.00000000e-01)
(-8.00000000e-01, -4.00000000e-01)
( 2.00000000e-01, -8.00000000e-01)
(-4.00000000e-01, -8.00000000e-01)
(-4.00000000e-01,  2.00000000e-01)
(-8.00000000e-01,  2.00000000e-01)};\label{line:tridemo2:orb3a}

\addplot [magenta, only marks, mark=triangle, mark size=2, mark options={solid}]
coordinates {
(-1.00000000e+00, -4.00000000e-01)
( 4.00000000e-01, -4.00000000e-01)
(-1.00000000e+00,  4.00000000e-01)
(-4.00000000e-01,  4.00000000e-01)
(-4.00000000e-01, -1.00000000e+00)
( 4.00000000e-01, -1.00000000e+00)};\label{line:tridemo2:orb3b}

\end{axis}
\end{tikzpicture}

%% file: _py/pyrmddemo_symorb.tikz
\begin{tikzpicture}
\begin{groupplot} [
group style={group size = 3 by 1, horizontal sep = 1.2cm}]
\nextgroupplot[axis equal image, xmin=-1, xmax=1, ymin=-1, ymax=1, view={25}{30}, zmin=-1, zmax=1, axis lines=none, width=0.65\textwidth]
\addplot3 [black, thick, forget plot]
coordinates {
(-1.00000000e+00, -1.00000000e+00, -1.00000000e+00)
( 1.00000000e+00, -1.00000000e+00, -1.00000000e+00)
( 1.00000000e+00,  1.00000000e+00, -1.00000000e+00)
(-1.00000000e+00,  1.00000000e+00, -1.00000000e+00)
(-1.00000000e+00, -1.00000000e+00, -1.00000000e+00)};

\addplot3 [black, thick, forget plot]
coordinates {
(-1.00000000e+00, -1.00000000e+00, -1.00000000e+00)
( 0.00000000e+00,  0.00000000e+00,  1.00000000e+00)};

\addplot3 [black, thick, forget plot]
coordinates {
( 1.00000000e+00, -1.00000000e+00, -1.00000000e+00)
( 0.00000000e+00,  0.00000000e+00,  1.00000000e+00)};

\addplot3 [black, thick, forget plot]
coordinates {
( 1.00000000e+00,  1.00000000e+00, -1.00000000e+00)
( 0.00000000e+00,  0.00000000e+00,  1.00000000e+00)};

\addplot3 [black, thick, forget plot]
coordinates {
(-1.00000000e+00,  1.00000000e+00, -1.00000000e+00)
( 0.00000000e+00,  0.00000000e+00,  1.00000000e+00)};

\addplot3 [blue, only marks, mark=o, mark size=2, mark options={solid}]
coordinates {
( 0.00000000e+00,  0.00000000e+00,  0.00000000e+00)};\label{line:pyrmddemo1:orb1}

\addplot3 [red, only marks, mark=square*, mark size=2, mark options={solid}]
coordinates {
( 7.50000000e-01,  0.00000000e+00, -8.00000000e-01)
(-7.50000000e-01,  0.00000000e+00, -8.00000000e-01)
( 0.00000000e+00,  7.50000000e-01, -8.00000000e-01)
( 0.00000000e+00, -7.50000000e-01, -8.00000000e-01)};\label{line:pyrmddemo1:orb2}

\addplot3 [magenta, only marks, mark=triangle*, mark size=2, mark options={solid}]
coordinates {
( 1.50000000e-01,  1.50000000e-01,  6.00000000e-01)
(-1.50000000e-01,  1.50000000e-01,  6.00000000e-01)
( 1.50000000e-01, -1.50000000e-01,  6.00000000e-01)
(-1.50000000e-01, -1.50000000e-01,  6.00000000e-01)};\label{line:pyrmddemo1:orb3}

\addplot3 [green, only marks, mark=pentagon*, mark size=2, mark options={solid}]
coordinates {
( 2.00000000e-01,  4.00000000e-01, -5.00000000e-01)
(-2.00000000e-01,  4.00000000e-01, -5.00000000e-01)
( 2.00000000e-01, -4.00000000e-01, -5.00000000e-01)
(-2.00000000e-01, -4.00000000e-01, -5.00000000e-01)
( 4.00000000e-01,  2.00000000e-01, -5.00000000e-01)
(-4.00000000e-01,  2.00000000e-01, -5.00000000e-01)
( 4.00000000e-01, -2.00000000e-01, -5.00000000e-01)
(-4.00000000e-01, -2.00000000e-01, -5.00000000e-01)};\label{line:pyrmddemo1:orb4}

\nextgroupplot[axis equal image, xmin=-1, xmax=1, ymin=-1, ymax=1, view={0}{90}, zmin=-1, zmax=1, axis lines=none, width=0.35\textwidth]
\addplot3 [black, thick, forget plot]
coordinates {
(-1.00000000e+00, -1.00000000e+00, -1.00000000e+00)
( 1.00000000e+00, -1.00000000e+00, -1.00000000e+00)
( 1.00000000e+00,  1.00000000e+00, -1.00000000e+00)
(-1.00000000e+00,  1.00000000e+00, -1.00000000e+00)
(-1.00000000e+00, -1.00000000e+00, -1.00000000e+00)};

\addplot3 [black, thick, forget plot]
coordinates {
(-1.00000000e+00, -1.00000000e+00, -1.00000000e+00)
( 0.00000000e+00,  0.00000000e+00,  1.00000000e+00)};

\addplot3 [black, thick, forget plot]
coordinates {
( 1.00000000e+00, -1.00000000e+00, -1.00000000e+00)
( 0.00000000e+00,  0.00000000e+00,  1.00000000e+00)};

\addplot3 [black, thick, forget plot]
coordinates {
( 1.00000000e+00,  1.00000000e+00, -1.00000000e+00)
( 0.00000000e+00,  0.00000000e+00,  1.00000000e+00)};

\addplot3 [black, thick, forget plot]
coordinates {
(-1.00000000e+00,  1.00000000e+00, -1.00000000e+00)
( 0.00000000e+00,  0.00000000e+00,  1.00000000e+00)};

\addplot3 [blue, only marks, mark=o, mark size=2, mark options={solid}]
coordinates {
( 0.00000000e+00,  0.00000000e+00,  0.00000000e+00)};\label{line:pyrmddemo1:orb1}

\addplot3 [red, only marks, mark=square*, mark size=2, mark options={solid}]
coordinates {
( 7.50000000e-01,  0.00000000e+00, -8.00000000e-01)
(-7.50000000e-01,  0.00000000e+00, -8.00000000e-01)
( 0.00000000e+00,  7.50000000e-01, -8.00000000e-01)
( 0.00000000e+00, -7.50000000e-01, -8.00000000e-01)};\label{line:pyrmddemo1:orb2}

\addplot3 [magenta, only marks, mark=triangle*, mark size=2, mark options={solid}]
coordinates {
( 1.50000000e-01,  1.50000000e-01,  6.00000000e-01)
(-1.50000000e-01,  1.50000000e-01,  6.00000000e-01)
( 1.50000000e-01, -1.50000000e-01,  6.00000000e-01)
(-1.50000000e-01, -1.50000000e-01,  6.00000000e-01)};\label{line:pyrmddemo1:orb3}

\addplot3 [green, only marks, mark=pentagon*, mark size=2, mark options={solid}]
coordinates {
( 2.00000000e-01,  4.00000000e-01, -5.00000000e-01)
(-2.00000000e-01,  4.00000000e-01, -5.00000000e-01)
( 2.00000000e-01, -4.00000000e-01, -5.00000000e-01)
(-2.00000000e-01, -4.00000000e-01, -5.00000000e-01)
( 4.00000000e-01,  2.00000000e-01, -5.00000000e-01)
(-4.00000000e-01,  2.00000000e-01, -5.00000000e-01)
( 4.00000000e-01, -2.00000000e-01, -5.00000000e-01)
(-4.00000000e-01, -2.00000000e-01, -5.00000000e-01)};\label{line:pyrmddemo1:orb4}

\nextgroupplot[axis equal image, xmin=-1, xmax=1, ymin=-1, ymax=1, view={0}{0}, zmin=-1, zmax=1, axis lines=none, width=0.35\textwidth]
\addplot3 [black, thick, forget plot]
coordinates {
(-1.00000000e+00, -1.00000000e+00, -1.00000000e+00)
( 1.00000000e+00, -1.00000000e+00, -1.00000000e+00)
( 1.00000000e+00,  1.00000000e+00, -1.00000000e+00)
(-1.00000000e+00,  1.00000000e+00, -1.00000000e+00)
(-1.00000000e+00, -1.00000000e+00, -1.00000000e+00)};

\addplot3 [black, thick, forget plot]
coordinates {
(-1.00000000e+00, -1.00000000e+00, -1.00000000e+00)
( 0.00000000e+00,  0.00000000e+00,  1.00000000e+00)};

\addplot3 [black, thick, forget plot]
coordinates {
( 1.00000000e+00, -1.00000000e+00, -1.00000000e+00)
( 0.00000000e+00,  0.00000000e+00,  1.00000000e+00)};

\addplot3 [black, thick, forget plot]
coordinates {
( 1.00000000e+00,  1.00000000e+00, -1.00000000e+00)
( 0.00000000e+00,  0.00000000e+00,  1.00000000e+00)};

\addplot3 [black, thick, forget plot]
coordinates {
(-1.00000000e+00,  1.00000000e+00, -1.00000000e+00)
( 0.00000000e+00,  0.00000000e+00,  1.00000000e+00)};

\addplot3 [blue, only marks, mark=o, mark size=2, mark options={solid}]
coordinates {
( 0.00000000e+00,  0.00000000e+00,  0.00000000e+00)};\label{line:pyrmddemo1:orb1}

\addplot3 [red, only marks, mark=square*, mark size=2, mark options={solid}]
coordinates {
( 7.50000000e-01,  0.00000000e+00, -8.00000000e-01)
(-7.50000000e-01,  0.00000000e+00, -8.00000000e-01)
( 0.00000000e+00,  7.50000000e-01, -8.00000000e-01)
( 0.00000000e+00, -7.50000000e-01, -8.00000000e-01)};\label{line:pyrmddemo1:orb2}

\addplot3 [magenta, only marks, mark=triangle*, mark size=2, mark options={solid}]
coordinates {
( 1.50000000e-01,  1.50000000e-01,  6.00000000e-01)
(-1.50000000e-01,  1.50000000e-01,  6.00000000e-01)
( 1.50000000e-01, -1.50000000e-01,  6.00000000e-01)
(-1.50000000e-01, -1.50000000e-01,  6.00000000e-01)};\label{line:pyrmddemo1:orb3}

\addplot3 [green, only marks, mark=pentagon*, mark size=2, mark options={solid}]
coordinates {
( 2.00000000e-01,  4.00000000e-01, -5.00000000e-01)
(-2.00000000e-01,  4.00000000e-01, -5.00000000e-01)
( 2.00000000e-01, -4.00000000e-01, -5.00000000e-01)
(-2.00000000e-01, -4.00000000e-01, -5.00000000e-01)
( 4.00000000e-01,  2.00000000e-01, -5.00000000e-01)
(-4.00000000e-01,  2.00000000e-01, -5.00000000e-01)
( 4.00000000e-01, -2.00000000e-01, -5.00000000e-01)
(-4.00000000e-01, -2.00000000e-01, -5.00000000e-01)};\label{line:pyrmddemo1:orb4}

\end{groupplot}\end{tikzpicture}

%% file: _optim.tex
From the symmetric nodal distributions introduced in Section~\ref{sec:absform}, we introduce the proposed optimization-based nodal distributions (Section~\ref{sec:optim:form}) and detail procedures to define admissible symmetry orbit collections (Section~\ref{sec:optim:admis}) and constraints (Section~\ref{sec:optim:compat}) for $H^1$-conforming finite elements.

\subsection{Optimization formulation}
\label{sec:optim:form}
Given a collection of symmetry orbits $\Ical$ and corresponding constraints (Section~\ref{sec:absform:symorbcoll}), we define an optimization-based symmetric nodal distribution as $\Nbb_{\bar\xibold_\star} \subset \Omega$, where $\bar{\xibold}_\star$ is the solution of the following optimization problem
\begin{equation} \label{eqn:optnode}
    \optunc{\bar\xibold\in\bar\Cbb_\xi}{f(\Nbb_{\bar\xibold})}
\end{equation}
and $f : \bar\Cbb_\xi \rightarrow \Rbb$ is the objective function. Ideally, the objective function would be the Lebesgue constant, $\Lambda : \bar\Cbb_\xi \rightarrow \Rbb$, defined as
\begin{equation}
\Lambda : \bar\xibold \mapsto \max_{x \in \Omega} \sum_{i=1}^n |\ell_i(x; \Nbb_{\bar\xibold})|,
\end{equation}
where $\ell_i(x; \Nbb)$ is the Lagrange basis function associated with the $i$th node in the nodal distribution $\Nbb$. However, this non-smooth function is not suitable for gradient-based optimization so we use a smooth surrogate for the Lebesgue constant as the objective function
\begin{equation} \label{eqn:lebobj}
f : \bar\xibold \mapsto \sum_{i=1}^n \int_{x\in\Omega} |\ell_i(x; \Nbb_{\bar\xibold})|^2 \, dx,
\end{equation}
The integral is discretized using a quadrature rule that exactly integrates a polynomial of degree $2p$, where $p$ is the polynomial degree of $\ell_i$, to exactly evaluate $f$.
The optimization problem in (\ref{eqn:optnode}) is initialized from the symmetry parameters corresponding to a known, symmetric nodal distribution (GLL nodes for tensor product elements, \cite{Isaac} for simplices, uniform nodes for other elements). 

\subsection{Admissible symmetry orbit collections}
\label{sec:optim:admis}
The optimization problem (\ref{eqn:optnode}) is well-defined once the symmetry orbit collection and constraints are defined. Given a complete polynomial space $\Pbb$ associated with finite element with geometry $\Omega$, we define the admissible symmetry orbit collections
to be any $\Ical = (k_1,\dots,k_n) \in \{1,\dots,M\}^n$ with cardinality equal to the dimension of $\Pbb$, i.e., the solution of
\begin{equation}
    \sum_{j=1}^n m_{k_j} = \dim\Pbb.
\end{equation}
This ensures the Lagrangian functions associated with the nodes $\Nbb_{\bar\xibold_\star}$ is a basis of $\Pbb$ provided the nodal coordinates are unique.

\begin{remark}
    This procedure generates many different admissible symmetry orbit collections. Some of these collections cannot result in optimal nodes, as they require too many nodes in the same symmetry orbits. This leads to nodes being concentrated along diagonal, vertical, or horizontal lines within the elements. To avoid such collections, the Lebesgue constant of the initial distribution is measured. For the non-unisolvent collections, the initial Lebesgue constant is extremely high (on the order of $10^{20}$) allowing for non-viable symmetry orbit collections to be quickly rejected. For simplex elements, \cite{unisolv} has shown that there is only one symmetry orbit collection which satisfies the unisolvency conditions. However, in general, unisolvency conditions are difficult to detect directly and the proposed brute-force approach has proven effective. 
\end{remark}


\subsection{Cross-element compatibility}
\label{sec:optim:compat}
Finally, we describe a procedure to construct constrained symmetry orbits of $\Omega\subset\Rbb^d$ that ensure the resulting nodal distribution on faces of the element $\partial\Omega_f \subset \Rbb^d$, $f=1,\dots, n_f$ ($n_f$ is the number of planar faces of the element $\Omega$) agrees with some prescribed distribution, $\Nbb_{d-1} \subset \Gamma$, of the face geometry, $\Gamma \subset \Rbb^{d-1}$. This is a desirable property because it ensures nodes of adjacent elements will exactly coincide when used in a conforming, mixed finite element mesh (multiple element geometries), which is practically necessary for $H^1$-conforming finite element discretizations. Naively generating optimal distributions for element geometries would not guarantee this property.

To enforce this constraint, we must identify which symmetry orbits of the collection $\Ical = (k_1,\dots,k_n)$ require constraints of the form (\ref{eqn:symcoord_addbnd}) and what constraint to impose, i.e., values of $(\Abm_{k_j}^{(j)}, \bbm_{\mathrm{l},k_j}^{(j)}, \bbm_{\mathrm{u},k_j}^{(j)})$ to use. Consider any $\rbm \in \Nbb_{d-1} \subset \Gamma$. First, we map this to any face of the element $\Omega$ and use $\Rcal^{-1}$ to obtain the corresponding natural coordinate, denoted $\hat\lambdabold$. By symmetry, the condition will be applied on all faces regardless of which face is chosen. Then, we impose the constraint on the symmetry parameter in (\ref{eqn:symcoord_addbnd}) with
\begin{equation} \label{eqn:cross1}
    \Abm_{k_j}^{(j)} = \Sbm_{k_j,1}, \qquad
    \bbm_{\mathrm{l},k_j}^{(j)} = \bbm_{\mathrm{u},k_j}^{(j)} = \hat\lambdabold - \sigmabold_{k_j,1}
\end{equation}
because this is equivalent to constraints on the natural coordinate of the form (\ref{eqn:nat_bnd}) with $\Bbm_\lambda = \Ibm$ and $\vbm_\mathrm{l} = \vbm_\mathrm{u} = \hat\lambdabold$.
This constraint is applied to $j$th entry of the symmetry orbit collection, corresponding to symmetry orbit $k_j$, where $j$ is the smallest number such that
\begin{equation} \label{eqn:cross2}
    \Pcal_{k_j,i}(\xibold) = \hat\lambdabold
\end{equation}
has a solution $\xibold$ for some $i \in \{1,\dots,m_{k_j}\}$ and the entry $j$ has not already been assigned a constraint. All points in $\Nbb_{d-1}$ that, after mapping to a face of $\Omega$, coincide with points in the $\bar{\Pcal}_{k_j}(\xibold)$ are removed from $\Nbb_{d-1}$ because the above construction ensures those points will be fixed at their appropriate locations in $\Omega$. This process is repeated until all points in $\Nbb_{d-1}$ are exhausted. Any $j$ that has not been assigned a constraint will use the standard symmetry orbit with $\Abm_{k_j}^{(j)} = \bbm_{\mathrm{l},k_j}^{(j)} = \bbm_{\mathrm{u},k_j}^{(j)} = \emptyset.$
The resulting constraints $\{(\Abm_{k_j}^{(j)}, \bbm_{\mathrm{l},k_j}^{(j)}, \bbm_{\mathrm{u},k_j}^{(j)})\}_{j=1}^n$ ensure the nodal distribution on each face of $\Omega$ matches the nodal distribution $\Nbb_{d-1}$ of $\Gamma$, which in turn ensures nodes on adjacent elements in a conforming mesh (mixed or homogeneous element geometries) will coincide. These constraints also ensure that all edges and faces have the same nodal distributions, making it possible to line up elements of different dimensions if desired.

%% file: _numexp.tex
In this section, we generate the proposed optimization-based nodal distributions for the seven elements considered in this work. We use (\ref{eqn:lebobj}), to be referred to as the Lebesgue objective, as the optimization objective because it is a smooth surrogate of the Lebesgue constant. Cross-element compatibility described in Section~\ref{sec:optim:compat} is enforced. We use the Julia interface to IPOPT \cite{Wächter2006} with default settings to solve the optimization problem in (\ref{eqn:optnode}). Convergence is obtained when the first-order optimality conditions are satisfied to a tolerance of $10^{-10}$ or $50$ major iterations is exceeded, whichever comes first. The optimizer is initialized from known nodal distributions (Section~\ref{sec:optim:form}).
For each element, we compare the Lebesgue constant and condition number of the finite element mass matrix associated with the proposed nodal distribution to those of standard distributions found in the literature. 

\subsection{Line}
\label{sec:numexp:line}
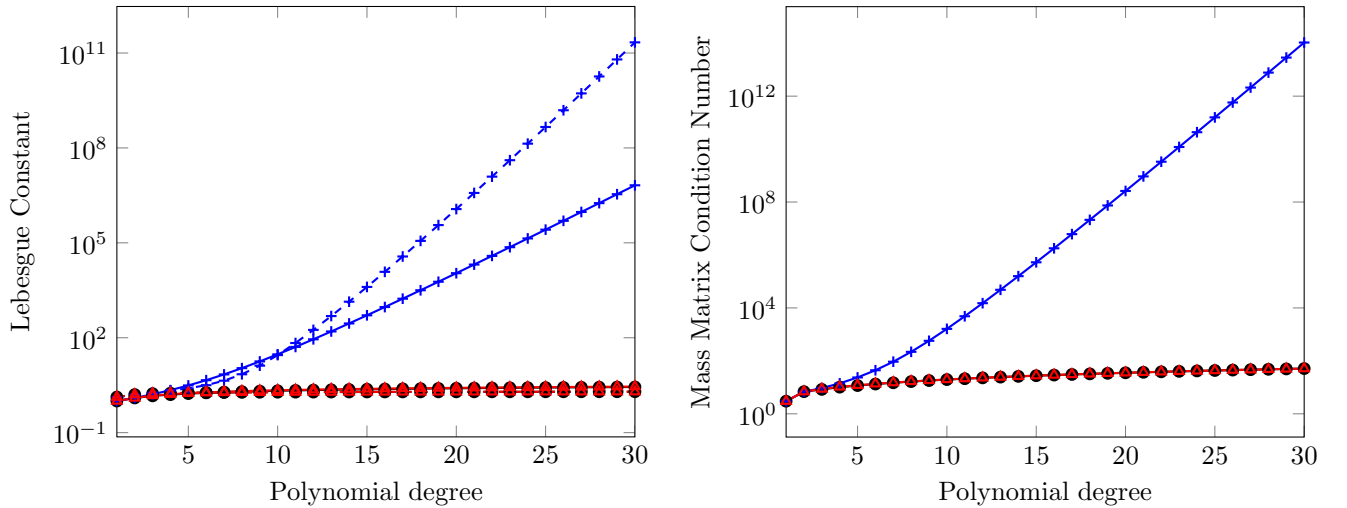
\begin{figure}[H]
    \centering
    \input{_py/compare_line.tikz}
    \caption{The Lebesgue constant (\textit{left}) and mass matrix condition number (\textit{right}) of the optimized nodes (\ref{line:lebcnst:0}), GLL nodes (\ref{line:lebcnst:2}), and uniform nodes (\ref{line:lebcnst:1}) for the \textit{line} element. The Lebesgue objective (\ref{eqn:lebobj}) for each nodal distribution is also included in the \textit{left} figure for each nodal distribution (\textit{dashed}).}
    \label{fig:compare_line}
\end{figure}

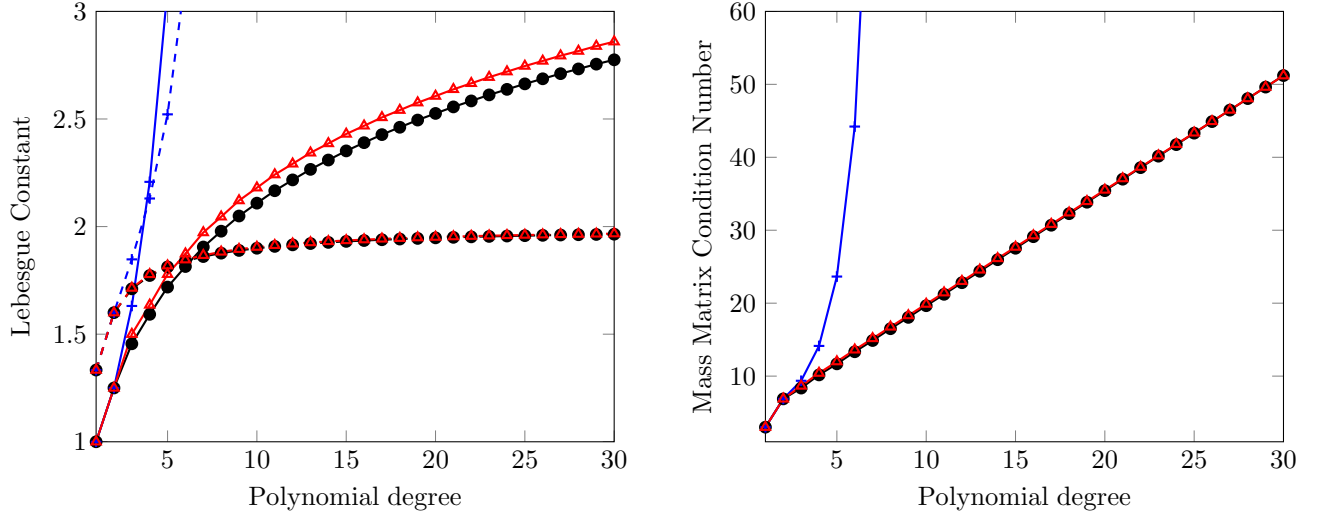
\begin{figure}[H]
    \centering
    \input{_py/compare_line_zoom.tikz}
    \caption{Zoom of Figure~\ref{fig:compare_line} to emphasize difference between optimized and GLL nodes.}
    \label{fig:compare_line_zoom}
\end{figure}

For the only 1D element, the line, the optimized distribution is compared with the GLL and uniform distributions up to polynomial degree 30 by means of the Lebesgue constant, Lebesgue objective, and mass matrix condition number (Figure \ref{fig:compare_line}-\ref{fig:compare_line_zoom}). The Lebesgue objective function is a good indicator of the overall trend in the Lebesgue constant. The optimized nodes achieve a lower Lebesgue objective than the GLL nodes, which is expected because this is the quantity being minimized. The optimized nodes and GLL nodes achieve similar Lebesgue constant and mass matrix condition number (slightly lower for the optimized nodes). This shows the optimized nodes are able to slightly improve on the near-optimal GLL distribution. Both the optimized and GLL nodes substantially outperform the uniform distribution, as expected. The GLL and uniform nodes satisfy the same constraints as the optimized nodes. All three nodal distributions are closed (include the vertices of the line) and obey the symmetries of the geometry.

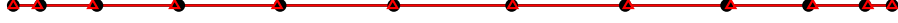
\begin{figure}[H]
    \centering
    \input{_py/nodes_line.tikz}
    \caption{Optimized nodal distribution (\ref{line:nodes:0}) and GLL nodes (\ref{line:nodes:2}) for line element with polynomial completeness $p = 11$.}
    \label{fig:nodes_line}
\end{figure}

The information in Figures~\ref{fig:compare_line}-\ref{fig:compare_line_zoom} is tabulated in Table~\ref{tab:line} and
Figure \ref{fig:nodes_line} shows a sample distribution of optimized nodes for $p=11$.

\subsection{Quadrilateral}
The optimized distribution is compared with the GLL and uniform distributions by means of the Lebesgue constant, Lebesgue objective, and mass matrix condition number up to polynomial degree 23 (Figure \ref{fig:compare_quad}-\ref{fig:compare_quad_zoom}). Similar to the line element, the Lebesgue objective and constant show similar trends with the optimized nodes achieving slightly lower values for all three metrics. This shows the optimized nodes are able to slightly improve on the near-optimal GLL distribution. The GLL and uniform nodes satisfy similar constraints as the optimized nodes. All three distributions obey the symmetries of the geometry and reduce to the corresponding 1D distribution on the edges of the quadrilateral.

\begin{figure}[H]
    \centering
    \input{_py/nodes_quad.tikz}
    \caption{Optimized nodal distribution (\ref{line:nodes:0}) and GLL nodes (\ref{line:nodes:2}) for quadrilateral element with polynomial completeness $p = 7$.}
    \label{fig:nodes_quad}
\end{figure}
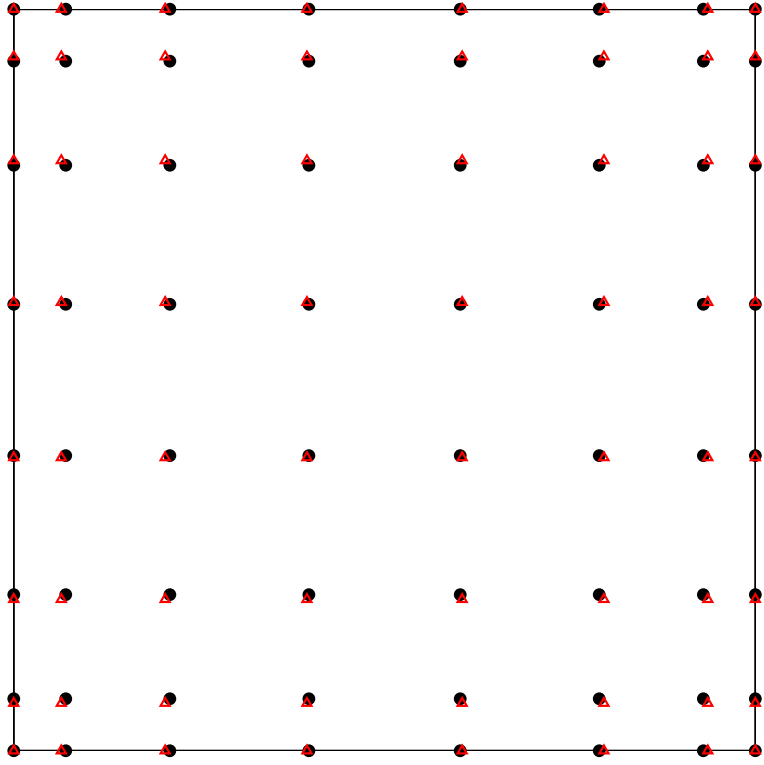

The information in Figures~\ref{fig:compare_quad}-\ref{fig:compare_quad_zoom} is tabulated in Table~\ref{tab:quad} and
Figure \ref{fig:nodes_quad} shows a sample distribution of optimized nodes for $p=7$.

\begin{figure}[H]
    \centering
    \input{_py/compare_quad.tikz}
    \caption{The Lebesgue constant (\textit{left}) and mass matrix condition number (\textit{right}) of the optimized nodes (\ref{line:lebcnst:0}), GLL nodes (\ref{line:lebcnst:2}), and uniform nodes (\ref{line:lebcnst:1}) for the \textit{quadrilateral} element. The Lebesgue objective (\ref{eqn:lebobj}) for each nodal distribution is also included in the \textit{left} figure for each nodal distribution (\textit{dashed}).}
    \label{fig:compare_quad}
\end{figure}
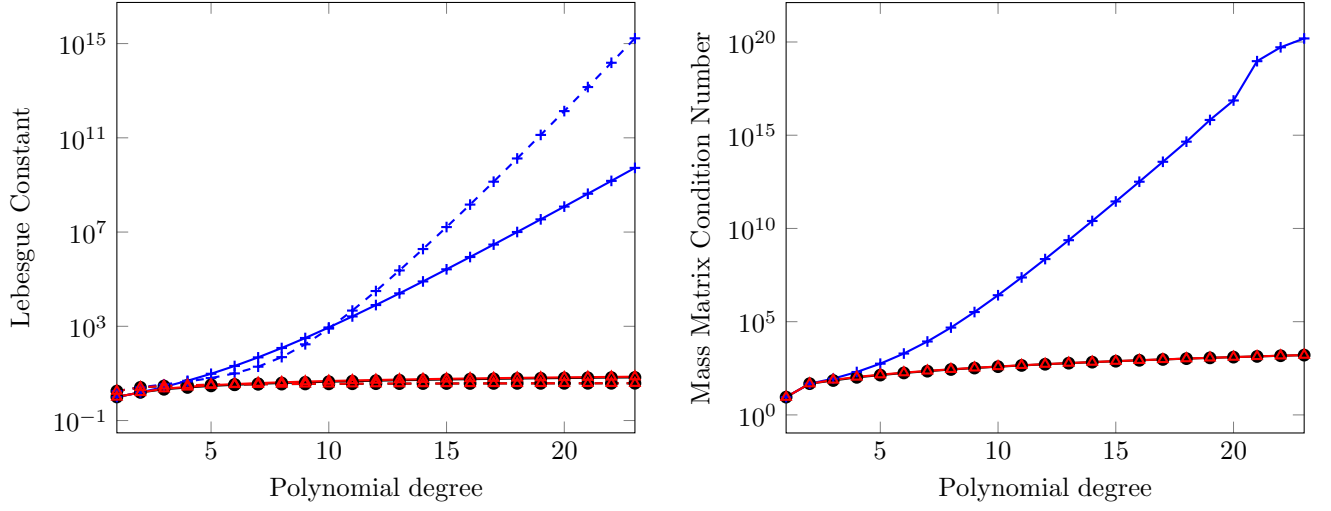


\begin{figure}[H]
    \centering
    \input{_py/compare_quad_zoom.tikz}
    \caption{Zoom of Figure~\ref{fig:compare_quad} to emphasize difference between optimized and GLL nodes.}
    \label{fig:compare_quad_zoom}
\end{figure}
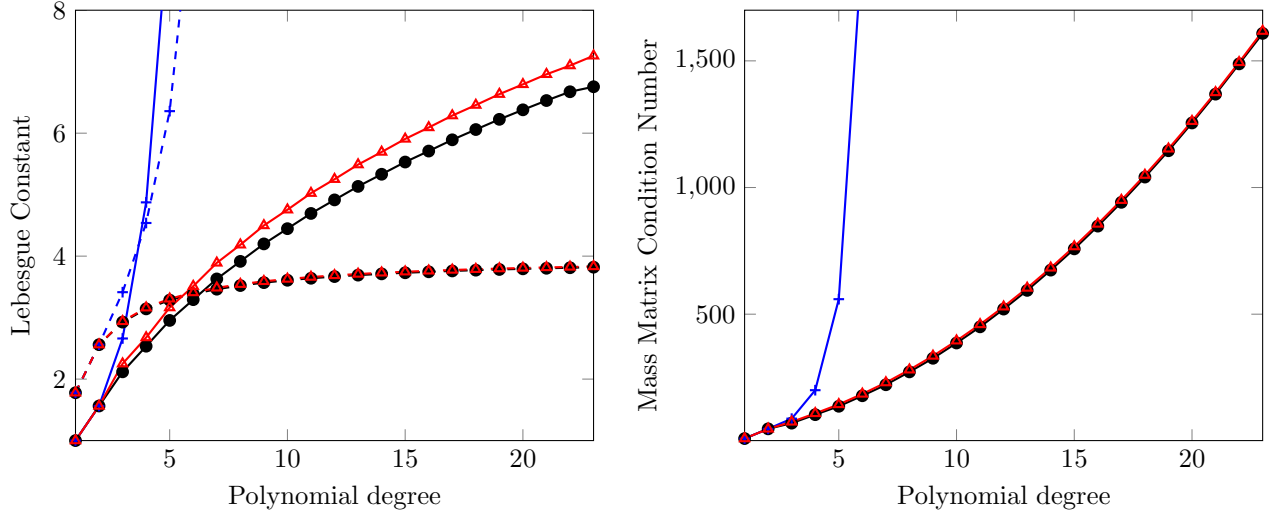


\subsection{Triangle}
\label{sec:numexp:tri}
\begin{figure}[H]
    \centering
    \input{_py/compare_tri.tikz}
    \caption{The Lebesgue constant (\textit{left}) and mass matrix condition number (\textit{right}) of the optimized nodes (\ref{line:lebcnst:0}), Isaac nodes (\ref{line:lebcnst:2}), Warburton nodes (\ref{line:lebcnst:3}), and uniform nodes (\ref{line:lebcnst:1}) for the \textit{triangle} element. The Lebesgue objective (\ref{eqn:lebobj}) for each nodal distribution is also included in the \textit{left} figure for each nodal distribution (\textit{dashed}).}
    \label{fig:compare_tri}
\end{figure}
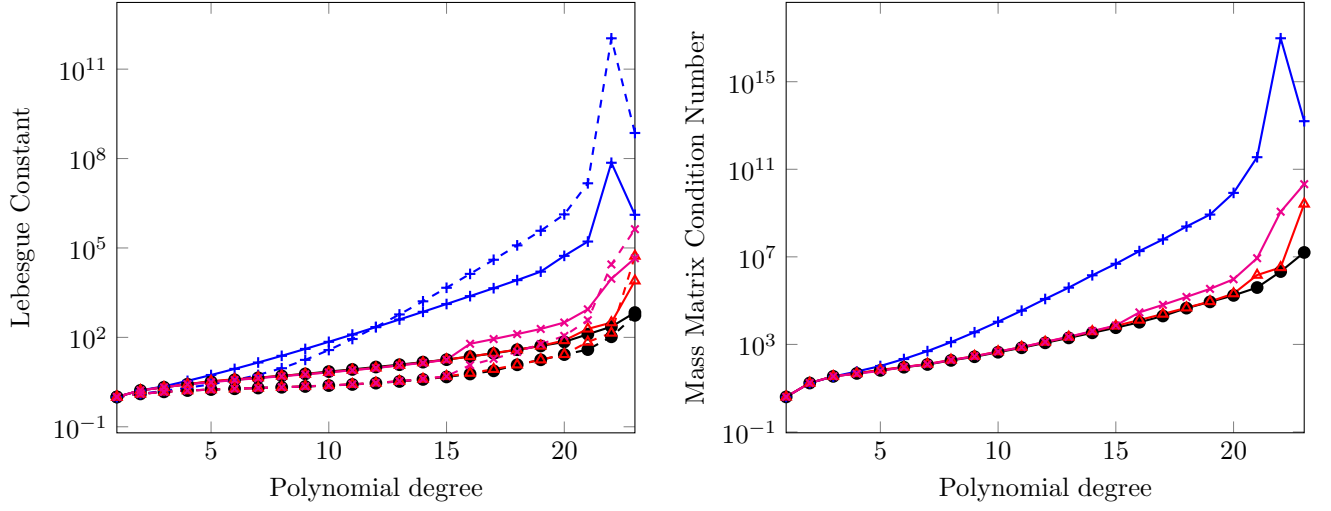

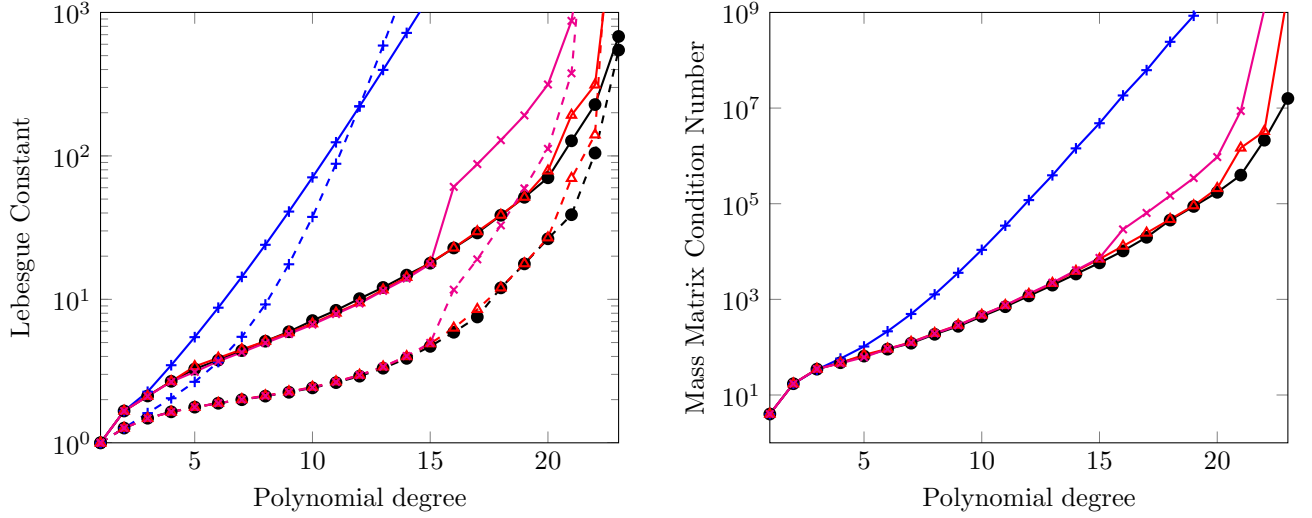
\begin{figure}[H]
    \centering
    \input{_py/compare_tri_zoom.tikz}
    \caption{Zoom of Figure~\ref{fig:compare_tri} to emphasize difference between optimized and Isaac nodes.}
    \label{fig:compare_tri_zoom}
\end{figure}

The optimized distribution is compared with commonly used nodes from the literature up to polynomial degree 23. In particular, the optimized distribution is compared with the uniform nodes, the Isaac nodes found in \cite{Isaac}, and the Warburton nodes found in \cite{warburton2006explicit} by means of the Lebesgue constant, Lebesgue objective, and mass matrix condition number (Figure \ref{fig:compare_tri}-\ref{fig:compare_tri_zoom}). These metrics are nearly identical for the Isaac nodes and optimized distributions until polynomial degree 15. Above degree 15, the optimized nodes achieve lower values for each of these metrics. The other distributions have larger Lebesgue constant, Lebesgue objective, and mass matrix condition number. The uniform nodes, Warburton nodes, optimized nodes, and Isaac nodes all satisfy the cross-element compatibility and obey the symmetries of the geometry. Each reduces to the corresponding 1D distribution on the edges of the triangle.

\begin{figure}[H]
    \centering
    \input{_py/nodes_tri.tikz}
    \caption{Optimized nodal distribution (\ref{line:nodes:0}) and Isaac nodes (\ref{line:nodes:2}) for triangle element with polynomial completeness $p = 7$.}
    \label{fig:nodes_tri}
\end{figure}
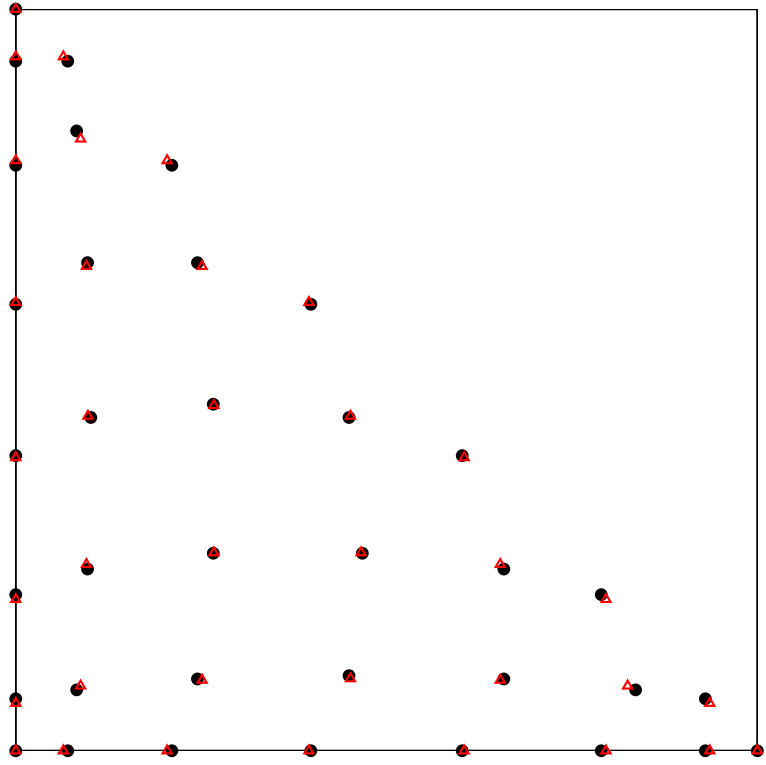

The information in Figures~\ref{fig:compare_tri}-\ref{fig:compare_tri_zoom} is tabulated in Table~\ref{tab:tri} and
Figure \ref{fig:nodes_tri} shows a sample distribution of optimized nodes for $p=7$.

\subsection{Hexahedron}
The optimized distribution is compared with the GLL and uniform distributions by means of the Lebesgue constant, Lebesgue objective, and mass matrix condition number up to polynomial degree 9 (Figure \ref{fig:compare_hex}-\ref{fig:compare_hex_zoom}).
Similar to the line element, the Lebesgue objective and constant show similar trends with the optimized
nodes achieving slightly lower values for all three metrics. This shows the optimized nodes are able to slightly improve on the near-optimal GLL distribution. The GLL and uniform nodes satisfy similar constraints as the optimized nodes. All three distributions obey the symmetries of the geometry and reduce to the corresponding 1D and 2D distributions on the edges and faces, respectively, of the hexahedron.

\begin{figure}[H]
    \centering
    \input{_py/nodes_hex.tikz}
    \caption{Optimized nodal distribution for hexahedron element with polynomial completeness $p = 4$.}
    \label{fig:nodes_hex}
\end{figure}
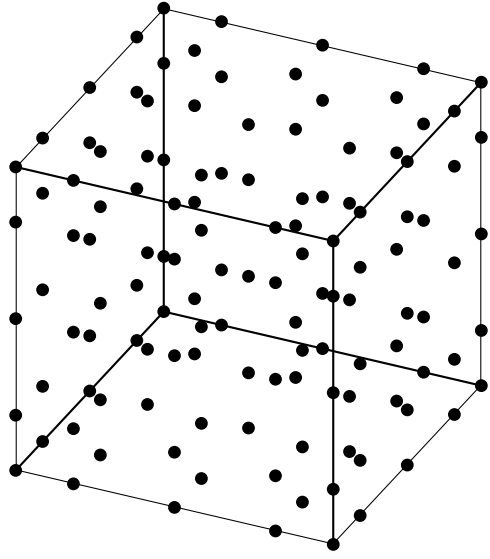

The information in Figures~\ref{fig:compare_hex}-\ref{fig:compare_hex_zoom} is tabulated in Table~\ref{tab:hex} and
Figure \ref{fig:nodes_hex} shows a sample distribution of optimized nodes for $p=4$.

\begin{figure}[H]
    \centering
    \input{_py/compare_hex.tikz}
    \caption{The Lebesgue constant (\textit{left}) and mass matrix condition number (\textit{right}) of the optimized nodes (\ref{line:lebcnst:0}), GLL nodes (\ref{line:lebcnst:2}), and uniform nodes (\ref{line:lebcnst:1}) for the \textit{hexahedron} element. The Lebesgue objective (\ref{eqn:lebobj}) for each nodal distribution is also included in the \textit{left} figure for each nodal distribution (\textit{dashed}).}
    \label{fig:compare_hex}
\end{figure}
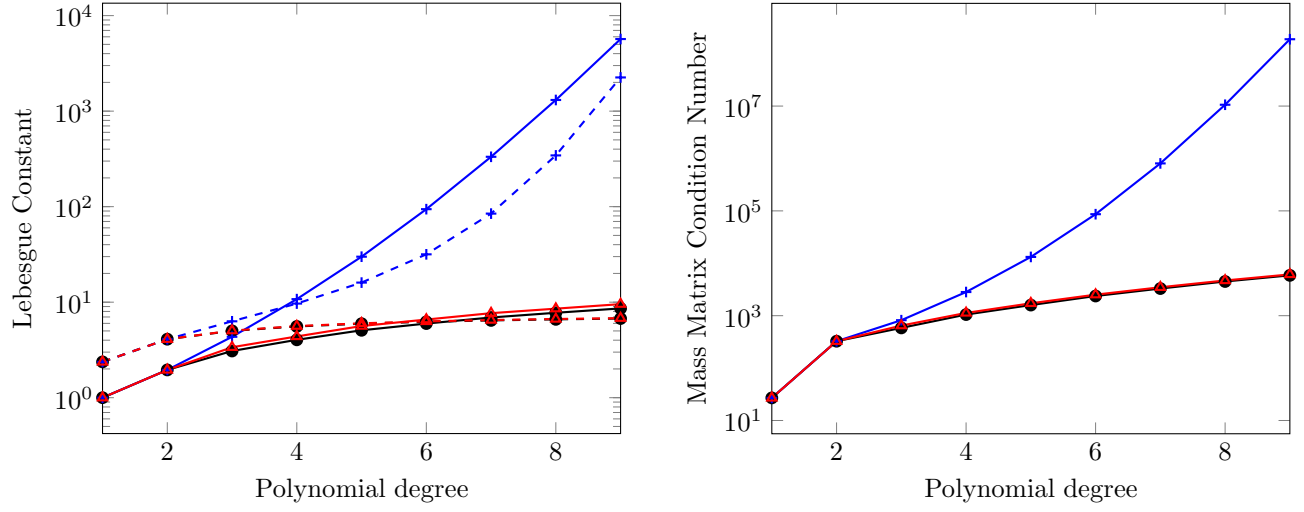


\begin{figure}[H]
    \centering
    \input{_py/compare_hex_zoom.tikz}
    \caption{Zoom of Figure~\ref{fig:compare_hex} to emphasize difference between optimized and GLL nodes.}
    \label{fig:compare_hex_zoom}
\end{figure}
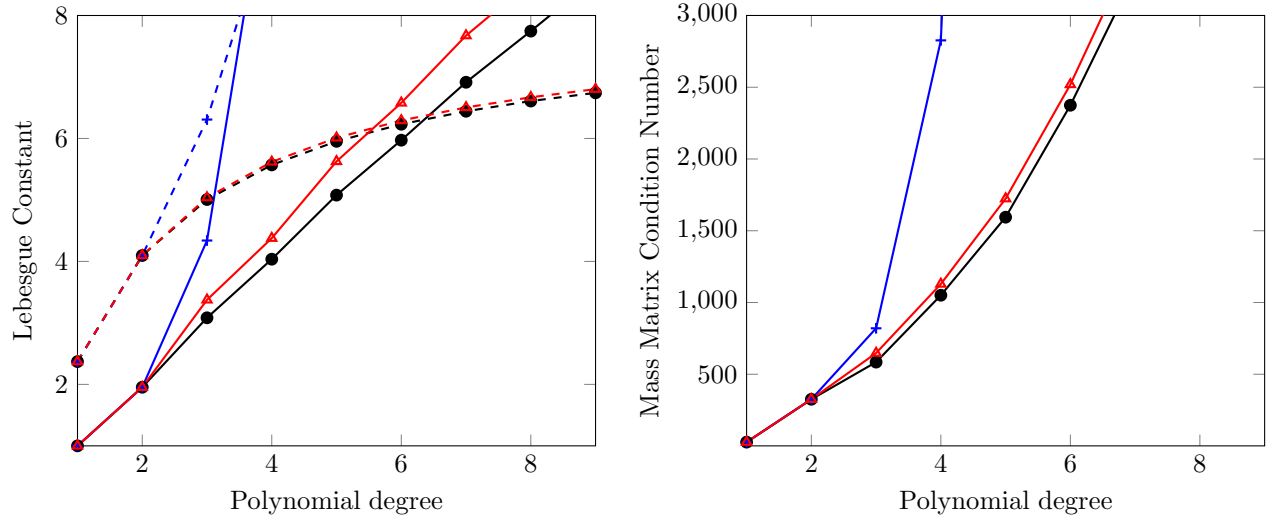


\subsection{Tetrahedron}
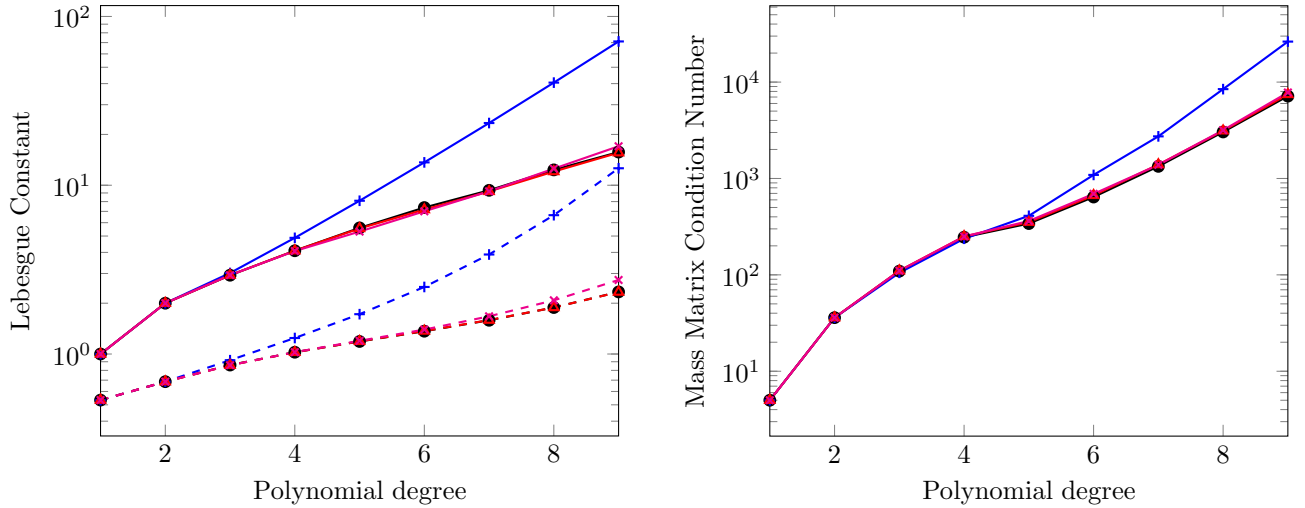
\begin{figure}[H]
    \centering
    \input{_py/compare_tet.tikz}
    \caption{The Lebesgue constant (\textit{left}) and mass matrix condition number (\textit{right}) of the optimized nodes (\ref{line:lebcnst:0}), Isaac nodes (\ref{line:lebcnst:2}), Warburton nodes (\ref{line:lebcnst:3}), and uniform nodes (\ref{line:lebcnst:1}) for the \textit{tetrahedron} element. The Lebesgue objective (\ref{eqn:lebobj}) for each nodal distribution is also included in the \textit{left} figure for each nodal distribution (\textit{dashed}).}
    \label{fig:compare_tet}
\end{figure}

The optimized distribution is compared with several distributions from the literature, namely the Issac \cite{Isaac} and the Warburton \cite{warburton2006explicit} nodes, as well as the uniform distributions by means of the Lebesgue constant, Lebesgue objective, and mass matrix condition number up to polynomial degree 9 (Figures \ref{fig:compare_tet}-\ref{fig:compare_tet_zoom}). The optimized distribution and Isaac nodes achieve the lowest metrics across all polynomial degrees, although the Warburton nodes are not much larger for the degrees considered. The Warburton, Isaac, and uniform nodes satisfy similar constraints as the optimized nodes. All four distributions obey the symmetries of the geometry and reduce to the corresponding 1D and 2D distributions on the edges and faces, respectively, of the tetrahedron.

\begin{figure}[H]
    \centering
    \input{_py/nodes_tet.tikz}
    \caption{Optimized nodal distribution for tetrahedron element with polynomial completeness $p = 4$.}
    \label{fig:nodes_tet}
\end{figure}
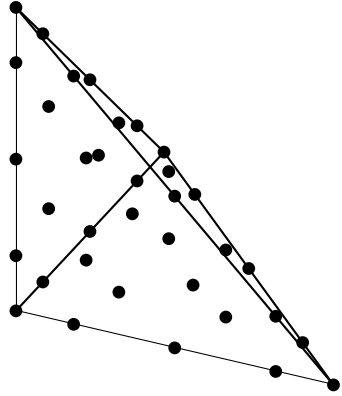

The information in Figures~\ref{fig:compare_tet}-\ref{fig:compare_tet_zoom} is tabulated in Table~\ref{tab:tet} and
Figure \ref{fig:nodes_tet} shows a sample distribution of optimized nodes for $p=4$.

\begin{figure}[H]
    \centering
    \input{_py/compare_tet_zoom.tikz}
    \caption{Zoom of Figure~\ref{fig:compare_tet} to emphasize difference between optimized and Isaac nodes.}
    \label{fig:compare_tet_zoom}
\end{figure}
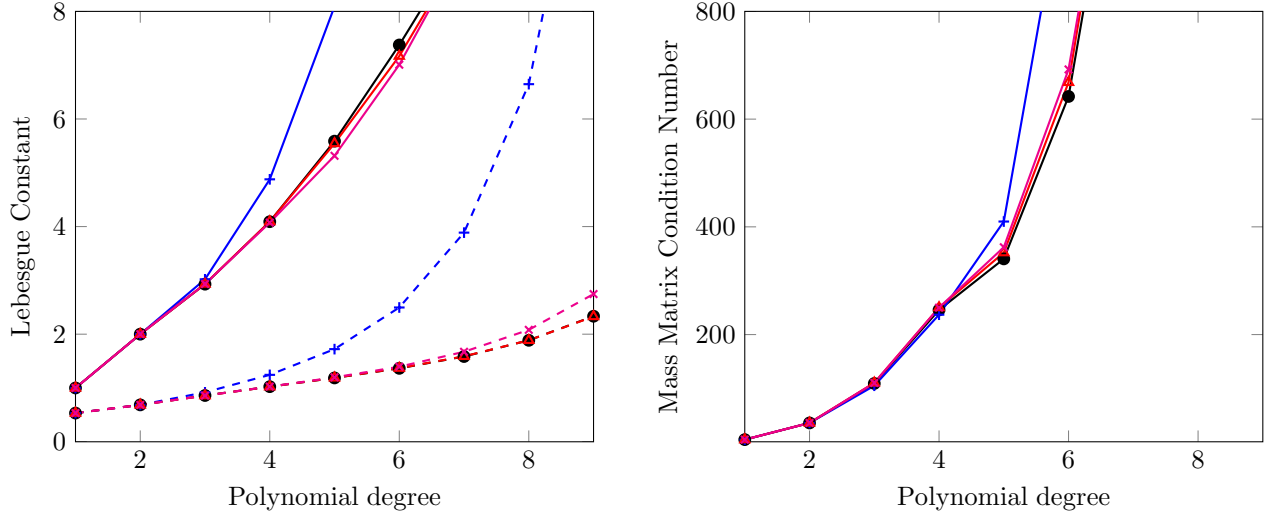


\subsection{Triangular Prism}
The optimized distribution is compared with the Cartesian product of popular distributions for the triangle with the GLL nodes up to polynomial degree 9 by means of the Lebesgue constant, Lebesgue objective, and mass matrix condition number (Figures \ref{fig:compare_triprism}-\ref{fig:compare_triprism_zoom}). Namely the Issac \cite{Isaac} and the Warburton \cite{warburton2006explicit} distributions are considered, and denoted as Isaac $\times$ GLL and Warburton $\times$ GLL. We also consider a tensor product of the optimized triangular nodes (Section~\ref{sec:numexp:tri}) with the optimized line nodes (Section~\ref{sec:numexp:line}).  All distributions perform similarly with the Warburton $\times$ GLL achieving the lowest Lebesgue constant and the optimized nodes reaching the lowest Lebesgue objective and mass matrix condition number. The tensor product of the optimized triangle and optimized line nodes leads to identical nodes as the direct optimization approach suggesting the optimal distribution maintains the tensor product structure. All nodal distributions considered satisfy similar constraints as the optimized nodes. All five distributions obey the symmetries of the geometry and reduce to the corresponding 1D and 2D distributions on the edges and faces, respectively, of the triangular prism.

\begin{figure}[H]
    \centering
    \input{_py/nodes_triprism.tikz}
    \caption{Optimized nodal distribution for triangular prism element with polynomial completeness $p = 4$.}
    \label{fig:nodes_triprism}
\end{figure}
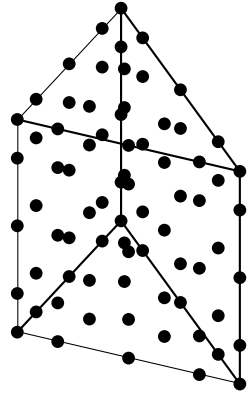

The information in Figures~\ref{fig:compare_triprism}-\ref{fig:compare_triprism_zoom} is tabulated in Table~\ref{tab:triprism} and
Figure \ref{fig:nodes_triprism} shows a sample distribution of optimized nodes for $p=4$.

\begin{figure}[H]
    \centering
    \input{_py/compare_triprism.tikz}
    \caption{The Lebesgue constant (\textit{left}) and mass matrix condition number (\textit{right}) of the optimized nodes (\ref{line:lebcnst:0}), Isaac $\times$ GLL nodes (\ref{line:lebcnst:2}), Warburton $\times$ GLL nodes (\ref{line:lebcnst:3}),
    optimize triangle $\times$ optimized line nodes (\ref{line:lebcnst:4}), and uniform nodes (\ref{line:lebcnst:1}) for the \textit{triangular prism} element. The Lebesgue objective (\ref{eqn:lebobj}) for each nodal distribution is also included in the \textit{left} figure for each nodal distribution (\textit{dashed}).}
    \label{fig:compare_triprism}
\end{figure}
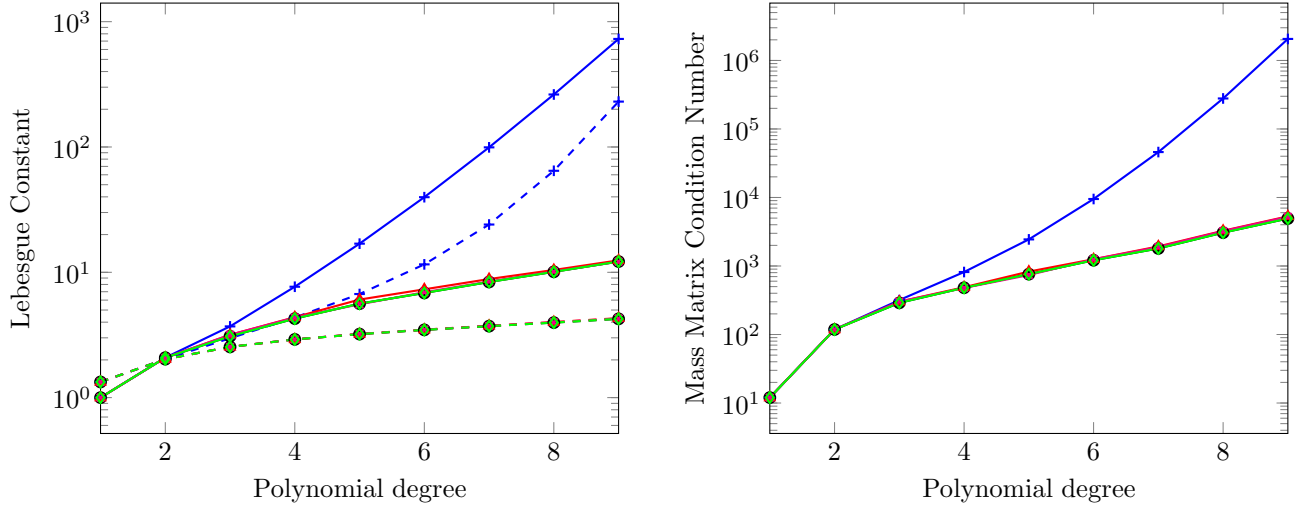


\begin{figure}[H]
    \centering
    \input{_py/compare_triprism_zoom.tikz}
    \caption{Zoom of Figure~\ref{fig:compare_triprism} to emphasize difference between optimized and Isaac $\times$ GLL nodes.}
    \label{fig:compare_triprism_zoom}
\end{figure}
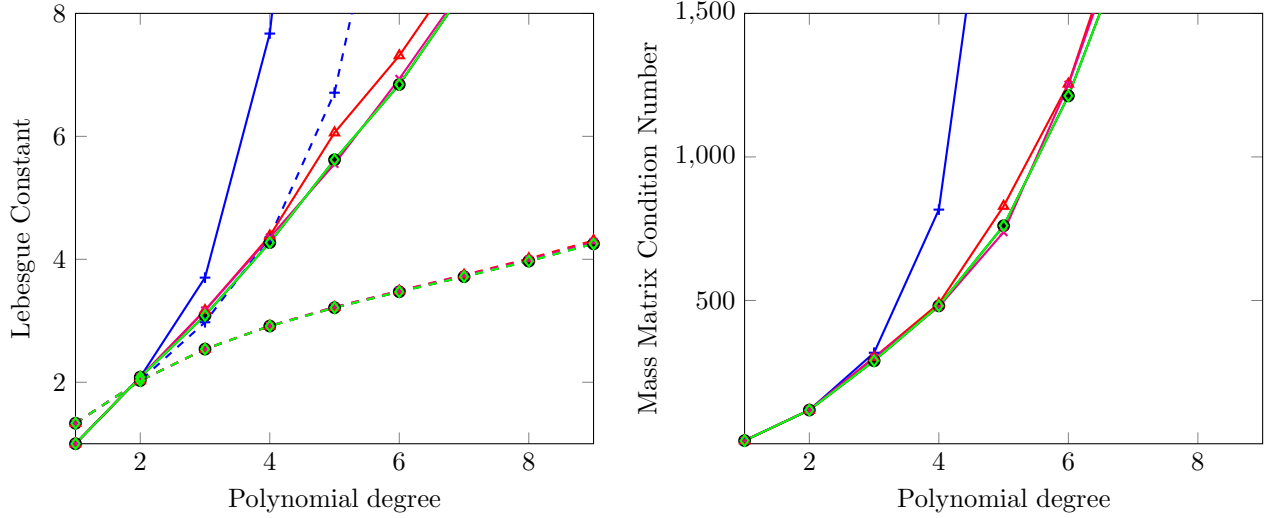


\subsection{Square Pyramid}
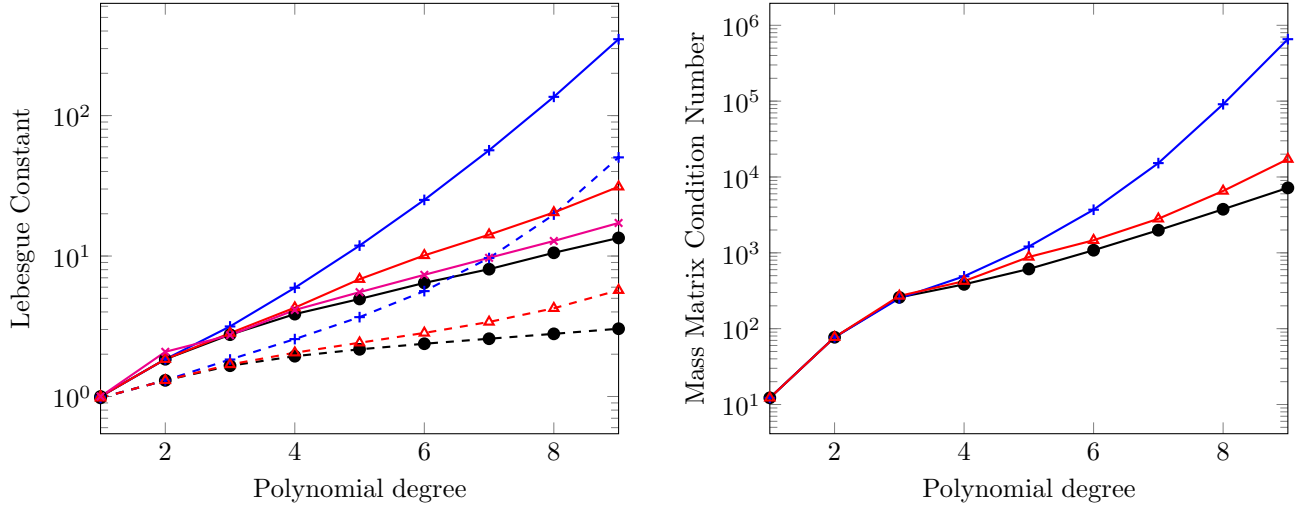
\begin{figure}[H]
    \centering
    \input{_py/compare_quadpyrmd.tikz}
    \caption{The Lebesgue constant (\textit{left}) and mass matrix condition number (\textit{right}) of the optimized nodes (\ref{line:lebcnst:0}), Stroud nodes (\ref{line:lebcnst:2}), Chan nodes (\ref{line:lebcnst:3}), and uniform nodes (\ref{line:lebcnst:1}) for the \textit{pyramid} element. The Lebesgue objective (\ref{eqn:lebobj}) for each nodal distribution is also included in the \textit{left} figure for each nodal distribution (\textit{dashed}).}
    \label{fig:compare_quadpyrmd}
\end{figure}

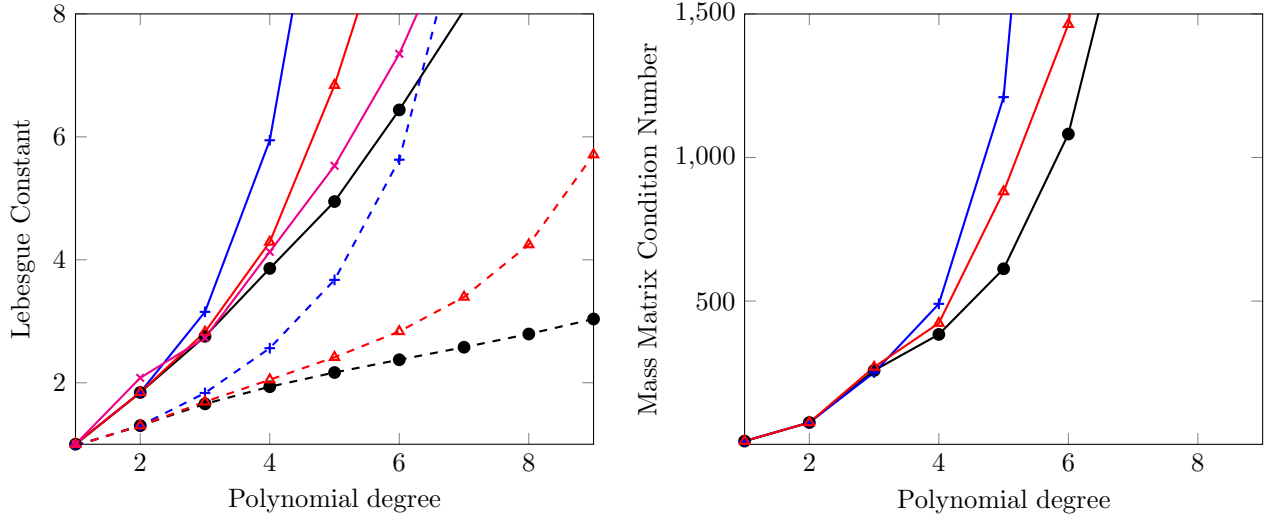
\begin{figure}[H]
    \centering
    \input{_py/compare_quadpyrmd_zoom.tikz}
    \caption{Zoom of Figure~\ref{fig:compare_quadpyrmd} to emphasize difference between optimized and Stroud nodes.}
    \label{fig:compare_quadpyrmd_zoom}
\end{figure}

The optimized distribution is compared with the Stroud \cite{hammer1956numerical}, Chan \cite{chan2015comparison}, and uniform distributions of the pyramid up to polynomial degree 9 by means of the Lebesgue constant, Lebesgue objective, and mass matrix condition number (Figures \ref{fig:compare_quadpyrmd}-\ref{fig:compare_quadpyrmd_zoom}). The optimized nodes achieve the lowest value for all these metrics across polynomial degrees. All nodal distributions satisfy similar constraints as the optimized nodes. All four distributions reduce to the corresponding 1D and 2D distributions on the edges and faces, respectively, of the square pyramid. Additionally, the Stroud-type, uniform, and Chan nodes are all symmetric.

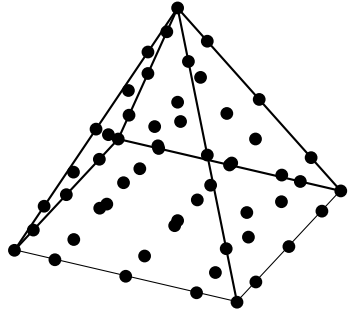
\begin{figure}[H]
    \centering
    \input{_py/nodes_quadpyrmd.tikz}
    \caption{Optimized nodal distribution for pyramid element with polynomial completeness $p = 4$.}
    \label{fig:nodes_quadpyrmd}
\end{figure}

The information in Figures~\ref{fig:compare_quadpyrmd}-\ref{fig:compare_quadpyrmd_zoom} is tabulated in Table~\ref{tab:quadpyrmd} and
Figure \ref{fig:nodes_quadpyrmd} shows a sample distribution of optimized nodes for $p=4$.

%% file: _py/compare_line.tikz
\begin{tikzpicture}
\begin{groupplot} [
group style={group size = 2 by 1, horizontal sep = 2.0cm}]
\nextgroupplot[xmin=1, xmax=30, ymode=log, xlabel={Polynomial degree}, ylabel={Lebesgue Constant}]
\addplot [black, thick, mark=*, mark size=2, mark options={solid}]
coordinates {
( 1.00000000e+00,  1.00000000e+00)
( 2.00000000e+00,  1.25000000e+00)
( 3.00000000e+00,  1.45541004e+00)
( 4.00000000e+00,  1.59218269e+00)
( 5.00000000e+00,  1.71883416e+00)
( 6.00000000e+00,  1.81440777e+00)
( 7.00000000e+00,  1.90505362e+00)
( 8.00000000e+00,  1.97857650e+00)
( 9.00000000e+00,  2.04914053e+00)
( 1.00000000e+01,  2.10888517e+00)
( 1.10000000e+01,  2.16665317e+00)
( 1.20000000e+01,  2.21692696e+00)
( 1.30000000e+01,  2.26586878e+00)
( 1.40000000e+01,  2.30928371e+00)
( 1.50000000e+01,  2.35171709e+00)
( 1.60000000e+01,  2.38995625e+00)
( 1.70000000e+01,  2.42737177e+00)
( 1.80000000e+01,  2.46151007e+00)
( 1.90000000e+01,  2.49499614e+00)
( 2.00000000e+01,  2.52581549e+00)
( 2.10000000e+01,  2.55613095e+00)
( 2.20000000e+01,  2.58391654e+00)
( 2.30000000e+01,  2.61191236e+00)
( 2.40000000e+01,  2.63762625e+00)
( 2.50000000e+01,  2.66320224e+00)
( 2.60000000e+01,  2.68678796e+00)
( 2.70000000e+01,  2.71066981e+00)
( 2.80000000e+01,  2.73261198e+00)
( 2.90000000e+01,  2.75484514e+00)
( 3.00000000e+01,  2.77497163e+00)};\label{line:lebcnst:0}

\addplot [black, thick, dashed, mark=*, mark size=2, mark options={solid}]
coordinates {
( 1.00000000e+00,  1.33333333e+00)
( 2.00000000e+00,  1.60000000e+00)
( 3.00000000e+00,  1.71075811e+00)
( 4.00000000e+00,  1.77260934e+00)
( 5.00000000e+00,  1.81236730e+00)
( 6.00000000e+00,  1.84016479e+00)
( 7.00000000e+00,  1.86072783e+00)
( 8.00000000e+00,  1.87657078e+00)
( 9.00000000e+00,  1.88915904e+00)
( 1.00000000e+01,  1.89940614e+00)
( 1.10000000e+01,  1.90791187e+00)
( 1.20000000e+01,  1.91508671e+00)
( 1.30000000e+01,  1.92122118e+00)
( 1.40000000e+01,  1.92652681e+00)
( 1.50000000e+01,  1.93116130e+00)
( 1.60000000e+01,  1.93524466e+00)
( 1.70000000e+01,  1.93886985e+00)
( 1.80000000e+01,  1.94211001e+00)
( 1.90000000e+01,  1.94502347e+00)
( 2.00000000e+01,  1.94765734e+00)
( 2.10000000e+01,  1.95005006e+00)
( 2.20000000e+01,  1.95223334e+00)
( 2.30000000e+01,  1.95423356e+00)
( 2.40000000e+01,  1.95607283e+00)
( 2.50000000e+01,  1.95776985e+00)
( 2.60000000e+01,  1.95934051e+00)
( 2.70000000e+01,  1.96079844e+00)
( 2.80000000e+01,  1.96215535e+00)
( 2.90000000e+01,  1.96342141e+00)
( 3.00000000e+01,  1.96460545e+00)};\label{line:lebobj:0}

\addplot [blue, thick, mark=+, mark size=2, mark options={solid}]
coordinates {
( 1.00000000e+00,  1.00000000e+00)
( 2.00000000e+00,  1.25000000e+00)
( 3.00000000e+00,  1.63112500e+00)
( 4.00000000e+00,  2.20782428e+00)
( 5.00000000e+00,  3.10627626e+00)
( 6.00000000e+00,  4.54934111e+00)
( 7.00000000e+00,  6.92944905e+00)
( 8.00000000e+00,  1.09450055e+01)
( 9.00000000e+00,  1.78476099e+01)
( 1.00000000e+01,  2.98981411e+01)
( 1.10000000e+01,  5.12141534e+01)
( 1.20000000e+01,  8.93237360e+01)
( 1.30000000e+01,  1.58086810e+02)
( 1.40000000e+01,  2.83180949e+02)
( 1.50000000e+01,  5.12347844e+02)
( 1.60000000e+01,  9.34273640e+02)
( 1.70000000e+01,  1.71575180e+03)
( 1.80000000e+01,  3.17033931e+03)
( 1.90000000e+01,  5.88654130e+03)
( 2.00000000e+01,  1.09792439e+04)
( 2.10000000e+01,  2.05736867e+04)
( 2.20000000e+01,  3.86668434e+04)
( 2.30000000e+01,  7.28406962e+04)
( 2.40000000e+01,  1.37851490e+05)
( 2.50000000e+01,  2.60966306e+05)
( 2.60000000e+01,  4.96482493e+05)
( 2.70000000e+01,  9.44687334e+05)
( 2.80000000e+01,  1.80244553e+06)
( 2.90000000e+01,  3.44766063e+06)
( 3.00000000e+01,  6.59251308e+06)};\label{line:lebcnst:1}

\addplot [blue, thick, dashed, mark=+, mark size=2, mark options={solid}]
coordinates {
( 1.00000000e+00,  1.33333333e+00)
( 2.00000000e+00,  1.60000000e+00)
( 3.00000000e+00,  1.84761905e+00)
( 4.00000000e+00,  2.13051146e+00)
( 5.00000000e+00,  2.52155483e+00)
( 6.00000000e+00,  3.16469930e+00)
( 7.00000000e+00,  4.38954141e+00)
( 8.00000000e+00,  7.00682271e+00)
( 9.00000000e+00,  1.31095178e+01)
( 1.00000000e+01,  2.83132471e+01)
( 1.10000000e+01,  6.81854792e+01)
( 1.20000000e+01,  1.77116363e+02)
( 1.30000000e+01,  4.84829486e+02)
( 1.40000000e+01,  1.37866245e+03)
( 1.50000000e+01,  4.03725924e+03)
( 1.60000000e+01,  1.21077779e+04)
( 1.70000000e+01,  3.70454799e+04)
( 1.80000000e+01,  1.15313289e+05)
( 1.90000000e+01,  3.64376426e+05)
( 2.00000000e+01,  1.16677291e+06)
( 2.10000000e+01,  3.78059076e+06)
( 2.20000000e+01,  1.23806571e+07)
( 2.30000000e+01,  4.09346539e+07)
( 2.40000000e+01,  1.36526858e+08)
( 2.50000000e+01,  4.58978536e+08)
( 2.60000000e+01,  1.55426212e+09)
( 2.70000000e+01,  5.29857992e+09)
( 2.80000000e+01,  1.81749686e+10)
( 2.90000000e+01,  6.26999890e+10)
( 3.00000000e+01,  2.17452432e+11)};\label{line:lebobj:1}

\addplot [red, thick, mark=triangle, mark size=2, mark options={solid}]
coordinates {
( 1.00000000e+00,  1.00000000e+00)
( 2.00000000e+00,  1.25000000e+00)
( 3.00000000e+00,  1.50000000e+00)
( 4.00000000e+00,  1.63587711e+00)
( 5.00000000e+00,  1.77859457e+00)
( 6.00000000e+00,  1.87373431e+00)
( 7.00000000e+00,  1.97236741e+00)
( 8.00000000e+00,  2.04563875e+00)
( 9.00000000e+00,  2.12096382e+00)
( 1.00000000e+01,  2.18049969e+00)
( 1.10000000e+01,  2.24147104e+00)
( 1.20000000e+01,  2.29163720e+00)
( 1.30000000e+01,  2.34281942e+00)
( 1.40000000e+01,  2.38618262e+00)
( 1.50000000e+01,  2.43026390e+00)
( 1.60000000e+01,  2.46828636e+00)
( 1.70000000e+01,  2.50715812e+00)
( 1.80000000e+01,  2.54106909e+00)
( 1.90000000e+01,  2.57577299e+00)
( 2.00000000e+01,  2.60646019e+00)
( 2.10000000e+01,  2.63771754e+00)
( 2.20000000e+01,  2.66577903e+00)
( 2.30000000e+01,  2.69417336e+00)
( 2.40000000e+01,  2.71997551e+00)
( 2.50000000e+01,  2.74603364e+00)
( 2.60000000e+01,  2.76957789e+00)
( 2.70000000e+01,  2.79399001e+00)
( 2.80000000e+01,  2.81562813e+00)
( 2.90000000e+01,  2.83858890e+00)
( 3.00000000e+01,  2.85913384e+00)};\label{line:lebcnst:2}

\addplot [red, thick, dashed, mark=triangle, mark size=2, mark options={solid}]
coordinates {
( 1.00000000e+00,  1.33333333e+00)
( 2.00000000e+00,  1.60000000e+00)
( 3.00000000e+00,  1.71428571e+00)
( 4.00000000e+00,  1.77777778e+00)
( 5.00000000e+00,  1.81818182e+00)
( 6.00000000e+00,  1.84615385e+00)
( 7.00000000e+00,  1.86666667e+00)
( 8.00000000e+00,  1.88235294e+00)
( 9.00000000e+00,  1.89473684e+00)
( 1.00000000e+01,  1.90476190e+00)
( 1.10000000e+01,  1.91304348e+00)
( 1.20000000e+01,  1.92000000e+00)
( 1.30000000e+01,  1.92592593e+00)
( 1.40000000e+01,  1.93103448e+00)
( 1.50000000e+01,  1.93548387e+00)
( 1.60000000e+01,  1.93939394e+00)
( 1.70000000e+01,  1.94285714e+00)
( 1.80000000e+01,  1.94594595e+00)
( 1.90000000e+01,  1.94871795e+00)
( 2.00000000e+01,  1.95121951e+00)
( 2.10000000e+01,  1.95348837e+00)
( 2.20000000e+01,  1.95555556e+00)
( 2.30000000e+01,  1.95744681e+00)
( 2.40000000e+01,  1.95918367e+00)
( 2.50000000e+01,  1.96078431e+00)
( 2.60000000e+01,  1.96226415e+00)
( 2.70000000e+01,  1.96363636e+00)
( 2.80000000e+01,  1.96491228e+00)
( 2.90000000e+01,  1.96610169e+00)
( 3.00000000e+01,  1.96721311e+00)};\label{line:lebobj:2}

\nextgroupplot[xmin=1, xmax=30, ymode=log, xlabel={Polynomial degree}, ylabel={Mass Matrix Condition Number}]
\addplot [black, thick, mark=*, mark size=2, mark options={solid}]
coordinates {
( 1.00000000e+00,  3.00000000e+00)
( 2.00000000e+00,  6.87964363e+00)
( 3.00000000e+00,  8.36166018e+00)
( 4.00000000e+00,  1.01661214e+01)
( 5.00000000e+00,  1.16825892e+01)
( 6.00000000e+00,  1.33409122e+01)
( 7.00000000e+00,  1.48818141e+01)
( 8.00000000e+00,  1.64992029e+01)
( 9.00000000e+00,  1.80525989e+01)
( 1.00000000e+01,  1.96535453e+01)
( 1.10000000e+01,  2.12141258e+01)
( 1.20000000e+01,  2.28070527e+01)
( 1.30000000e+01,  2.43721525e+01)
( 1.40000000e+01,  2.59606530e+01)
( 1.50000000e+01,  2.75287880e+01)
( 1.60000000e+01,  2.91146255e+01)
( 1.70000000e+01,  3.06849020e+01)
( 1.80000000e+01,  3.22690319e+01)
( 1.90000000e+01,  3.38408787e+01)
( 2.00000000e+01,  3.54238585e+01)
( 2.10000000e+01,  3.69968928e+01)
( 2.20000000e+01,  3.85790673e+01)
( 2.30000000e+01,  4.01530225e+01)
( 2.40000000e+01,  4.17346149e+01)
( 2.50000000e+01,  4.33092994e+01)
( 2.60000000e+01,  4.48904599e+01)
( 2.70000000e+01,  4.64657323e+01)
( 2.80000000e+01,  4.80465651e+01)
( 2.90000000e+01,  4.96223187e+01)
( 3.00000000e+01,  5.12028981e+01)};\label{line:Mcond:0}

\addplot [blue, thick, mark=+, mark size=2, mark options={solid}]
coordinates {
( 1.00000000e+00,  3.00000000e+00)
( 2.00000000e+00,  6.87964363e+00)
( 3.00000000e+00,  9.36100322e+00)
( 4.00000000e+00,  1.41390607e+01)
( 5.00000000e+00,  2.36484280e+01)
( 6.00000000e+00,  4.42129182e+01)
( 7.00000000e+00,  9.30355777e+01)
( 8.00000000e+00,  2.19667302e+02)
( 9.00000000e+00,  5.73439492e+02)
( 1.00000000e+01,  1.61965520e+03)
( 1.10000000e+01,  4.84506749e+03)
( 1.20000000e+01,  1.50950408e+04)
( 1.30000000e+01,  4.84183858e+04)
( 1.40000000e+01,  1.58713113e+05)
( 1.50000000e+01,  5.29233819e+05)
( 1.60000000e+01,  1.79009634e+06)
( 1.70000000e+01,  6.13045459e+06)
( 1.80000000e+01,  2.12280267e+07)
( 1.90000000e+01,  7.42386577e+07)
( 2.00000000e+01,  2.61925749e+08)
( 2.10000000e+01,  9.31280508e+08)
( 2.20000000e+01,  3.33337841e+09)
( 2.30000000e+01,  1.20003465e+10)
( 2.40000000e+01,  4.34205220e+10)
( 2.50000000e+01,  1.57821210e+11)
( 2.60000000e+01,  5.76043295e+11)
( 2.70000000e+01,  2.11086010e+12)
( 2.80000000e+01,  7.76394211e+12)
( 2.90000000e+01,  2.86569279e+13)
( 3.00000000e+01,  1.06109319e+14)};\label{line:Mcond:1}

\addplot [red, thick, mark=triangle, mark size=2, mark options={solid}]
coordinates {
( 1.00000000e+00,  3.00000000e+00)
( 2.00000000e+00,  6.87964363e+00)
( 3.00000000e+00,  8.65148372e+00)
( 4.00000000e+00,  1.04121778e+01)
( 5.00000000e+00,  1.19898257e+01)
( 6.00000000e+00,  1.36065948e+01)
( 7.00000000e+00,  1.51640408e+01)
( 8.00000000e+00,  1.67477726e+01)
( 9.00000000e+00,  1.83034409e+01)
( 1.00000000e+01,  1.98767363e+01)
( 1.10000000e+01,  2.14338487e+01)
( 1.20000000e+01,  2.30033532e+01)
( 1.30000000e+01,  2.45623368e+01)
( 1.40000000e+01,  2.61304224e+01)
( 1.50000000e+01,  2.76911055e+01)
( 1.60000000e+01,  2.92587205e+01)
( 1.70000000e+01,  3.08208306e+01)
( 1.80000000e+01,  3.23883776e+01)
( 1.90000000e+01,  3.39516614e+01)
( 2.00000000e+01,  3.55193136e+01)
( 2.10000000e+01,  3.70835607e+01)
( 2.20000000e+01,  3.86513887e+01)
( 2.30000000e+01,  4.02164304e+01)
( 2.40000000e+01,  4.17844573e+01)
( 2.50000000e+01,  4.33501593e+01)
( 2.60000000e+01,  4.49183858e+01)
( 2.70000000e+01,  4.64846413e+01)
( 2.80000000e+01,  4.80530580e+01)
( 2.90000000e+01,  4.96197813e+01)
( 3.00000000e+01,  5.11883747e+01)};\label{line:Mcond:2}

\end{groupplot}\end{tikzpicture}

%% file: _py/compare_line_zoom.tikz
\begin{tikzpicture}
\begin{groupplot} [
group style={group size = 2 by 1, horizontal sep = 2.0cm}]
\nextgroupplot[xmin=1, xmax=30, ymode=linear, xlabel={Polynomial degree}, ylabel={Lebesgue Constant}, ymin=1, ymax=3]
\addplot [black, thick, mark=*, mark size=2, mark options={solid}]
coordinates {
( 1.00000000e+00,  1.00000000e+00)
( 2.00000000e+00,  1.25000000e+00)
( 3.00000000e+00,  1.45541004e+00)
( 4.00000000e+00,  1.59218269e+00)
( 5.00000000e+00,  1.71883416e+00)
( 6.00000000e+00,  1.81440777e+00)
( 7.00000000e+00,  1.90505362e+00)
( 8.00000000e+00,  1.97857650e+00)
( 9.00000000e+00,  2.04914053e+00)
( 1.00000000e+01,  2.10888517e+00)
( 1.10000000e+01,  2.16665317e+00)
( 1.20000000e+01,  2.21692696e+00)
( 1.30000000e+01,  2.26586878e+00)
( 1.40000000e+01,  2.30928371e+00)
( 1.50000000e+01,  2.35171709e+00)
( 1.60000000e+01,  2.38995625e+00)
( 1.70000000e+01,  2.42737177e+00)
( 1.80000000e+01,  2.46151007e+00)
( 1.90000000e+01,  2.49499614e+00)
( 2.00000000e+01,  2.52581549e+00)
( 2.10000000e+01,  2.55613095e+00)
( 2.20000000e+01,  2.58391654e+00)
( 2.30000000e+01,  2.61191236e+00)
( 2.40000000e+01,  2.63762625e+00)
( 2.50000000e+01,  2.66320224e+00)
( 2.60000000e+01,  2.68678796e+00)
( 2.70000000e+01,  2.71066981e+00)
( 2.80000000e+01,  2.73261198e+00)
( 2.90000000e+01,  2.75484514e+00)
( 3.00000000e+01,  2.77497163e+00)};\label{line:lebcnst:0}

\addplot [black, thick, dashed, mark=*, mark size=2, mark options={solid}]
coordinates {
( 1.00000000e+00,  1.33333333e+00)
( 2.00000000e+00,  1.60000000e+00)
( 3.00000000e+00,  1.71075811e+00)
( 4.00000000e+00,  1.77260934e+00)
( 5.00000000e+00,  1.81236730e+00)
( 6.00000000e+00,  1.84016479e+00)
( 7.00000000e+00,  1.86072783e+00)
( 8.00000000e+00,  1.87657078e+00)
( 9.00000000e+00,  1.88915904e+00)
( 1.00000000e+01,  1.89940614e+00)
( 1.10000000e+01,  1.90791187e+00)
( 1.20000000e+01,  1.91508671e+00)
( 1.30000000e+01,  1.92122118e+00)
( 1.40000000e+01,  1.92652681e+00)
( 1.50000000e+01,  1.93116130e+00)
( 1.60000000e+01,  1.93524466e+00)
( 1.70000000e+01,  1.93886985e+00)
( 1.80000000e+01,  1.94211001e+00)
( 1.90000000e+01,  1.94502347e+00)
( 2.00000000e+01,  1.94765734e+00)
( 2.10000000e+01,  1.95005006e+00)
( 2.20000000e+01,  1.95223334e+00)
( 2.30000000e+01,  1.95423356e+00)
( 2.40000000e+01,  1.95607283e+00)
( 2.50000000e+01,  1.95776985e+00)
( 2.60000000e+01,  1.95934051e+00)
( 2.70000000e+01,  1.96079844e+00)
( 2.80000000e+01,  1.96215535e+00)
( 2.90000000e+01,  1.96342141e+00)
( 3.00000000e+01,  1.96460545e+00)};\label{line:lebobj:0}

\addplot [blue, thick, mark=+, mark size=2, mark options={solid}]
coordinates {
( 1.00000000e+00,  1.00000000e+00)
( 2.00000000e+00,  1.25000000e+00)
( 3.00000000e+00,  1.63112500e+00)
( 4.00000000e+00,  2.20782428e+00)
( 5.00000000e+00,  3.10627626e+00)};\label{line:lebcnst:1}

\addplot [blue, thick, dashed, mark=+, mark size=2, mark options={solid}]
coordinates {
( 1.00000000e+00,  1.33333333e+00)
( 2.00000000e+00,  1.60000000e+00)
( 3.00000000e+00,  1.84761905e+00)
( 4.00000000e+00,  2.13051146e+00)
( 5.00000000e+00,  2.52155483e+00)
( 6.00000000e+00,  3.16469930e+00)};\label{line:lebobj:1}

\addplot [red, thick, mark=triangle, mark size=2, mark options={solid}]
coordinates {
( 1.00000000e+00,  1.00000000e+00)
( 2.00000000e+00,  1.25000000e+00)
( 3.00000000e+00,  1.50000000e+00)
( 4.00000000e+00,  1.63587711e+00)
( 5.00000000e+00,  1.77859457e+00)
( 6.00000000e+00,  1.87373431e+00)
( 7.00000000e+00,  1.97236741e+00)
( 8.00000000e+00,  2.04563875e+00)
( 9.00000000e+00,  2.12096382e+00)
( 1.00000000e+01,  2.18049969e+00)
( 1.10000000e+01,  2.24147104e+00)
( 1.20000000e+01,  2.29163720e+00)
( 1.30000000e+01,  2.34281942e+00)
( 1.40000000e+01,  2.38618262e+00)
( 1.50000000e+01,  2.43026390e+00)
( 1.60000000e+01,  2.46828636e+00)
( 1.70000000e+01,  2.50715812e+00)
( 1.80000000e+01,  2.54106909e+00)
( 1.90000000e+01,  2.57577299e+00)
( 2.00000000e+01,  2.60646019e+00)
( 2.10000000e+01,  2.63771754e+00)
( 2.20000000e+01,  2.66577903e+00)
( 2.30000000e+01,  2.69417336e+00)
( 2.40000000e+01,  2.71997551e+00)
( 2.50000000e+01,  2.74603364e+00)
( 2.60000000e+01,  2.76957789e+00)
( 2.70000000e+01,  2.79399001e+00)
( 2.80000000e+01,  2.81562813e+00)
( 2.90000000e+01,  2.83858890e+00)
( 3.00000000e+01,  2.85913384e+00)};\label{line:lebcnst:2}

\addplot [red, thick, dashed, mark=triangle, mark size=2, mark options={solid}]
coordinates {
( 1.00000000e+00,  1.33333333e+00)
( 2.00000000e+00,  1.60000000e+00)
( 3.00000000e+00,  1.71428571e+00)
( 4.00000000e+00,  1.77777778e+00)
( 5.00000000e+00,  1.81818182e+00)
( 6.00000000e+00,  1.84615385e+00)
( 7.00000000e+00,  1.86666667e+00)
( 8.00000000e+00,  1.88235294e+00)
( 9.00000000e+00,  1.89473684e+00)
( 1.00000000e+01,  1.90476190e+00)
( 1.10000000e+01,  1.91304348e+00)
( 1.20000000e+01,  1.92000000e+00)
( 1.30000000e+01,  1.92592593e+00)
( 1.40000000e+01,  1.93103448e+00)
( 1.50000000e+01,  1.93548387e+00)
( 1.60000000e+01,  1.93939394e+00)
( 1.70000000e+01,  1.94285714e+00)
( 1.80000000e+01,  1.94594595e+00)
( 1.90000000e+01,  1.94871795e+00)
( 2.00000000e+01,  1.95121951e+00)
( 2.10000000e+01,  1.95348837e+00)
( 2.20000000e+01,  1.95555556e+00)
( 2.30000000e+01,  1.95744681e+00)
( 2.40000000e+01,  1.95918367e+00)
( 2.50000000e+01,  1.96078431e+00)
( 2.60000000e+01,  1.96226415e+00)
( 2.70000000e+01,  1.96363636e+00)
( 2.80000000e+01,  1.96491228e+00)
( 2.90000000e+01,  1.96610169e+00)
( 3.00000000e+01,  1.96721311e+00)};\label{line:lebobj:2}

\nextgroupplot[xmin=1, xmax=30, ymode=linear, xlabel={Polynomial degree}, ylabel={Mass Matrix Condition Number}, ymin=1, ymax=60]
\addplot [black, thick, mark=*, mark size=2, mark options={solid}]
coordinates {
( 1.00000000e+00,  3.00000000e+00)
( 2.00000000e+00,  6.87964363e+00)
( 3.00000000e+00,  8.36166018e+00)
( 4.00000000e+00,  1.01661214e+01)
( 5.00000000e+00,  1.16825892e+01)
( 6.00000000e+00,  1.33409122e+01)
( 7.00000000e+00,  1.48818141e+01)
( 8.00000000e+00,  1.64992029e+01)
( 9.00000000e+00,  1.80525989e+01)
( 1.00000000e+01,  1.96535453e+01)
( 1.10000000e+01,  2.12141258e+01)
( 1.20000000e+01,  2.28070527e+01)
( 1.30000000e+01,  2.43721525e+01)
( 1.40000000e+01,  2.59606530e+01)
( 1.50000000e+01,  2.75287880e+01)
( 1.60000000e+01,  2.91146255e+01)
( 1.70000000e+01,  3.06849020e+01)
( 1.80000000e+01,  3.22690319e+01)
( 1.90000000e+01,  3.38408787e+01)
( 2.00000000e+01,  3.54238585e+01)
( 2.10000000e+01,  3.69968928e+01)
( 2.20000000e+01,  3.85790673e+01)
( 2.30000000e+01,  4.01530225e+01)
( 2.40000000e+01,  4.17346149e+01)
( 2.50000000e+01,  4.33092994e+01)
( 2.60000000e+01,  4.48904599e+01)
( 2.70000000e+01,  4.64657323e+01)
( 2.80000000e+01,  4.80465651e+01)
( 2.90000000e+01,  4.96223187e+01)
( 3.00000000e+01,  5.12028981e+01)};\label{line:Mcond:0}

\addplot [blue, thick, mark=+, mark size=2, mark options={solid}]
coordinates {
( 1.00000000e+00,  3.00000000e+00)
( 2.00000000e+00,  6.87964363e+00)
( 3.00000000e+00,  9.36100322e+00)
( 4.00000000e+00,  1.41390607e+01)
( 5.00000000e+00,  2.36484280e+01)
( 6.00000000e+00,  4.42129182e+01)
( 7.00000000e+00,  9.30355777e+01)};\label{line:Mcond:1}

\addplot [red, thick, mark=triangle, mark size=2, mark options={solid}]
coordinates {
( 1.00000000e+00,  3.00000000e+00)
( 2.00000000e+00,  6.87964363e+00)
( 3.00000000e+00,  8.65148372e+00)
( 4.00000000e+00,  1.04121778e+01)
( 5.00000000e+00,  1.19898257e+01)
( 6.00000000e+00,  1.36065948e+01)
( 7.00000000e+00,  1.51640408e+01)
( 8.00000000e+00,  1.67477726e+01)
( 9.00000000e+00,  1.83034409e+01)
( 1.00000000e+01,  1.98767363e+01)
( 1.10000000e+01,  2.14338487e+01)
( 1.20000000e+01,  2.30033532e+01)
( 1.30000000e+01,  2.45623368e+01)
( 1.40000000e+01,  2.61304224e+01)
( 1.50000000e+01,  2.76911055e+01)
( 1.60000000e+01,  2.92587205e+01)
( 1.70000000e+01,  3.08208306e+01)
( 1.80000000e+01,  3.23883776e+01)
( 1.90000000e+01,  3.39516614e+01)
( 2.00000000e+01,  3.55193136e+01)
( 2.10000000e+01,  3.70835607e+01)
( 2.20000000e+01,  3.86513887e+01)
( 2.30000000e+01,  4.02164304e+01)
( 2.40000000e+01,  4.17844573e+01)
( 2.50000000e+01,  4.33501593e+01)
( 2.60000000e+01,  4.49183858e+01)
( 2.70000000e+01,  4.64846413e+01)
( 2.80000000e+01,  4.80530580e+01)
( 2.90000000e+01,  4.96197813e+01)
( 3.00000000e+01,  5.11883747e+01)};\label{line:Mcond:2}

\end{groupplot}\end{tikzpicture}

%% file: _py/nodes_line.tikz
\begin{tikzpicture}
\begin{axis}[
axis equal image,
xmin=-1,
xmax=1,
ymin=-0.2,
ymax=0.2,
axis lines=none,
width=0.8\textwidth]
\addplot [black, thick, forget plot]
coordinates {
(-1.00000000e+00,  0.00000000e+00)
( 1.00000000e+00,  0.00000000e+00)};

\addplot [black, thick, mark=*, mark size=2, mark options={solid}]
coordinates {
(-1.00000000e+00, -0.00000000e+00)
(-9.38430237e-01, -0.00000000e+00)
(-8.09737030e-01, -0.00000000e+00)
(-6.23406965e-01, -0.00000000e+00)
(-3.92717280e-01, -0.00000000e+00)
(-1.34085726e-01, -0.00000000e+00)
( 1.34085726e-01,  0.00000000e+00)
( 3.92717280e-01,  0.00000000e+00)
( 6.23406965e-01,  0.00000000e+00)
( 8.09737030e-01,  0.00000000e+00)
( 9.38430237e-01,  0.00000000e+00)
( 1.00000000e+00,  0.00000000e+00)};\label{line:nodes:0}

\addplot [red, thick, mark=triangle, mark size=2, mark options={solid}]
coordinates {
(-1.00000000e+00, -0.00000000e+00)
(-9.44899272e-01, -0.00000000e+00)
(-8.19279322e-01, -0.00000000e+00)
(-6.32876153e-01, -0.00000000e+00)
(-3.99530941e-01, -0.00000000e+00)
(-1.36552933e-01, -0.00000000e+00)
( 1.36552933e-01,  0.00000000e+00)
( 3.99530941e-01,  0.00000000e+00)
( 6.32876153e-01,  0.00000000e+00)
( 8.19279322e-01,  0.00000000e+00)
( 9.44899272e-01,  0.00000000e+00)
( 1.00000000e+00,  0.00000000e+00)};\label{line:nodes:2}

\end{axis}
\end{tikzpicture}

%% file: _py/compare_quad.tikz
\begin{tikzpicture}
\begin{groupplot} [
group style={group size = 2 by 1, horizontal sep = 2.0cm}]
\nextgroupplot[xmin=1, xmax=23, ymode=log, xlabel={Polynomial degree}, ylabel={Lebesgue Constant}]
\addplot [black, thick, mark=*, mark size=2, mark options={solid}]
coordinates {
( 1.00000000e+00,  1.00000000e+00)
( 2.00000000e+00,  1.56250000e+00)
( 3.00000000e+00,  2.11821839e+00)
( 4.00000000e+00,  2.53504572e+00)
( 5.00000000e+00,  2.95439086e+00)
( 6.00000000e+00,  3.29207556e+00)
( 7.00000000e+00,  3.62922930e+00)
( 8.00000000e+00,  3.91476495e+00)
( 9.00000000e+00,  4.19897690e+00)
( 1.00000000e+01,  4.44664583e+00)
( 1.10000000e+01,  4.69438598e+00)
( 1.20000000e+01,  4.91476514e+00)
( 1.30000000e+01,  5.13416132e+00)
( 1.40000000e+01,  5.33173389e+00)
( 1.50000000e+01,  5.53057326e+00)
( 1.60000000e+01,  5.70895052e+00)
( 1.70000000e+01,  5.89213367e+00)
( 1.80000000e+01,  6.05903149e+00)
( 1.90000000e+01,  6.22500568e+00)
( 2.00000000e+01,  6.37980464e+00)
( 2.10000000e+01,  6.53074854e+00)
( 2.20000000e+01,  6.67369730e+00)
( 2.30000000e+01,  6.75509350e+00)};\label{line:lebcnst:0}

\addplot [black, thick, dashed, mark=*, mark size=2, mark options={solid}]
coordinates {
( 1.00000000e+00,  1.77777778e+00)
( 2.00000000e+00,  2.56000000e+00)
( 3.00000000e+00,  2.92669330e+00)
( 4.00000000e+00,  3.14214387e+00)
( 5.00000000e+00,  3.28467525e+00)
( 6.00000000e+00,  3.38620646e+00)
( 7.00000000e+00,  3.46230805e+00)
( 8.00000000e+00,  3.52151788e+00)
( 9.00000000e+00,  3.56892189e+00)
( 1.00000000e+01,  3.60774368e+00)
( 1.10000000e+01,  3.64012772e+00)
( 1.20000000e+01,  3.66755711e+00)
( 1.30000000e+01,  3.69109083e+00)
( 1.40000000e+01,  3.71150556e+00)
( 1.50000000e+01,  3.72938397e+00)
( 1.60000000e+01,  3.74517190e+00)
( 1.70000000e+01,  3.75921631e+00)
( 1.80000000e+01,  3.77179130e+00)
( 1.90000000e+01,  3.78311631e+00)
( 2.00000000e+01,  3.79336909e+00)
( 2.10000000e+01,  3.80270545e+00)
( 2.20000000e+01,  3.81121667e+00)
( 2.30000000e+01,  3.82045784e+00)};\label{line:lebobj:0}

\addplot [blue, thick, mark=+, mark size=2, mark options={solid}]
coordinates {
( 1.00000000e+00,  1.00000000e+00)
( 2.00000000e+00,  1.56250000e+00)
( 3.00000000e+00,  2.66056877e+00)
( 4.00000000e+00,  4.87448804e+00)
( 5.00000000e+00,  9.64895219e+00)
( 6.00000000e+00,  2.06908509e+01)
( 7.00000000e+00,  4.80172641e+01)
( 8.00000000e+00,  1.19778816e+02)
( 9.00000000e+00,  3.18441477e+02)
( 1.00000000e+01,  8.93453667e+02)
( 1.10000000e+01,  2.61898081e+03)
( 1.20000000e+01,  7.97872981e+03)
( 1.30000000e+01,  2.49914394e+04)
( 1.40000000e+01,  8.01914501e+04)
( 1.50000000e+01,  2.62500314e+05)
( 1.60000000e+01,  8.71819498e+05)
( 1.70000000e+01,  2.94122225e+06)
( 1.80000000e+01,  1.00510539e+07)
( 1.90000000e+01,  3.46074358e+07)
( 2.00000000e+01,  1.20418449e+08)
( 2.10000000e+01,  4.23271809e+08)
( 2.20000000e+01,  1.48622304e+09)
( 2.30000000e+01,  5.29426363e+09)};\label{line:lebcnst:1}

\addplot [blue, thick, dashed, mark=+, mark size=2, mark options={solid}]
coordinates {
( 1.00000000e+00,  1.77777778e+00)
( 2.00000000e+00,  2.56000000e+00)
( 3.00000000e+00,  3.41369615e+00)
( 4.00000000e+00,  4.53907910e+00)
( 5.00000000e+00,  6.35823878e+00)
( 6.00000000e+00,  1.00153217e+01)
( 7.00000000e+00,  1.92680738e+01)
( 8.00000000e+00,  4.90955646e+01)
( 9.00000000e+00,  1.71859456e+02)
( 1.00000000e+01,  8.01639961e+02)
( 1.10000000e+01,  4.64925957e+03)
( 1.20000000e+01,  3.13702060e+04)
( 1.30000000e+01,  2.35059630e+05)
( 1.40000000e+01,  1.90071015e+06)
( 1.50000000e+01,  1.62994622e+07)
( 1.60000000e+01,  1.46598287e+08)
( 1.70000000e+01,  1.37236762e+09)
( 1.80000000e+01,  1.32971593e+10)
( 1.90000000e+01,  1.32769969e+11)
( 2.00000000e+01,  1.36135922e+12)
( 2.10000000e+01,  1.42920073e+13)
( 2.20000000e+01,  1.53269370e+14)
( 2.30000000e+01,  1.67194455e+15)};\label{line:lebobj:1}

\addplot [red, thick, mark=triangle, mark size=2, mark options={solid}]
coordinates {
( 1.00000000e+00,  1.00000000e+00)
( 2.00000000e+00,  1.56250000e+00)
( 3.00000000e+00,  2.25000000e+00)
( 4.00000000e+00,  2.67609393e+00)
( 5.00000000e+00,  3.16339864e+00)
( 6.00000000e+00,  3.51088027e+00)
( 7.00000000e+00,  3.89023322e+00)
( 8.00000000e+00,  4.18463789e+00)
( 9.00000000e+00,  4.49848755e+00)
( 1.00000000e+01,  4.75440498e+00)
( 1.10000000e+01,  5.02419244e+00)
( 1.20000000e+01,  5.25085803e+00)
( 1.30000000e+01,  5.48880283e+00)
( 1.40000000e+01,  5.69386749e+00)
( 1.50000000e+01,  5.90618264e+00)
( 1.60000000e+01,  6.09243755e+00)
( 1.70000000e+01,  6.28584183e+00)
( 1.80000000e+01,  6.45703218e+00)
( 1.90000000e+01,  6.63460660e+00)
( 2.00000000e+01,  6.79363474e+00)
( 2.10000000e+01,  6.95755365e+00)
( 2.20000000e+01,  7.09970120e+00)
( 2.30000000e+01,  7.25857072e+00)};\label{line:lebcnst:2}

\addplot [red, thick, dashed, mark=triangle, mark size=2, mark options={solid}]
coordinates {
( 1.00000000e+00,  1.77777778e+00)
( 2.00000000e+00,  2.56000000e+00)
( 3.00000000e+00,  2.93877551e+00)
( 4.00000000e+00,  3.16049383e+00)
( 5.00000000e+00,  3.30578512e+00)
( 6.00000000e+00,  3.40828402e+00)
( 7.00000000e+00,  3.48444444e+00)
( 8.00000000e+00,  3.54325260e+00)
( 9.00000000e+00,  3.59002770e+00)
( 1.00000000e+01,  3.62811791e+00)
( 1.10000000e+01,  3.65973535e+00)
( 1.20000000e+01,  3.68640000e+00)
( 1.30000000e+01,  3.70919067e+00)
( 1.40000000e+01,  3.72889417e+00)
( 1.50000000e+01,  3.74609781e+00)
( 1.60000000e+01,  3.76124885e+00)
( 1.70000000e+01,  3.77469388e+00)
( 1.80000000e+01,  3.78670562e+00)
( 1.90000000e+01,  3.79750164e+00)
( 2.00000000e+01,  3.80725758e+00)
( 2.10000000e+01,  3.81611680e+00)
( 2.20000000e+01,  3.82419754e+00)
( 2.30000000e+01,  3.83159800e+00)};\label{line:lebobj:2}

\nextgroupplot[xmin=1, xmax=23, ymode=log, xlabel={Polynomial degree}, ylabel={Mass Matrix Condition Number}]
\addplot [black, thick, mark=*, mark size=2, mark options={solid}]
coordinates {
( 1.00000000e+00,  9.00000000e+00)
( 2.00000000e+00,  4.73294965e+01)
( 3.00000000e+00,  6.99173610e+01)
( 4.00000000e+00,  1.03350025e+02)
( 5.00000000e+00,  1.36482889e+02)
( 6.00000000e+00,  1.77979939e+02)
( 7.00000000e+00,  2.21468392e+02)
( 8.00000000e+00,  2.72223697e+02)
( 9.00000000e+00,  3.25896327e+02)
( 1.00000000e+01,  3.86261842e+02)
( 1.10000000e+01,  4.50039134e+02)
( 1.20000000e+01,  5.20161652e+02)
( 1.30000000e+01,  5.94001816e+02)
( 1.40000000e+01,  6.73955505e+02)
( 1.50000000e+01,  7.57834167e+02)
( 1.60000000e+01,  8.47661415e+02)
( 1.70000000e+01,  9.41563210e+02)
( 1.80000000e+01,  1.04129040e+03)
( 1.90000000e+01,  1.14520505e+03)
( 2.00000000e+01,  1.25485534e+03)
( 2.10000000e+01,  1.36858329e+03)
( 2.20000000e+01,  1.48809427e+03)
( 2.30000000e+01,  1.60843605e+03)};\label{line:Mcond:0}

\addplot [blue, thick, mark=+, mark size=2, mark options={solid}]
coordinates {
( 1.00000000e+00,  9.00000000e+00)
( 2.00000000e+00,  4.73294965e+01)
( 3.00000000e+00,  8.76283814e+01)
( 4.00000000e+00,  1.99913038e+02)
( 5.00000000e+00,  5.59248145e+02)
( 6.00000000e+00,  1.95478214e+03)
( 7.00000000e+00,  8.65561872e+03)
( 8.00000000e+00,  4.82537234e+04)
( 9.00000000e+00,  3.28832851e+05)
( 1.00000000e+01,  2.62328298e+06)
( 1.10000000e+01,  2.34746790e+07)
( 1.20000000e+01,  2.27860256e+08)
( 1.30000000e+01,  2.34434010e+09)
( 1.40000000e+01,  2.51898568e+10)
( 1.50000000e+01,  2.80088506e+11)
( 1.60000000e+01,  3.20443896e+12)
( 1.70000000e+01,  3.75855665e+13)
( 1.80000000e+01,  4.50644952e+14)
( 1.90000000e+01,  6.58295380e+15)
( 2.00000000e+01,  7.28716130e+16)
( 2.10000000e+01,  9.35846925e+18)
( 2.20000000e+01,  5.19716627e+19)
( 2.30000000e+01,  1.54997506e+20)};\label{line:Mcond:1}

\addplot [red, thick, mark=triangle, mark size=2, mark options={solid}]
coordinates {
( 1.00000000e+00,  9.00000000e+00)
( 2.00000000e+00,  4.73294965e+01)
( 3.00000000e+00,  7.48481705e+01)
( 4.00000000e+00,  1.08413446e+02)
( 5.00000000e+00,  1.43755920e+02)
( 6.00000000e+00,  1.85139421e+02)
( 7.00000000e+00,  2.29948133e+02)
( 8.00000000e+00,  2.80487888e+02)
( 9.00000000e+00,  3.35015948e+02)
( 1.00000000e+01,  3.95084644e+02)
( 1.10000000e+01,  4.59409872e+02)
( 1.20000000e+01,  5.29154258e+02)
( 1.30000000e+01,  6.03308387e+02)
( 1.40000000e+01,  6.82798974e+02)
( 1.50000000e+01,  7.66797325e+02)
( 1.60000000e+01,  8.56072724e+02)
( 1.70000000e+01,  9.49923601e+02)
( 1.80000000e+01,  1.04900702e+03)
( 1.90000000e+01,  1.15271524e+03)
( 2.00000000e+01,  1.26162239e+03)
( 2.10000000e+01,  1.37518810e+03)
( 2.20000000e+01,  1.49392914e+03)
( 2.30000000e+01,  1.61737172e+03)};\label{line:Mcond:2}

\end{groupplot}\end{tikzpicture}

%% file: _py/compare_quad_zoom.tikz
\begin{tikzpicture}
\begin{groupplot} [
group style={group size = 2 by 1, horizontal sep = 2.0cm}]
\nextgroupplot[xmin=1, xmax=23, ymode=linear, xlabel={Polynomial degree}, ylabel={Lebesgue Constant}, ymin=1, ymax=8]
\addplot [black, thick, mark=*, mark size=2, mark options={solid}]
coordinates {
( 1.00000000e+00,  1.00000000e+00)
( 2.00000000e+00,  1.56250000e+00)
( 3.00000000e+00,  2.11821839e+00)
( 4.00000000e+00,  2.53504572e+00)
( 5.00000000e+00,  2.95439086e+00)
( 6.00000000e+00,  3.29207556e+00)
( 7.00000000e+00,  3.62922930e+00)
( 8.00000000e+00,  3.91476495e+00)
( 9.00000000e+00,  4.19897690e+00)
( 1.00000000e+01,  4.44664583e+00)
( 1.10000000e+01,  4.69438598e+00)
( 1.20000000e+01,  4.91476514e+00)
( 1.30000000e+01,  5.13416132e+00)
( 1.40000000e+01,  5.33173389e+00)
( 1.50000000e+01,  5.53057326e+00)
( 1.60000000e+01,  5.70895052e+00)
( 1.70000000e+01,  5.89213367e+00)
( 1.80000000e+01,  6.05903149e+00)
( 1.90000000e+01,  6.22500568e+00)
( 2.00000000e+01,  6.37980464e+00)
( 2.10000000e+01,  6.53074854e+00)
( 2.20000000e+01,  6.67369730e+00)
( 2.30000000e+01,  6.75509350e+00)};\label{line:lebcnst:0}

\addplot [black, thick, dashed, mark=*, mark size=2, mark options={solid}]
coordinates {
( 1.00000000e+00,  1.77777778e+00)
( 2.00000000e+00,  2.56000000e+00)
( 3.00000000e+00,  2.92669330e+00)
( 4.00000000e+00,  3.14214387e+00)
( 5.00000000e+00,  3.28467525e+00)
( 6.00000000e+00,  3.38620646e+00)
( 7.00000000e+00,  3.46230805e+00)
( 8.00000000e+00,  3.52151788e+00)
( 9.00000000e+00,  3.56892189e+00)
( 1.00000000e+01,  3.60774368e+00)
( 1.10000000e+01,  3.64012772e+00)
( 1.20000000e+01,  3.66755711e+00)
( 1.30000000e+01,  3.69109083e+00)
( 1.40000000e+01,  3.71150556e+00)
( 1.50000000e+01,  3.72938397e+00)
( 1.60000000e+01,  3.74517190e+00)
( 1.70000000e+01,  3.75921631e+00)
( 1.80000000e+01,  3.77179130e+00)
( 1.90000000e+01,  3.78311631e+00)
( 2.00000000e+01,  3.79336909e+00)
( 2.10000000e+01,  3.80270545e+00)
( 2.20000000e+01,  3.81121667e+00)
( 2.30000000e+01,  3.82045784e+00)};\label{line:lebobj:0}

\addplot [blue, thick, mark=+, mark size=2, mark options={solid}]
coordinates {
( 1.00000000e+00,  1.00000000e+00)
( 2.00000000e+00,  1.56250000e+00)
( 3.00000000e+00,  2.66056877e+00)
( 4.00000000e+00,  4.87448804e+00)
( 5.00000000e+00,  9.64895219e+00)};\label{line:lebcnst:1}

\addplot [blue, thick, dashed, mark=+, mark size=2, mark options={solid}]
coordinates {
( 1.00000000e+00,  1.77777778e+00)
( 2.00000000e+00,  2.56000000e+00)
( 3.00000000e+00,  3.41369615e+00)
( 4.00000000e+00,  4.53907910e+00)
( 5.00000000e+00,  6.35823878e+00)
( 6.00000000e+00,  1.00153217e+01)};\label{line:lebobj:1}

\addplot [red, thick, mark=triangle, mark size=2, mark options={solid}]
coordinates {
( 1.00000000e+00,  1.00000000e+00)
( 2.00000000e+00,  1.56250000e+00)
( 3.00000000e+00,  2.25000000e+00)
( 4.00000000e+00,  2.67609393e+00)
( 5.00000000e+00,  3.16339864e+00)
( 6.00000000e+00,  3.51088027e+00)
( 7.00000000e+00,  3.89023322e+00)
( 8.00000000e+00,  4.18463789e+00)
( 9.00000000e+00,  4.49848755e+00)
( 1.00000000e+01,  4.75440498e+00)
( 1.10000000e+01,  5.02419244e+00)
( 1.20000000e+01,  5.25085803e+00)
( 1.30000000e+01,  5.48880283e+00)
( 1.40000000e+01,  5.69386749e+00)
( 1.50000000e+01,  5.90618264e+00)
( 1.60000000e+01,  6.09243755e+00)
( 1.70000000e+01,  6.28584183e+00)
( 1.80000000e+01,  6.45703218e+00)
( 1.90000000e+01,  6.63460660e+00)
( 2.00000000e+01,  6.79363474e+00)
( 2.10000000e+01,  6.95755365e+00)
( 2.20000000e+01,  7.09970120e+00)
( 2.30000000e+01,  7.25857072e+00)};\label{line:lebcnst:2}

\addplot [red, thick, dashed, mark=triangle, mark size=2, mark options={solid}]
coordinates {
( 1.00000000e+00,  1.77777778e+00)
( 2.00000000e+00,  2.56000000e+00)
( 3.00000000e+00,  2.93877551e+00)
( 4.00000000e+00,  3.16049383e+00)
( 5.00000000e+00,  3.30578512e+00)
( 6.00000000e+00,  3.40828402e+00)
( 7.00000000e+00,  3.48444444e+00)
( 8.00000000e+00,  3.54325260e+00)
( 9.00000000e+00,  3.59002770e+00)
( 1.00000000e+01,  3.62811791e+00)
( 1.10000000e+01,  3.65973535e+00)
( 1.20000000e+01,  3.68640000e+00)
( 1.30000000e+01,  3.70919067e+00)
( 1.40000000e+01,  3.72889417e+00)
( 1.50000000e+01,  3.74609781e+00)
( 1.60000000e+01,  3.76124885e+00)
( 1.70000000e+01,  3.77469388e+00)
( 1.80000000e+01,  3.78670562e+00)
( 1.90000000e+01,  3.79750164e+00)
( 2.00000000e+01,  3.80725758e+00)
( 2.10000000e+01,  3.81611680e+00)
( 2.20000000e+01,  3.82419754e+00)
( 2.30000000e+01,  3.83159800e+00)};\label{line:lebobj:2}

\nextgroupplot[xmin=1, xmax=23, ymode=linear, xlabel={Polynomial degree}, ylabel={Mass Matrix Condition Number}, ymin=1, ymax=1700]
\addplot [black, thick, mark=*, mark size=2, mark options={solid}]
coordinates {
( 1.00000000e+00,  9.00000000e+00)
( 2.00000000e+00,  4.73294965e+01)
( 3.00000000e+00,  6.99173610e+01)
( 4.00000000e+00,  1.03350025e+02)
( 5.00000000e+00,  1.36482889e+02)
( 6.00000000e+00,  1.77979939e+02)
( 7.00000000e+00,  2.21468392e+02)
( 8.00000000e+00,  2.72223697e+02)
( 9.00000000e+00,  3.25896327e+02)
( 1.00000000e+01,  3.86261842e+02)
( 1.10000000e+01,  4.50039134e+02)
( 1.20000000e+01,  5.20161652e+02)
( 1.30000000e+01,  5.94001816e+02)
( 1.40000000e+01,  6.73955505e+02)
( 1.50000000e+01,  7.57834167e+02)
( 1.60000000e+01,  8.47661415e+02)
( 1.70000000e+01,  9.41563210e+02)
( 1.80000000e+01,  1.04129040e+03)
( 1.90000000e+01,  1.14520505e+03)
( 2.00000000e+01,  1.25485534e+03)
( 2.10000000e+01,  1.36858329e+03)
( 2.20000000e+01,  1.48809427e+03)
( 2.30000000e+01,  1.60843605e+03)};\label{line:Mcond:0}

\addplot [blue, thick, mark=+, mark size=2, mark options={solid}]
coordinates {
( 1.00000000e+00,  9.00000000e+00)
( 2.00000000e+00,  4.73294965e+01)
( 3.00000000e+00,  8.76283814e+01)
( 4.00000000e+00,  1.99913038e+02)
( 5.00000000e+00,  5.59248145e+02)
( 6.00000000e+00,  1.95478214e+03)};\label{line:Mcond:1}

\addplot [red, thick, mark=triangle, mark size=2, mark options={solid}]
coordinates {
( 1.00000000e+00,  9.00000000e+00)
( 2.00000000e+00,  4.73294965e+01)
( 3.00000000e+00,  7.48481705e+01)
( 4.00000000e+00,  1.08413446e+02)
( 5.00000000e+00,  1.43755920e+02)
( 6.00000000e+00,  1.85139421e+02)
( 7.00000000e+00,  2.29948133e+02)
( 8.00000000e+00,  2.80487888e+02)
( 9.00000000e+00,  3.35015948e+02)
( 1.00000000e+01,  3.95084644e+02)
( 1.10000000e+01,  4.59409872e+02)
( 1.20000000e+01,  5.29154258e+02)
( 1.30000000e+01,  6.03308387e+02)
( 1.40000000e+01,  6.82798974e+02)
( 1.50000000e+01,  7.66797325e+02)
( 1.60000000e+01,  8.56072724e+02)
( 1.70000000e+01,  9.49923601e+02)
( 1.80000000e+01,  1.04900702e+03)
( 1.90000000e+01,  1.15271524e+03)
( 2.00000000e+01,  1.26162239e+03)
( 2.10000000e+01,  1.37518810e+03)
( 2.20000000e+01,  1.49392914e+03)
( 2.30000000e+01,  1.61737172e+03)};\label{line:Mcond:2}

\end{groupplot}\end{tikzpicture}

%% file: _py/nodes_quad.tikz
\begin{tikzpicture}
\begin{axis}[
axis equal image,
xmin=-1,
xmax=1,
ymin=-1,
ymax=1,
axis lines=none,
width=0.8\textwidth]
\addplot [black, thick, forget plot]
coordinates {
(-1.00000000e+00, -1.00000000e+00)
( 1.00000000e+00, -1.00000000e+00)
( 1.00000000e+00,  1.00000000e+00)
(-1.00000000e+00,  1.00000000e+00)
(-1.00000000e+00, -1.00000000e+00)};

\addplot [black, thick, mark=*, mark size=2, mark options={solid}, only marks]
coordinates {
(-1.00000000e+00, -1.00000000e+00)
(-8.59806936e-01, -1.00000000e+00)
(-5.79014510e-01, -1.00000000e+00)
(-2.04062295e-01, -1.00000000e+00)
( 2.04062295e-01, -1.00000000e+00)
( 5.79014510e-01, -1.00000000e+00)
( 8.59806936e-01, -1.00000000e+00)
( 1.00000000e+00, -1.00000000e+00)
(-1.00000000e+00, -8.59806936e-01)
(-8.59806936e-01, -8.59806936e-01)
(-5.79014510e-01, -8.59806936e-01)
(-2.04062295e-01, -8.59806936e-01)
( 2.04062295e-01, -8.59806936e-01)
( 5.79014510e-01, -8.59806936e-01)
( 8.59806936e-01, -8.59806936e-01)
( 1.00000000e+00, -8.59806936e-01)
(-1.00000000e+00, -5.79014510e-01)
(-8.59806936e-01, -5.79014510e-01)
(-5.79014510e-01, -5.79014510e-01)
(-2.04062295e-01, -5.79014510e-01)
( 2.04062295e-01, -5.79014510e-01)
( 5.79014510e-01, -5.79014510e-01)
( 8.59806936e-01, -5.79014510e-01)
( 1.00000000e+00, -5.79014510e-01)
(-1.00000000e+00, -2.04062295e-01)
(-8.59806936e-01, -2.04062295e-01)
(-5.79014510e-01, -2.04062295e-01)
(-2.04062295e-01, -2.04062295e-01)
( 2.04062295e-01, -2.04062295e-01)
( 5.79014510e-01, -2.04062295e-01)
( 8.59806936e-01, -2.04062295e-01)
( 1.00000000e+00, -2.04062295e-01)
(-1.00000000e+00,  2.04062295e-01)
(-8.59806936e-01,  2.04062295e-01)
(-5.79014510e-01,  2.04062295e-01)
(-2.04062295e-01,  2.04062295e-01)
( 2.04062295e-01,  2.04062295e-01)
( 5.79014510e-01,  2.04062295e-01)
( 8.59806936e-01,  2.04062295e-01)
( 1.00000000e+00,  2.04062295e-01)
(-1.00000000e+00,  5.79014510e-01)
(-8.59806936e-01,  5.79014510e-01)
(-5.79014510e-01,  5.79014510e-01)
(-2.04062295e-01,  5.79014510e-01)
( 2.04062295e-01,  5.79014510e-01)
( 5.79014510e-01,  5.79014510e-01)
( 8.59806936e-01,  5.79014510e-01)
( 1.00000000e+00,  5.79014510e-01)
(-1.00000000e+00,  8.59806936e-01)
(-8.59806936e-01,  8.59806936e-01)
(-5.79014510e-01,  8.59806936e-01)
(-2.04062295e-01,  8.59806936e-01)
( 2.04062295e-01,  8.59806936e-01)
( 5.79014510e-01,  8.59806936e-01)
( 8.59806936e-01,  8.59806936e-01)
( 1.00000000e+00,  8.59806936e-01)
(-1.00000000e+00,  1.00000000e+00)
(-8.59806936e-01,  1.00000000e+00)
(-5.79014510e-01,  1.00000000e+00)
(-2.04062295e-01,  1.00000000e+00)
( 2.04062295e-01,  1.00000000e+00)
( 5.79014510e-01,  1.00000000e+00)
( 8.59806936e-01,  1.00000000e+00)
( 1.00000000e+00,  1.00000000e+00)};\label{line:nodes:0}

\addplot [red, thick, mark=triangle, mark size=2, mark options={solid}, only marks]
coordinates {
(-1.00000000e+00, -1.00000000e+00)
(-8.71740149e-01, -1.00000000e+00)
(-5.91700181e-01, -1.00000000e+00)
(-2.09299218e-01, -1.00000000e+00)
( 2.09299218e-01, -1.00000000e+00)
( 5.91700181e-01, -1.00000000e+00)
( 8.71740149e-01, -1.00000000e+00)
( 1.00000000e+00, -1.00000000e+00)
(-1.00000000e+00, -8.71740149e-01)
(-8.71740149e-01, -8.71740149e-01)
(-5.91700181e-01, -8.71740149e-01)
(-2.09299218e-01, -8.71740149e-01)
( 2.09299218e-01, -8.71740149e-01)
( 5.91700181e-01, -8.71740149e-01)
( 8.71740149e-01, -8.71740149e-01)
( 1.00000000e+00, -8.71740149e-01)
(-1.00000000e+00, -5.91700181e-01)
(-8.71740149e-01, -5.91700181e-01)
(-5.91700181e-01, -5.91700181e-01)
(-2.09299218e-01, -5.91700181e-01)
( 2.09299218e-01, -5.91700181e-01)
( 5.91700181e-01, -5.91700181e-01)
( 8.71740149e-01, -5.91700181e-01)
( 1.00000000e+00, -5.91700181e-01)
(-1.00000000e+00, -2.09299218e-01)
(-8.71740149e-01, -2.09299218e-01)
(-5.91700181e-01, -2.09299218e-01)
(-2.09299218e-01, -2.09299218e-01)
( 2.09299218e-01, -2.09299218e-01)
( 5.91700181e-01, -2.09299218e-01)
( 8.71740149e-01, -2.09299218e-01)
( 1.00000000e+00, -2.09299218e-01)
(-1.00000000e+00,  2.09299218e-01)
(-8.71740149e-01,  2.09299218e-01)
(-5.91700181e-01,  2.09299218e-01)
(-2.09299218e-01,  2.09299218e-01)
( 2.09299218e-01,  2.09299218e-01)
( 5.91700181e-01,  2.09299218e-01)
( 8.71740149e-01,  2.09299218e-01)
( 1.00000000e+00,  2.09299218e-01)
(-1.00000000e+00,  5.91700181e-01)
(-8.71740149e-01,  5.91700181e-01)
(-5.91700181e-01,  5.91700181e-01)
(-2.09299218e-01,  5.91700181e-01)
( 2.09299218e-01,  5.91700181e-01)
( 5.91700181e-01,  5.91700181e-01)
( 8.71740149e-01,  5.91700181e-01)
( 1.00000000e+00,  5.91700181e-01)
(-1.00000000e+00,  8.71740149e-01)
(-8.71740149e-01,  8.71740149e-01)
(-5.91700181e-01,  8.71740149e-01)
(-2.09299218e-01,  8.71740149e-01)
( 2.09299218e-01,  8.71740149e-01)
( 5.91700181e-01,  8.71740149e-01)
( 8.71740149e-01,  8.71740149e-01)
( 1.00000000e+00,  8.71740149e-01)
(-1.00000000e+00,  1.00000000e+00)
(-8.71740149e-01,  1.00000000e+00)
(-5.91700181e-01,  1.00000000e+00)
(-2.09299218e-01,  1.00000000e+00)
( 2.09299218e-01,  1.00000000e+00)
( 5.91700181e-01,  1.00000000e+00)
( 8.71740149e-01,  1.00000000e+00)
( 1.00000000e+00,  1.00000000e+00)};\label{line:nodes:2}

\end{axis}
\end{tikzpicture}

%% file: _py/compare_tri.tikz
\begin{tikzpicture}
\begin{groupplot} [
group style={group size = 2 by 1, horizontal sep = 2.0cm}]
\nextgroupplot[xmin=1, xmax=23, ymode=log, xlabel={Polynomial degree}, ylabel={Lebesgue Constant}]
\addplot [black, thick, mark=*, mark size=2, mark options={solid}]
coordinates {
( 1.00000000e+00,  1.00000000e+00)
( 2.00000000e+00,  1.66665600e+00)
( 3.00000000e+00,  2.12053917e+00)
( 4.00000000e+00,  2.68329914e+00)
( 5.00000000e+00,  3.26957407e+00)
( 6.00000000e+00,  3.77312297e+00)
( 7.00000000e+00,  4.39460769e+00)
( 8.00000000e+00,  5.10948909e+00)
( 9.00000000e+00,  5.94567682e+00)
( 1.00000000e+01,  7.11383048e+00)
( 1.10000000e+01,  8.38672666e+00)
( 1.20000000e+01,  1.01157432e+01)
( 1.30000000e+01,  1.21057922e+01)
( 1.40000000e+01,  1.47570928e+01)
( 1.50000000e+01,  1.78720400e+01)
( 1.60000000e+01,  2.28168957e+01)
( 1.70000000e+01,  2.90493512e+01)
( 1.80000000e+01,  3.85580126e+01)
( 1.90000000e+01,  5.14494912e+01)
( 2.00000000e+01,  7.04085926e+01)
( 2.10000000e+01,  1.27417088e+02)
( 2.20000000e+01,  2.27998141e+02)
( 2.30000000e+01,  6.80188896e+02)};\label{line:lebcnst:0}

\addplot [black, thick, dashed, mark=*, mark size=2, mark options={solid}]
coordinates {
( 1.00000000e+00,  1.00000000e+00)
( 2.00000000e+00,  1.26666667e+00)
( 3.00000000e+00,  1.48413839e+00)
( 4.00000000e+00,  1.64395526e+00)
( 5.00000000e+00,  1.77266746e+00)
( 6.00000000e+00,  1.88685663e+00)
( 7.00000000e+00,  1.99791676e+00)
( 8.00000000e+00,  2.11560993e+00)
( 9.00000000e+00,  2.25037229e+00)
( 1.00000000e+01,  2.41524936e+00)
( 1.10000000e+01,  2.62807500e+00)
( 1.20000000e+01,  2.91446340e+00)
( 1.30000000e+01,  3.31228772e+00)
( 1.40000000e+01,  3.87863299e+00)
( 1.50000000e+01,  4.70068541e+00)
( 1.60000000e+01,  5.89647830e+00)
( 1.70000000e+01,  7.53719050e+00)
( 1.80000000e+01,  1.20244873e+01)
( 1.90000000e+01,  1.76745391e+01)
( 2.00000000e+01,  2.64352433e+01)
( 2.10000000e+01,  3.89194232e+01)
( 2.20000000e+01,  1.04665263e+02)
( 2.30000000e+01,  5.46875900e+02)};\label{line:lebobj:0}

\addplot [blue, thick, mark=+, mark size=2, mark options={solid}]
coordinates {
( 1.00000000e+00,  1.00000000e+00)
( 2.00000000e+00,  1.66665600e+00)
( 3.00000000e+00,  2.26974357e+00)
( 4.00000000e+00,  3.47467034e+00)
( 5.00000000e+00,  5.45184000e+00)
( 6.00000000e+00,  8.74662857e+00)
( 7.00000000e+00,  1.43443011e+01)
( 8.00000000e+00,  2.40052679e+01)
( 9.00000000e+00,  4.09156167e+01)
( 1.00000000e+01,  7.08679805e+01)
( 1.10000000e+01,  1.24434087e+02)
( 1.20000000e+01,  2.21406921e+02)
( 1.30000000e+01,  3.97530469e+02)
( 1.40000000e+01,  7.20188440e+02)
( 1.50000000e+01,  1.31544927e+03)
( 1.60000000e+01,  2.41845440e+03)
( 1.70000000e+01,  4.46404265e+03)
( 1.80000000e+01,  8.34043971e+03)
( 1.90000000e+01,  1.60504343e+04)
( 2.00000000e+01,  5.42998623e+04)
( 2.10000000e+01,  1.64173684e+05)
( 2.20000000e+01,  7.22287818e+07)
( 2.30000000e+01,  1.29604349e+06)};\label{line:lebcnst:1}

\addplot [blue, thick, dashed, mark=+, mark size=2, mark options={solid}]
coordinates {
( 1.00000000e+00,  1.00000000e+00)
( 2.00000000e+00,  1.26666667e+00)
( 3.00000000e+00,  1.61071429e+00)
( 4.00000000e+00,  2.04550265e+00)
( 5.00000000e+00,  2.66092347e+00)
( 6.00000000e+00,  3.65555078e+00)
( 7.00000000e+00,  5.48062990e+00)
( 8.00000000e+00,  9.21573548e+00)
( 9.00000000e+00,  1.75632669e+01)
( 1.00000000e+01,  3.75566613e+01)
( 1.10000000e+01,  8.81310397e+01)
( 1.20000000e+01,  2.21788881e+02)
( 1.30000000e+01,  5.87918454e+02)
( 1.40000000e+01,  1.62135958e+03)
( 1.50000000e+01,  4.61354428e+03)
( 1.60000000e+01,  1.34688120e+04)
( 1.70000000e+01,  4.01843732e+04)
( 1.80000000e+01,  1.22140433e+05)
( 1.90000000e+01,  3.77173537e+05)
( 2.00000000e+01,  1.33589387e+06)
( 2.10000000e+01,  1.47792297e+07)
( 2.20000000e+01,  1.08186438e+12)
( 2.30000000e+01,  7.20504987e+08)};\label{line:lebobj:1}

\addplot [red, thick, mark=triangle, mark size=2, mark options={solid}]
coordinates {
( 1.00000000e+00,  1.00000000e+00)
( 2.00000000e+00,  1.66665600e+00)
( 3.00000000e+00,  2.11253387e+00)
( 4.00000000e+00,  2.67853115e+00)
( 5.00000000e+00,  3.40729685e+00)
( 6.00000000e+00,  3.90426304e+00)
( 7.00000000e+00,  4.47889830e+00)
( 8.00000000e+00,  5.10347297e+00)
( 9.00000000e+00,  5.87217923e+00)
( 1.00000000e+01,  6.77176125e+00)
( 1.10000000e+01,  8.04159220e+00)
( 1.20000000e+01,  9.49012599e+00)
( 1.30000000e+01,  1.16617858e+01)
( 1.40000000e+01,  1.42673273e+01)
( 1.50000000e+01,  1.80047775e+01)
( 1.60000000e+01,  2.28355813e+01)
( 1.70000000e+01,  2.96857231e+01)
( 1.80000000e+01,  3.86867506e+01)
( 1.90000000e+01,  5.16053538e+01)
( 2.00000000e+01,  7.86773769e+01)
( 2.10000000e+01,  1.92490242e+02)
( 2.20000000e+01,  3.12841374e+02)
( 2.30000000e+01,  7.91545865e+03)};\label{line:lebcnst:2}

\addplot [red, thick, dashed, mark=triangle, mark size=2, mark options={solid}]
coordinates {
( 1.00000000e+00,  1.00000000e+00)
( 2.00000000e+00,  1.26666667e+00)
( 3.00000000e+00,  1.48095238e+00)
( 4.00000000e+00,  1.63974203e+00)
( 5.00000000e+00,  1.77066623e+00)
( 6.00000000e+00,  1.88938413e+00)
( 7.00000000e+00,  2.00640394e+00)
( 8.00000000e+00,  2.13089033e+00)
( 9.00000000e+00,  2.27291376e+00)
( 1.00000000e+01,  2.44550372e+00)
( 1.10000000e+01,  2.66725793e+00)
( 1.20000000e+01,  2.96624698e+00)
( 1.30000000e+01,  3.38626987e+00)
( 1.40000000e+01,  3.99718343e+00)
( 1.50000000e+01,  4.91223340e+00)
( 1.60000000e+01,  6.31749938e+00)
( 1.70000000e+01,  8.52225478e+00)
( 1.80000000e+01,  1.20438320e+01)
( 1.90000000e+01,  1.77544840e+01)
( 2.00000000e+01,  2.68875806e+01)
( 2.10000000e+01,  6.98075726e+01)
( 2.20000000e+01,  1.40345434e+02)
( 2.30000000e+01,  5.35829637e+04)};\label{line:lebobj:2}

\addplot [magenta, thick, mark=x, mark size=2, mark options={solid}]
coordinates {
( 1.00000000e+00,  1.00000000e+00)
( 2.00000000e+00,  1.66665600e+00)
( 3.00000000e+00,  2.11253387e+00)
( 4.00000000e+00,  2.66220987e+00)
( 5.00000000e+00,  3.12114005e+00)
( 6.00000000e+00,  3.70175888e+00)
( 7.00000000e+00,  4.27451795e+00)
( 8.00000000e+00,  4.96231207e+00)
( 9.00000000e+00,  5.73642015e+00)
( 1.00000000e+01,  6.67058411e+00)
( 1.10000000e+01,  7.90172184e+00)
( 1.20000000e+01,  9.35902717e+00)
( 1.30000000e+01,  1.14660639e+01)
( 1.40000000e+01,  1.39676329e+01)
( 1.50000000e+01,  1.76344209e+01)
( 1.60000000e+01,  6.08492636e+01)
( 1.70000000e+01,  8.77067874e+01)
( 1.80000000e+01,  1.28685419e+02)
( 1.90000000e+01,  1.91620530e+02)
( 2.00000000e+01,  3.14685991e+02)
( 2.10000000e+01,  8.73580598e+02)
( 2.20000000e+01,  9.26631808e+03)
( 2.30000000e+01,  4.51752760e+04)};\label{line:lebcnst:3}

\addplot [magenta, thick, dashed, mark=x, mark size=2, mark options={solid}]
coordinates {
( 1.00000000e+00,  1.00000000e+00)
( 2.00000000e+00,  1.26666667e+00)
( 3.00000000e+00,  1.48095238e+00)
( 4.00000000e+00,  1.64133255e+00)
( 5.00000000e+00,  1.77392291e+00)
( 6.00000000e+00,  1.88773100e+00)
( 7.00000000e+00,  2.00199232e+00)
( 8.00000000e+00,  2.12190158e+00)
( 9.00000000e+00,  2.26059177e+00)
( 1.00000000e+01,  2.43330041e+00)
( 1.10000000e+01,  2.65265537e+00)
( 1.20000000e+01,  2.96225400e+00)
( 1.30000000e+01,  3.38481759e+00)
( 1.40000000e+01,  4.02493205e+00)
( 1.50000000e+01,  4.96774398e+00)
( 1.60000000e+01,  1.16828237e+01)
( 1.70000000e+01,  1.90044021e+01)
( 1.80000000e+01,  3.27834545e+01)
( 1.90000000e+01,  5.93331401e+01)
( 2.00000000e+01,  1.11991540e+02)
( 2.10000000e+01,  3.76547183e+02)
( 2.20000000e+01,  2.80199399e+04)
( 2.30000000e+01,  4.27060215e+05)};\label{line:lebobj:3}

\nextgroupplot[xmin=1, xmax=23, ymode=log, xlabel={Polynomial degree}, ylabel={Mass Matrix Condition Number}]
\addplot [black, thick, mark=*, mark size=2, mark options={solid}]
coordinates {
( 1.00000000e+00,  4.00000000e+00)
( 2.00000000e+00,  1.72085560e+01)
( 3.00000000e+00,  3.46658213e+01)
( 4.00000000e+00,  4.73433814e+01)
( 5.00000000e+00,  6.50769220e+01)
( 6.00000000e+00,  9.09019846e+01)
( 7.00000000e+00,  1.21203978e+02)
( 8.00000000e+00,  1.85925398e+02)
( 9.00000000e+00,  2.73197709e+02)
( 1.00000000e+01,  4.39390397e+02)
( 1.10000000e+01,  7.02058116e+02)
( 1.20000000e+01,  1.18189718e+03)
( 1.30000000e+01,  1.98413017e+03)
( 1.40000000e+01,  3.39059322e+03)
( 1.50000000e+01,  5.79642927e+03)
( 1.60000000e+01,  1.02974282e+04)
( 1.70000000e+01,  1.98526451e+04)
( 1.80000000e+01,  4.54279901e+04)
( 1.90000000e+01,  8.77235856e+04)
( 2.00000000e+01,  1.73242425e+05)
( 2.10000000e+01,  3.96104655e+05)
( 2.20000000e+01,  2.12958433e+06)
( 2.30000000e+01,  1.58852819e+07)};\label{line:Mcond:0}

\addplot [blue, thick, mark=+, mark size=2, mark options={solid}]
coordinates {
( 1.00000000e+00,  4.00000000e+00)
( 2.00000000e+00,  1.72085560e+01)
( 3.00000000e+00,  3.39692070e+01)
( 4.00000000e+00,  5.77440829e+01)
( 5.00000000e+00,  1.03365670e+02)
( 6.00000000e+00,  2.14866277e+02)
( 7.00000000e+00,  4.94568764e+02)
( 8.00000000e+00,  1.26925724e+03)
( 9.00000000e+00,  3.59387776e+03)
( 1.00000000e+01,  1.08500508e+04)
( 1.10000000e+01,  3.47992033e+04)
( 1.20000000e+01,  1.19009132e+05)
( 1.30000000e+01,  3.92743490e+05)
( 1.40000000e+01,  1.43677486e+06)
( 1.50000000e+01,  4.81530833e+06)
( 1.60000000e+01,  1.83524862e+07)
( 1.70000000e+01,  6.20754984e+07)
( 1.80000000e+01,  2.43013057e+08)
( 1.90000000e+01,  8.52548425e+08)
( 2.00000000e+01,  8.34598833e+09)
( 2.10000000e+01,  3.56433899e+11)
( 2.20000000e+01,  9.65862310e+16)
( 2.30000000e+01,  1.57202769e+13)};\label{line:Mcond:1}

\addplot [red, thick, mark=triangle, mark size=2, mark options={solid}]
coordinates {
( 1.00000000e+00,  4.00000000e+00)
( 2.00000000e+00,  1.72085560e+01)
( 3.00000000e+00,  3.48435211e+01)
( 4.00000000e+00,  4.70013374e+01)
( 5.00000000e+00,  6.91067172e+01)
( 6.00000000e+00,  9.21760479e+01)
( 7.00000000e+00,  1.27291932e+02)
( 8.00000000e+00,  1.95097425e+02)
( 9.00000000e+00,  2.89844741e+02)
( 1.00000000e+01,  4.69928392e+02)
( 1.10000000e+01,  7.48994061e+02)
( 1.20000000e+01,  1.27824945e+03)
( 1.30000000e+01,  2.17006751e+03)
( 1.40000000e+01,  3.88705166e+03)
( 1.50000000e+01,  6.98107749e+03)
( 1.60000000e+01,  1.30301287e+04)
( 1.70000000e+01,  2.44448352e+04)
( 1.80000000e+01,  4.72157065e+04)
( 1.90000000e+01,  9.13226940e+04)
( 2.00000000e+01,  2.09380424e+05)
( 2.10000000e+01,  1.46966403e+06)
( 2.20000000e+01,  3.27709386e+06)
( 2.30000000e+01,  2.60952967e+09)};\label{line:Mcond:2}

\addplot [magenta, thick, mark=x, mark size=2, mark options={solid}]
coordinates {
( 1.00000000e+00,  4.00000000e+00)
( 2.00000000e+00,  1.72085560e+01)
( 3.00000000e+00,  3.48435211e+01)
( 4.00000000e+00,  4.59259863e+01)
( 5.00000000e+00,  6.15440177e+01)
( 6.00000000e+00,  9.19927512e+01)
( 7.00000000e+00,  1.24539136e+02)
( 8.00000000e+00,  1.92824864e+02)
( 9.00000000e+00,  2.85488630e+02)
( 1.00000000e+01,  4.69591768e+02)
( 1.10000000e+01,  7.50792708e+02)
( 1.20000000e+01,  1.30554039e+03)
( 1.30000000e+01,  2.22964716e+03)
( 1.40000000e+01,  4.05505729e+03)
( 1.50000000e+01,  7.34319398e+03)
( 1.60000000e+01,  2.90486811e+04)
( 1.70000000e+01,  6.42142214e+04)
( 1.80000000e+01,  1.47620340e+05)
( 1.90000000e+01,  3.44396772e+05)
( 2.00000000e+01,  9.46551994e+05)
( 2.10000000e+01,  8.69455256e+06)
( 2.20000000e+01,  1.15175534e+09)
( 2.30000000e+01,  2.09004027e+10)};\label{line:Mcond:3}

\end{groupplot}\end{tikzpicture}

%% file: _py/compare_tri_zoom.tikz
\begin{tikzpicture}
\begin{groupplot} [
group style={group size = 2 by 1, horizontal sep = 2.0cm}]
\nextgroupplot[xmin=1, xmax=23, ymode=log, xlabel={Polynomial degree}, ylabel={Lebesgue Constant}, ymin=1, ymax=1000]
\addplot [black, thick, mark=*, mark size=2, mark options={solid}]
coordinates {
( 1.00000000e+00,  1.00000000e+00)
( 2.00000000e+00,  1.66665600e+00)
( 3.00000000e+00,  2.12053917e+00)
( 4.00000000e+00,  2.68329914e+00)
( 5.00000000e+00,  3.26957407e+00)
( 6.00000000e+00,  3.77312297e+00)
( 7.00000000e+00,  4.39460769e+00)
( 8.00000000e+00,  5.10948909e+00)
( 9.00000000e+00,  5.94567682e+00)
( 1.00000000e+01,  7.11383048e+00)
( 1.10000000e+01,  8.38672666e+00)
( 1.20000000e+01,  1.01157432e+01)
( 1.30000000e+01,  1.21057922e+01)
( 1.40000000e+01,  1.47570928e+01)
( 1.50000000e+01,  1.78720400e+01)
( 1.60000000e+01,  2.28168957e+01)
( 1.70000000e+01,  2.90493512e+01)
( 1.80000000e+01,  3.85580126e+01)
( 1.90000000e+01,  5.14494912e+01)
( 2.00000000e+01,  7.04085926e+01)
( 2.10000000e+01,  1.27417088e+02)
( 2.20000000e+01,  2.27998141e+02)
( 2.30000000e+01,  6.80188896e+02)};\label{line:lebcnst:0}

\addplot [black, thick, dashed, mark=*, mark size=2, mark options={solid}]
coordinates {
( 1.00000000e+00,  1.00000000e+00)
( 2.00000000e+00,  1.26666667e+00)
( 3.00000000e+00,  1.48413839e+00)
( 4.00000000e+00,  1.64395526e+00)
( 5.00000000e+00,  1.77266746e+00)
( 6.00000000e+00,  1.88685663e+00)
( 7.00000000e+00,  1.99791676e+00)
( 8.00000000e+00,  2.11560993e+00)
( 9.00000000e+00,  2.25037229e+00)
( 1.00000000e+01,  2.41524936e+00)
( 1.10000000e+01,  2.62807500e+00)
( 1.20000000e+01,  2.91446340e+00)
( 1.30000000e+01,  3.31228772e+00)
( 1.40000000e+01,  3.87863299e+00)
( 1.50000000e+01,  4.70068541e+00)
( 1.60000000e+01,  5.89647830e+00)
( 1.70000000e+01,  7.53719050e+00)
( 1.80000000e+01,  1.20244873e+01)
( 1.90000000e+01,  1.76745391e+01)
( 2.00000000e+01,  2.64352433e+01)
( 2.10000000e+01,  3.89194232e+01)
( 2.20000000e+01,  1.04665263e+02)
( 2.30000000e+01,  5.46875900e+02)};\label{line:lebobj:0}

\addplot [blue, thick, mark=+, mark size=2, mark options={solid}]
coordinates {
( 1.00000000e+00,  1.00000000e+00)
( 2.00000000e+00,  1.66665600e+00)
( 3.00000000e+00,  2.26974357e+00)
( 4.00000000e+00,  3.47467034e+00)
( 5.00000000e+00,  5.45184000e+00)
( 6.00000000e+00,  8.74662857e+00)
( 7.00000000e+00,  1.43443011e+01)
( 8.00000000e+00,  2.40052679e+01)
( 9.00000000e+00,  4.09156167e+01)
( 1.00000000e+01,  7.08679805e+01)
( 1.10000000e+01,  1.24434087e+02)
( 1.20000000e+01,  2.21406921e+02)
( 1.30000000e+01,  3.97530469e+02)
( 1.40000000e+01,  7.20188440e+02)
( 1.50000000e+01,  1.31544927e+03)};\label{line:lebcnst:1}

\addplot [blue, thick, dashed, mark=+, mark size=2, mark options={solid}]
coordinates {
( 1.00000000e+00,  1.00000000e+00)
( 2.00000000e+00,  1.26666667e+00)
( 3.00000000e+00,  1.61071429e+00)
( 4.00000000e+00,  2.04550265e+00)
( 5.00000000e+00,  2.66092347e+00)
( 6.00000000e+00,  3.65555078e+00)
( 7.00000000e+00,  5.48062990e+00)
( 8.00000000e+00,  9.21573548e+00)
( 9.00000000e+00,  1.75632669e+01)
( 1.00000000e+01,  3.75566613e+01)
( 1.10000000e+01,  8.81310397e+01)
( 1.20000000e+01,  2.21788881e+02)
( 1.30000000e+01,  5.87918454e+02)
( 1.40000000e+01,  1.62135958e+03)};\label{line:lebobj:1}

\addplot [red, thick, mark=triangle, mark size=2, mark options={solid}]
coordinates {
( 1.00000000e+00,  1.00000000e+00)
( 2.00000000e+00,  1.66665600e+00)
( 3.00000000e+00,  2.11253387e+00)
( 4.00000000e+00,  2.67853115e+00)
( 5.00000000e+00,  3.40729685e+00)
( 6.00000000e+00,  3.90426304e+00)
( 7.00000000e+00,  4.47889830e+00)
( 8.00000000e+00,  5.10347297e+00)
( 9.00000000e+00,  5.87217923e+00)
( 1.00000000e+01,  6.77176125e+00)
( 1.10000000e+01,  8.04159220e+00)
( 1.20000000e+01,  9.49012599e+00)
( 1.30000000e+01,  1.16617858e+01)
( 1.40000000e+01,  1.42673273e+01)
( 1.50000000e+01,  1.80047775e+01)
( 1.60000000e+01,  2.28355813e+01)
( 1.70000000e+01,  2.96857231e+01)
( 1.80000000e+01,  3.86867506e+01)
( 1.90000000e+01,  5.16053538e+01)
( 2.00000000e+01,  7.86773769e+01)
( 2.10000000e+01,  1.92490242e+02)
( 2.20000000e+01,  3.12841374e+02)
( 2.30000000e+01,  7.91545865e+03)};\label{line:lebcnst:2}

\addplot [red, thick, dashed, mark=triangle, mark size=2, mark options={solid}]
coordinates {
( 1.00000000e+00,  1.00000000e+00)
( 2.00000000e+00,  1.26666667e+00)
( 3.00000000e+00,  1.48095238e+00)
( 4.00000000e+00,  1.63974203e+00)
( 5.00000000e+00,  1.77066623e+00)
( 6.00000000e+00,  1.88938413e+00)
( 7.00000000e+00,  2.00640394e+00)
( 8.00000000e+00,  2.13089033e+00)
( 9.00000000e+00,  2.27291376e+00)
( 1.00000000e+01,  2.44550372e+00)
( 1.10000000e+01,  2.66725793e+00)
( 1.20000000e+01,  2.96624698e+00)
( 1.30000000e+01,  3.38626987e+00)
( 1.40000000e+01,  3.99718343e+00)
( 1.50000000e+01,  4.91223340e+00)
( 1.60000000e+01,  6.31749938e+00)
( 1.70000000e+01,  8.52225478e+00)
( 1.80000000e+01,  1.20438320e+01)
( 1.90000000e+01,  1.77544840e+01)
( 2.00000000e+01,  2.68875806e+01)
( 2.10000000e+01,  6.98075726e+01)
( 2.20000000e+01,  1.40345434e+02)
( 2.30000000e+01,  5.35829637e+04)};\label{line:lebobj:2}

\addplot [magenta, thick, mark=x, mark size=2, mark options={solid}]
coordinates {
( 1.00000000e+00,  1.00000000e+00)
( 2.00000000e+00,  1.66665600e+00)
( 3.00000000e+00,  2.11253387e+00)
( 4.00000000e+00,  2.66220987e+00)
( 5.00000000e+00,  3.12114005e+00)
( 6.00000000e+00,  3.70175888e+00)
( 7.00000000e+00,  4.27451795e+00)
( 8.00000000e+00,  4.96231207e+00)
( 9.00000000e+00,  5.73642015e+00)
( 1.00000000e+01,  6.67058411e+00)
( 1.10000000e+01,  7.90172184e+00)
( 1.20000000e+01,  9.35902717e+00)
( 1.30000000e+01,  1.14660639e+01)
( 1.40000000e+01,  1.39676329e+01)
( 1.50000000e+01,  1.76344209e+01)
( 1.60000000e+01,  6.08492636e+01)
( 1.70000000e+01,  8.77067874e+01)
( 1.80000000e+01,  1.28685419e+02)
( 1.90000000e+01,  1.91620530e+02)
( 2.00000000e+01,  3.14685991e+02)
( 2.10000000e+01,  8.73580598e+02)
( 2.20000000e+01,  9.26631808e+03)};\label{line:lebcnst:3}

\addplot [magenta, thick, dashed, mark=x, mark size=2, mark options={solid}]
coordinates {
( 1.00000000e+00,  1.00000000e+00)
( 2.00000000e+00,  1.26666667e+00)
( 3.00000000e+00,  1.48095238e+00)
( 4.00000000e+00,  1.64133255e+00)
( 5.00000000e+00,  1.77392291e+00)
( 6.00000000e+00,  1.88773100e+00)
( 7.00000000e+00,  2.00199232e+00)
( 8.00000000e+00,  2.12190158e+00)
( 9.00000000e+00,  2.26059177e+00)
( 1.00000000e+01,  2.43330041e+00)
( 1.10000000e+01,  2.65265537e+00)
( 1.20000000e+01,  2.96225400e+00)
( 1.30000000e+01,  3.38481759e+00)
( 1.40000000e+01,  4.02493205e+00)
( 1.50000000e+01,  4.96774398e+00)
( 1.60000000e+01,  1.16828237e+01)
( 1.70000000e+01,  1.90044021e+01)
( 1.80000000e+01,  3.27834545e+01)
( 1.90000000e+01,  5.93331401e+01)
( 2.00000000e+01,  1.11991540e+02)
( 2.10000000e+01,  3.76547183e+02)
( 2.20000000e+01,  2.80199399e+04)};\label{line:lebobj:3}

\nextgroupplot[xmin=1, xmax=23, ymode=log, xlabel={Polynomial degree}, ylabel={Mass Matrix Condition Number}, ymin=1, ymax=1000000000.0]
\addplot [black, thick, mark=*, mark size=2, mark options={solid}]
coordinates {
( 1.00000000e+00,  4.00000000e+00)
( 2.00000000e+00,  1.72085560e+01)
( 3.00000000e+00,  3.46658213e+01)
( 4.00000000e+00,  4.73433814e+01)
( 5.00000000e+00,  6.50769220e+01)
( 6.00000000e+00,  9.09019846e+01)
( 7.00000000e+00,  1.21203978e+02)
( 8.00000000e+00,  1.85925398e+02)
( 9.00000000e+00,  2.73197709e+02)
( 1.00000000e+01,  4.39390397e+02)
( 1.10000000e+01,  7.02058116e+02)
( 1.20000000e+01,  1.18189718e+03)
( 1.30000000e+01,  1.98413017e+03)
( 1.40000000e+01,  3.39059322e+03)
( 1.50000000e+01,  5.79642927e+03)
( 1.60000000e+01,  1.02974282e+04)
( 1.70000000e+01,  1.98526451e+04)
( 1.80000000e+01,  4.54279901e+04)
( 1.90000000e+01,  8.77235856e+04)
( 2.00000000e+01,  1.73242425e+05)
( 2.10000000e+01,  3.96104655e+05)
( 2.20000000e+01,  2.12958433e+06)
( 2.30000000e+01,  1.58852819e+07)};\label{line:Mcond:0}

\addplot [blue, thick, mark=+, mark size=2, mark options={solid}]
coordinates {
( 1.00000000e+00,  4.00000000e+00)
( 2.00000000e+00,  1.72085560e+01)
( 3.00000000e+00,  3.39692070e+01)
( 4.00000000e+00,  5.77440829e+01)
( 5.00000000e+00,  1.03365670e+02)
( 6.00000000e+00,  2.14866277e+02)
( 7.00000000e+00,  4.94568764e+02)
( 8.00000000e+00,  1.26925724e+03)
( 9.00000000e+00,  3.59387776e+03)
( 1.00000000e+01,  1.08500508e+04)
( 1.10000000e+01,  3.47992033e+04)
( 1.20000000e+01,  1.19009132e+05)
( 1.30000000e+01,  3.92743490e+05)
( 1.40000000e+01,  1.43677486e+06)
( 1.50000000e+01,  4.81530833e+06)
( 1.60000000e+01,  1.83524862e+07)
( 1.70000000e+01,  6.20754984e+07)
( 1.80000000e+01,  2.43013057e+08)
( 1.90000000e+01,  8.52548425e+08)
( 2.00000000e+01,  8.34598833e+09)};\label{line:Mcond:1}

\addplot [red, thick, mark=triangle, mark size=2, mark options={solid}]
coordinates {
( 1.00000000e+00,  4.00000000e+00)
( 2.00000000e+00,  1.72085560e+01)
( 3.00000000e+00,  3.48435211e+01)
( 4.00000000e+00,  4.70013374e+01)
( 5.00000000e+00,  6.91067172e+01)
( 6.00000000e+00,  9.21760479e+01)
( 7.00000000e+00,  1.27291932e+02)
( 8.00000000e+00,  1.95097425e+02)
( 9.00000000e+00,  2.89844741e+02)
( 1.00000000e+01,  4.69928392e+02)
( 1.10000000e+01,  7.48994061e+02)
( 1.20000000e+01,  1.27824945e+03)
( 1.30000000e+01,  2.17006751e+03)
( 1.40000000e+01,  3.88705166e+03)
( 1.50000000e+01,  6.98107749e+03)
( 1.60000000e+01,  1.30301287e+04)
( 1.70000000e+01,  2.44448352e+04)
( 1.80000000e+01,  4.72157065e+04)
( 1.90000000e+01,  9.13226940e+04)
( 2.00000000e+01,  2.09380424e+05)
( 2.10000000e+01,  1.46966403e+06)
( 2.20000000e+01,  3.27709386e+06)
( 2.30000000e+01,  2.60952967e+09)};\label{line:Mcond:2}

\addplot [magenta, thick, mark=x, mark size=2, mark options={solid}]
coordinates {
( 1.00000000e+00,  4.00000000e+00)
( 2.00000000e+00,  1.72085560e+01)
( 3.00000000e+00,  3.48435211e+01)
( 4.00000000e+00,  4.59259863e+01)
( 5.00000000e+00,  6.15440177e+01)
( 6.00000000e+00,  9.19927512e+01)
( 7.00000000e+00,  1.24539136e+02)
( 8.00000000e+00,  1.92824864e+02)
( 9.00000000e+00,  2.85488630e+02)
( 1.00000000e+01,  4.69591768e+02)
( 1.10000000e+01,  7.50792708e+02)
( 1.20000000e+01,  1.30554039e+03)
( 1.30000000e+01,  2.22964716e+03)
( 1.40000000e+01,  4.05505729e+03)
( 1.50000000e+01,  7.34319398e+03)
( 1.60000000e+01,  2.90486811e+04)
( 1.70000000e+01,  6.42142214e+04)
( 1.80000000e+01,  1.47620340e+05)
( 1.90000000e+01,  3.44396772e+05)
( 2.00000000e+01,  9.46551994e+05)
( 2.10000000e+01,  8.69455256e+06)
( 2.20000000e+01,  1.15175534e+09)};\label{line:Mcond:3}

\end{groupplot}\end{tikzpicture}

%% file: _py/nodes_tri.tikz
\begin{tikzpicture}
\begin{axis}[
axis equal image,
xmin=-1,
xmax=1,
ymin=-1,
ymax=1,
axis lines=none,
width=0.8\textwidth]
\addplot [black, thick, forget plot]
coordinates {
(-1.00000000e+00, -1.00000000e+00)
( 1.00000000e+00, -1.00000000e+00)
( 1.00000000e+00,  1.00000000e+00)
(-1.00000000e+00,  1.00000000e+00)
(-1.00000000e+00, -1.00000000e+00)};

\addplot [black, thick, mark=*, mark size=2, mark options={solid}, only marks]
coordinates {
(-1.00000000e+00, -1.00000000e+00)
(-8.59806936e-01, -1.00000000e+00)
(-5.79014510e-01, -1.00000000e+00)
(-2.04062295e-01, -1.00000000e+00)
( 2.04062295e-01, -1.00000000e+00)
( 5.79014510e-01, -1.00000000e+00)
( 8.59806936e-01, -1.00000000e+00)
( 1.00000000e+00, -1.00000000e+00)
(-1.00000000e+00, -8.59806936e-01)
(-8.35819865e-01, -8.35819865e-01)
(-5.09981269e-01, -8.06399803e-01)
(-1.01220867e-01, -7.97558266e-01)
( 3.16381072e-01, -8.06399803e-01)
( 6.71639729e-01, -8.35819865e-01)
( 8.59806936e-01, -8.59806936e-01)
(-1.00000000e+00, -5.79014510e-01)
(-8.06399803e-01, -5.09981269e-01)
(-4.67340918e-01, -4.67340918e-01)
(-6.53181647e-02, -4.67340918e-01)
( 3.16381072e-01, -5.09981269e-01)
( 5.79014510e-01, -5.79014510e-01)
(-1.00000000e+00, -2.04062295e-01)
(-7.97558266e-01, -1.01220867e-01)
(-4.67340918e-01, -6.53181647e-02)
(-1.01220867e-01, -1.01220867e-01)
( 2.04062295e-01, -2.04062295e-01)
(-1.00000000e+00,  2.04062295e-01)
(-8.06399803e-01,  3.16381072e-01)
(-5.09981269e-01,  3.16381072e-01)
(-2.04062295e-01,  2.04062295e-01)
(-1.00000000e+00,  5.79014510e-01)
(-8.35819865e-01,  6.71639729e-01)
(-5.79014510e-01,  5.79014510e-01)
(-1.00000000e+00,  8.59806936e-01)
(-8.59806936e-01,  8.59806936e-01)
(-1.00000000e+00,  1.00000000e+00)};\label{line:nodes:0}

\addplot [red, thick, mark=triangle, mark size=2, mark options={solid}, only marks]
coordinates {
(-1.00000000e+00, -1.00000000e+00)
(-8.71740149e-01, -1.00000000e+00)
(-5.91700181e-01, -1.00000000e+00)
(-2.09299218e-01, -1.00000000e+00)
( 2.09299218e-01, -1.00000000e+00)
( 5.91700181e-01, -1.00000000e+00)
( 8.71740149e-01, -1.00000000e+00)
( 1.00000000e+00, -1.00000000e+00)
(-1.00000000e+00, -8.71740149e-01)
(-8.25116806e-01, -8.25116806e-01)
(-4.97315697e-01, -8.09350441e-01)
(-9.73412792e-02, -8.05317442e-01)
( 3.06666138e-01, -8.09350441e-01)
( 6.50233612e-01, -8.25116806e-01)
( 8.71740149e-01, -8.71740149e-01)
(-1.00000000e+00, -5.91700181e-01)
(-8.09350441e-01, -4.97315697e-01)
(-4.65705246e-01, -4.65705246e-01)
(-6.85895084e-02, -4.65705246e-01)
( 3.06666138e-01, -4.97315697e-01)
( 5.91700181e-01, -5.91700181e-01)
(-1.00000000e+00, -2.09299218e-01)
(-8.05317442e-01, -9.73412792e-02)
(-4.65705246e-01, -6.85895084e-02)
(-9.73412792e-02, -9.73412792e-02)
( 2.09299218e-01, -2.09299218e-01)
(-1.00000000e+00,  2.09299218e-01)
(-8.09350441e-01,  3.06666138e-01)
(-4.97315697e-01,  3.06666138e-01)
(-2.09299218e-01,  2.09299218e-01)
(-1.00000000e+00,  5.91700181e-01)
(-8.25116806e-01,  6.50233612e-01)
(-5.91700181e-01,  5.91700181e-01)
(-1.00000000e+00,  8.71740149e-01)
(-8.71740149e-01,  8.71740149e-01)
(-1.00000000e+00,  1.00000000e+00)};\label{line:nodes:2}

\end{axis}
\end{tikzpicture}

%% file: _py/compare_hex.tikz
\begin{tikzpicture}
\begin{groupplot} [
group style={group size = 2 by 1, horizontal sep = 2.0cm}]
\nextgroupplot[xmin=1, xmax=9, ymode=log, xlabel={Polynomial degree}, ylabel={Lebesgue Constant}]
\addplot [black, thick, mark=*, mark size=2, mark options={solid}]
coordinates {
( 1.00000000e+00,  1.00000000e+00)
( 2.00000000e+00,  1.95312500e+00)
( 3.00000000e+00,  3.08287631e+00)
( 4.00000000e+00,  4.03605478e+00)
( 5.00000000e+00,  5.07810792e+00)
( 6.00000000e+00,  5.97207064e+00)
( 7.00000000e+00,  6.91387643e+00)
( 8.00000000e+00,  7.74566192e+00)
( 9.00000000e+00,  8.60429374e+00)};\label{line:lebcnst:0}

\addplot [black, thick, dashed, mark=*, mark size=2, mark options={solid}]
coordinates {
( 1.00000000e+00,  2.37037037e+00)
( 2.00000000e+00,  4.09600000e+00)
( 3.00000000e+00,  5.00686429e+00)
( 4.00000000e+00,  5.56979358e+00)
( 5.00000000e+00,  5.95303802e+00)
( 6.00000000e+00,  6.23117791e+00)
( 7.00000000e+00,  6.44241293e+00)
( 8.00000000e+00,  6.60837753e+00)
( 9.00000000e+00,  6.74226106e+00)};\label{line:lebobj:0}

\addplot [blue, thick, mark=+, mark size=2, mark options={solid}]
coordinates {
( 1.00000000e+00,  1.00000000e+00)
( 2.00000000e+00,  1.95312500e+00)
( 3.00000000e+00,  4.33972023e+00)
( 4.00000000e+00,  1.07608456e+01)
( 5.00000000e+00,  2.99342257e+01)
( 6.00000000e+00,  9.41168814e+01)
( 7.00000000e+00,  3.32733185e+02)
( 8.00000000e+00,  1.30863274e+03)
( 9.00000000e+00,  5.68512730e+03)};\label{line:lebcnst:1}

\addplot [blue, thick, dashed, mark=+, mark size=2, mark options={solid}]
coordinates {
( 1.00000000e+00,  2.37037037e+00)
( 2.00000000e+00,  4.09600000e+00)
( 3.00000000e+00,  6.30721002e+00)
( 4.00000000e+00,  9.67056005e+00)
( 5.00000000e+00,  1.60326477e+01)
( 6.00000000e+00,  3.16954815e+01)
( 7.00000000e+00,  8.45780079e+01)
( 8.00000000e+00,  3.44003917e+02)
( 9.00000000e+00,  2.25299459e+03)};\label{line:lebobj:1}

\addplot [red, thick, mark=triangle, mark size=2, mark options={solid}]
coordinates {
( 1.00000000e+00,  1.00000000e+00)
( 2.00000000e+00,  1.95312500e+00)
( 3.00000000e+00,  3.37500000e+00)
( 4.00000000e+00,  4.37776081e+00)
( 5.00000000e+00,  5.62640365e+00)
( 6.00000000e+00,  6.57845684e+00)
( 7.00000000e+00,  7.67296923e+00)
( 8.00000000e+00,  8.55353877e+00)
( 9.00000000e+00,  9.54112935e+00)};\label{line:lebcnst:2}

\addplot [red, thick, dashed, mark=triangle, mark size=2, mark options={solid}]
coordinates {
( 1.00000000e+00,  2.37037037e+00)
( 2.00000000e+00,  4.09600000e+00)
( 3.00000000e+00,  5.03790087e+00)
( 4.00000000e+00,  5.61865569e+00)
( 5.00000000e+00,  6.01051841e+00)
( 6.00000000e+00,  6.29221666e+00)
( 7.00000000e+00,  6.50429630e+00)
( 8.00000000e+00,  6.66965194e+00)
( 9.00000000e+00,  6.80215775e+00)};\label{line:lebobj:2}

\nextgroupplot[xmin=1, xmax=9, ymode=log, xlabel={Polynomial degree}, ylabel={Mass Matrix Condition Number}]
\addplot [black, thick, mark=*, mark size=2, mark options={solid}]
coordinates {
( 1.00000000e+00,  2.70000000e+01)
( 2.00000000e+00,  3.25610069e+02)
( 3.00000000e+00,  5.84625214e+02)
( 4.00000000e+00,  1.05066891e+03)
( 5.00000000e+00,  1.59447352e+03)
( 6.00000000e+00,  2.37441474e+03)
( 7.00000000e+00,  3.29585144e+03)
( 8.00000000e+00,  4.49147401e+03)
( 9.00000000e+00,  5.88327566e+03)};\label{line:Mcond:0}

\addplot [blue, thick, mark=+, mark size=2, mark options={solid}]
coordinates {
( 1.00000000e+00,  2.70000000e+01)
( 2.00000000e+00,  3.25610069e+02)
( 3.00000000e+00,  8.20289561e+02)
( 4.00000000e+00,  2.82658258e+03)
( 5.00000000e+00,  1.32253395e+04)
( 6.00000000e+00,  8.64266228e+04)
( 7.00000000e+00,  8.05280487e+05)
( 8.00000000e+00,  1.05997652e+07)
( 9.00000000e+00,  1.88565743e+08)};\label{line:Mcond:1}

\addplot [red, thick, mark=triangle, mark size=2, mark options={solid}]
coordinates {
( 1.00000000e+00,  2.70000000e+01)
( 2.00000000e+00,  3.25610069e+02)
( 3.00000000e+00,  6.47547728e+02)
( 4.00000000e+00,  1.12882007e+03)
( 5.00000000e+00,  1.72360843e+03)
( 6.00000000e+00,  2.51911708e+03)
( 7.00000000e+00,  3.48694288e+03)
( 8.00000000e+00,  4.69754736e+03)
( 9.00000000e+00,  6.13194459e+03)};\label{line:Mcond:2}

\end{groupplot}\end{tikzpicture}

%% file: _py/compare_hex_zoom.tikz
\begin{tikzpicture}
\begin{groupplot} [
group style={group size = 2 by 1, horizontal sep = 2.0cm}]
\nextgroupplot[xmin=1, xmax=9, ymode=linear, xlabel={Polynomial degree}, ylabel={Lebesgue Constant}, ymin=1, ymax=8]
\addplot [black, thick, mark=*, mark size=2, mark options={solid}]
coordinates {
( 1.00000000e+00,  1.00000000e+00)
( 2.00000000e+00,  1.95312500e+00)
( 3.00000000e+00,  3.08287631e+00)
( 4.00000000e+00,  4.03605478e+00)
( 5.00000000e+00,  5.07810792e+00)
( 6.00000000e+00,  5.97207064e+00)
( 7.00000000e+00,  6.91387643e+00)
( 8.00000000e+00,  7.74566192e+00)
( 9.00000000e+00,  8.60429374e+00)};\label{line:lebcnst:0}

\addplot [black, thick, dashed, mark=*, mark size=2, mark options={solid}]
coordinates {
( 1.00000000e+00,  2.37037037e+00)
( 2.00000000e+00,  4.09600000e+00)
( 3.00000000e+00,  5.00686429e+00)
( 4.00000000e+00,  5.56979358e+00)
( 5.00000000e+00,  5.95303802e+00)
( 6.00000000e+00,  6.23117791e+00)
( 7.00000000e+00,  6.44241293e+00)
( 8.00000000e+00,  6.60837753e+00)
( 9.00000000e+00,  6.74226106e+00)};\label{line:lebobj:0}

\addplot [blue, thick, mark=+, mark size=2, mark options={solid}]
coordinates {
( 1.00000000e+00,  1.00000000e+00)
( 2.00000000e+00,  1.95312500e+00)
( 3.00000000e+00,  4.33972023e+00)
( 4.00000000e+00,  1.07608456e+01)};\label{line:lebcnst:1}

\addplot [blue, thick, dashed, mark=+, mark size=2, mark options={solid}]
coordinates {
( 1.00000000e+00,  2.37037037e+00)
( 2.00000000e+00,  4.09600000e+00)
( 3.00000000e+00,  6.30721002e+00)
( 4.00000000e+00,  9.67056005e+00)};\label{line:lebobj:1}

\addplot [red, thick, mark=triangle, mark size=2, mark options={solid}]
coordinates {
( 1.00000000e+00,  1.00000000e+00)
( 2.00000000e+00,  1.95312500e+00)
( 3.00000000e+00,  3.37500000e+00)
( 4.00000000e+00,  4.37776081e+00)
( 5.00000000e+00,  5.62640365e+00)
( 6.00000000e+00,  6.57845684e+00)
( 7.00000000e+00,  7.67296923e+00)
( 8.00000000e+00,  8.55353877e+00)};\label{line:lebcnst:2}

\addplot [red, thick, dashed, mark=triangle, mark size=2, mark options={solid}]
coordinates {
( 1.00000000e+00,  2.37037037e+00)
( 2.00000000e+00,  4.09600000e+00)
( 3.00000000e+00,  5.03790087e+00)
( 4.00000000e+00,  5.61865569e+00)
( 5.00000000e+00,  6.01051841e+00)
( 6.00000000e+00,  6.29221666e+00)
( 7.00000000e+00,  6.50429630e+00)
( 8.00000000e+00,  6.66965194e+00)
( 9.00000000e+00,  6.80215775e+00)};\label{line:lebobj:2}

\nextgroupplot[xmin=1, xmax=9, ymode=linear, xlabel={Polynomial degree}, ylabel={Mass Matrix Condition Number}, ymin=1, ymax=3000]
\addplot [black, thick, mark=*, mark size=2, mark options={solid}]
coordinates {
( 1.00000000e+00,  2.70000000e+01)
( 2.00000000e+00,  3.25610069e+02)
( 3.00000000e+00,  5.84625214e+02)
( 4.00000000e+00,  1.05066891e+03)
( 5.00000000e+00,  1.59447352e+03)
( 6.00000000e+00,  2.37441474e+03)
( 7.00000000e+00,  3.29585144e+03)};\label{line:Mcond:0}

\addplot [blue, thick, mark=+, mark size=2, mark options={solid}]
coordinates {
( 1.00000000e+00,  2.70000000e+01)
( 2.00000000e+00,  3.25610069e+02)
( 3.00000000e+00,  8.20289561e+02)
( 4.00000000e+00,  2.82658258e+03)
( 5.00000000e+00,  1.32253395e+04)};\label{line:Mcond:1}

\addplot [red, thick, mark=triangle, mark size=2, mark options={solid}]
coordinates {
( 1.00000000e+00,  2.70000000e+01)
( 2.00000000e+00,  3.25610069e+02)
( 3.00000000e+00,  6.47547728e+02)
( 4.00000000e+00,  1.12882007e+03)
( 5.00000000e+00,  1.72360843e+03)
( 6.00000000e+00,  2.51911708e+03)
( 7.00000000e+00,  3.48694288e+03)};\label{line:Mcond:2}

\end{groupplot}\end{tikzpicture}

%% file: _py/nodes_hex.tikz
\begin{tikzpicture}
\begin{axis}[
axis equal image,
xmin=-1,
xmax=1,
ymin=-1,
ymax=1,
zmin=-1,
zmax=1,
axis lines=none,
width=0.8\textwidth]
\addplot3 [black, thick, forget plot]
coordinates {
(-1.00000000e+00, -1.00000000e+00, -1.00000000e+00)
( 1.00000000e+00, -1.00000000e+00, -1.00000000e+00)};

\addplot3 [black, thick, forget plot]
coordinates {
(-1.00000000e+00, -1.00000000e+00, -1.00000000e+00)
(-1.00000000e+00,  1.00000000e+00, -1.00000000e+00)};

\addplot3 [black, thick, forget plot]
coordinates {
(-1.00000000e+00,  1.00000000e+00, -1.00000000e+00)
( 1.00000000e+00,  1.00000000e+00, -1.00000000e+00)};

\addplot3 [black, thick, forget plot]
coordinates {
( 1.00000000e+00, -1.00000000e+00, -1.00000000e+00)
( 1.00000000e+00,  1.00000000e+00, -1.00000000e+00)};

\addplot3 [black, thick, forget plot]
coordinates {
(-1.00000000e+00, -1.00000000e+00,  1.00000000e+00)
( 1.00000000e+00, -1.00000000e+00,  1.00000000e+00)};

\addplot3 [black, thick, forget plot]
coordinates {
(-1.00000000e+00, -1.00000000e+00,  1.00000000e+00)
(-1.00000000e+00,  1.00000000e+00,  1.00000000e+00)};

\addplot3 [black, thick, forget plot]
coordinates {
(-1.00000000e+00,  1.00000000e+00,  1.00000000e+00)
( 1.00000000e+00,  1.00000000e+00,  1.00000000e+00)};

\addplot3 [black, thick, forget plot]
coordinates {
( 1.00000000e+00, -1.00000000e+00,  1.00000000e+00)
( 1.00000000e+00,  1.00000000e+00,  1.00000000e+00)};

\addplot3 [black, thick, forget plot]
coordinates {
(-1.00000000e+00, -1.00000000e+00, -1.00000000e+00)
(-1.00000000e+00, -1.00000000e+00,  1.00000000e+00)};

\addplot3 [black, thick, forget plot]
coordinates {
( 1.00000000e+00, -1.00000000e+00, -1.00000000e+00)
( 1.00000000e+00, -1.00000000e+00,  1.00000000e+00)};

\addplot3 [black, thick, forget plot]
coordinates {
(-1.00000000e+00,  1.00000000e+00, -1.00000000e+00)
(-1.00000000e+00,  1.00000000e+00,  1.00000000e+00)};

\addplot3 [black, thick, forget plot]
coordinates {
( 1.00000000e+00,  1.00000000e+00, -1.00000000e+00)
( 1.00000000e+00,  1.00000000e+00,  1.00000000e+00)};

\addplot3 [black, thick, mark=*, mark size=2, mark options={solid}, only marks]
coordinates {
(-1.00000000e+00, -1.00000000e+00, -1.00000000e+00)
(-6.36326016e-01, -1.00000000e+00, -1.00000000e+00)
(-0.00000000e+00, -1.00000000e+00, -1.00000000e+00)
( 6.36326016e-01, -1.00000000e+00, -1.00000000e+00)
( 1.00000000e+00, -1.00000000e+00, -1.00000000e+00)
(-1.00000000e+00, -6.36326016e-01, -1.00000000e+00)
(-6.36326016e-01, -6.36326016e-01, -1.00000000e+00)
(-0.00000000e+00, -6.36326016e-01, -1.00000000e+00)
( 6.36326016e-01, -6.36326016e-01, -1.00000000e+00)
( 1.00000000e+00, -6.36326016e-01, -1.00000000e+00)
(-1.00000000e+00, -0.00000000e+00, -1.00000000e+00)
(-6.36326016e-01, -0.00000000e+00, -1.00000000e+00)
(-0.00000000e+00, -0.00000000e+00, -1.00000000e+00)
( 6.36326016e-01, -0.00000000e+00, -1.00000000e+00)
( 1.00000000e+00, -0.00000000e+00, -1.00000000e+00)
(-1.00000000e+00,  6.36326016e-01, -1.00000000e+00)
(-6.36326016e-01,  6.36326016e-01, -1.00000000e+00)
(-0.00000000e+00,  6.36326016e-01, -1.00000000e+00)
( 6.36326016e-01,  6.36326016e-01, -1.00000000e+00)
( 1.00000000e+00,  6.36326016e-01, -1.00000000e+00)
(-1.00000000e+00,  1.00000000e+00, -1.00000000e+00)
(-6.36326016e-01,  1.00000000e+00, -1.00000000e+00)
(-0.00000000e+00,  1.00000000e+00, -1.00000000e+00)
( 6.36326016e-01,  1.00000000e+00, -1.00000000e+00)
( 1.00000000e+00,  1.00000000e+00, -1.00000000e+00)
(-1.00000000e+00, -1.00000000e+00, -6.36326016e-01)
(-6.36326016e-01, -1.00000000e+00, -6.36326016e-01)
(-0.00000000e+00, -1.00000000e+00, -6.36326016e-01)
( 6.36326016e-01, -1.00000000e+00, -6.36326016e-01)
( 1.00000000e+00, -1.00000000e+00, -6.36326016e-01)
(-1.00000000e+00, -6.36326016e-01, -6.36326016e-01)
(-6.36326016e-01, -6.36326016e-01, -6.36326016e-01)
(-0.00000000e+00, -6.36326016e-01, -6.36326016e-01)
( 6.36326016e-01, -6.36326016e-01, -6.36326016e-01)
( 1.00000000e+00, -6.36326016e-01, -6.36326016e-01)
(-1.00000000e+00, -0.00000000e+00, -6.36326016e-01)
(-6.36326016e-01, -0.00000000e+00, -6.36326016e-01)
(-0.00000000e+00, -0.00000000e+00, -6.36326016e-01)
( 6.36326016e-01, -0.00000000e+00, -6.36326016e-01)
( 1.00000000e+00, -0.00000000e+00, -6.36326016e-01)
(-1.00000000e+00,  6.36326016e-01, -6.36326016e-01)
(-6.36326016e-01,  6.36326016e-01, -6.36326016e-01)
(-0.00000000e+00,  6.36326016e-01, -6.36326016e-01)
( 6.36326016e-01,  6.36326016e-01, -6.36326016e-01)
( 1.00000000e+00,  6.36326016e-01, -6.36326016e-01)
(-1.00000000e+00,  1.00000000e+00, -6.36326016e-01)
(-6.36326016e-01,  1.00000000e+00, -6.36326016e-01)
(-0.00000000e+00,  1.00000000e+00, -6.36326016e-01)
( 6.36326016e-01,  1.00000000e+00, -6.36326016e-01)
( 1.00000000e+00,  1.00000000e+00, -6.36326016e-01)
(-1.00000000e+00, -1.00000000e+00, -0.00000000e+00)
(-6.36326016e-01, -1.00000000e+00, -0.00000000e+00)
(-0.00000000e+00, -1.00000000e+00, -0.00000000e+00)
( 6.36326016e-01, -1.00000000e+00, -0.00000000e+00)
( 1.00000000e+00, -1.00000000e+00, -0.00000000e+00)
(-1.00000000e+00, -6.36326016e-01, -0.00000000e+00)
(-6.36326016e-01, -6.36326016e-01, -0.00000000e+00)
(-0.00000000e+00, -6.36326016e-01, -0.00000000e+00)
( 6.36326016e-01, -6.36326016e-01, -0.00000000e+00)
( 1.00000000e+00, -6.36326016e-01, -0.00000000e+00)
(-1.00000000e+00, -0.00000000e+00, -0.00000000e+00)
(-6.36326016e-01, -0.00000000e+00, -0.00000000e+00)
( 0.00000000e+00,  0.00000000e+00,  0.00000000e+00)
( 6.36326016e-01, -0.00000000e+00, -0.00000000e+00)
( 1.00000000e+00, -0.00000000e+00, -0.00000000e+00)
(-1.00000000e+00,  6.36326016e-01, -0.00000000e+00)
(-6.36326016e-01,  6.36326016e-01, -0.00000000e+00)
(-0.00000000e+00,  6.36326016e-01, -0.00000000e+00)
( 6.36326016e-01,  6.36326016e-01, -0.00000000e+00)
( 1.00000000e+00,  6.36326016e-01, -0.00000000e+00)
(-1.00000000e+00,  1.00000000e+00, -0.00000000e+00)
(-6.36326016e-01,  1.00000000e+00, -0.00000000e+00)
(-0.00000000e+00,  1.00000000e+00, -0.00000000e+00)
( 6.36326016e-01,  1.00000000e+00, -0.00000000e+00)
( 1.00000000e+00,  1.00000000e+00, -0.00000000e+00)
(-1.00000000e+00, -1.00000000e+00,  6.36326016e-01)
(-6.36326016e-01, -1.00000000e+00,  6.36326016e-01)
(-0.00000000e+00, -1.00000000e+00,  6.36326016e-01)
( 6.36326016e-01, -1.00000000e+00,  6.36326016e-01)
( 1.00000000e+00, -1.00000000e+00,  6.36326016e-01)
(-1.00000000e+00, -6.36326016e-01,  6.36326016e-01)
(-6.36326016e-01, -6.36326016e-01,  6.36326016e-01)
(-0.00000000e+00, -6.36326016e-01,  6.36326016e-01)
( 6.36326016e-01, -6.36326016e-01,  6.36326016e-01)
( 1.00000000e+00, -6.36326016e-01,  6.36326016e-01)
(-1.00000000e+00, -0.00000000e+00,  6.36326016e-01)
(-6.36326016e-01, -0.00000000e+00,  6.36326016e-01)
(-0.00000000e+00, -0.00000000e+00,  6.36326016e-01)
( 6.36326016e-01, -0.00000000e+00,  6.36326016e-01)
( 1.00000000e+00, -0.00000000e+00,  6.36326016e-01)
(-1.00000000e+00,  6.36326016e-01,  6.36326016e-01)
(-6.36326016e-01,  6.36326016e-01,  6.36326016e-01)
(-0.00000000e+00,  6.36326016e-01,  6.36326016e-01)
( 6.36326016e-01,  6.36326016e-01,  6.36326016e-01)
( 1.00000000e+00,  6.36326016e-01,  6.36326016e-01)
(-1.00000000e+00,  1.00000000e+00,  6.36326016e-01)
(-6.36326016e-01,  1.00000000e+00,  6.36326016e-01)
(-0.00000000e+00,  1.00000000e+00,  6.36326016e-01)
( 6.36326016e-01,  1.00000000e+00,  6.36326016e-01)
( 1.00000000e+00,  1.00000000e+00,  6.36326016e-01)
(-1.00000000e+00, -1.00000000e+00,  1.00000000e+00)
(-6.36326016e-01, -1.00000000e+00,  1.00000000e+00)
(-0.00000000e+00, -1.00000000e+00,  1.00000000e+00)
( 6.36326016e-01, -1.00000000e+00,  1.00000000e+00)
( 1.00000000e+00, -1.00000000e+00,  1.00000000e+00)
(-1.00000000e+00, -6.36326016e-01,  1.00000000e+00)
(-6.36326016e-01, -6.36326016e-01,  1.00000000e+00)
(-0.00000000e+00, -6.36326016e-01,  1.00000000e+00)
( 6.36326016e-01, -6.36326016e-01,  1.00000000e+00)
( 1.00000000e+00, -6.36326016e-01,  1.00000000e+00)
(-1.00000000e+00, -0.00000000e+00,  1.00000000e+00)
(-6.36326016e-01, -0.00000000e+00,  1.00000000e+00)
(-0.00000000e+00, -0.00000000e+00,  1.00000000e+00)
( 6.36326016e-01, -0.00000000e+00,  1.00000000e+00)
( 1.00000000e+00, -0.00000000e+00,  1.00000000e+00)
(-1.00000000e+00,  6.36326016e-01,  1.00000000e+00)
(-6.36326016e-01,  6.36326016e-01,  1.00000000e+00)
(-0.00000000e+00,  6.36326016e-01,  1.00000000e+00)
( 6.36326016e-01,  6.36326016e-01,  1.00000000e+00)
( 1.00000000e+00,  6.36326016e-01,  1.00000000e+00)
(-1.00000000e+00,  1.00000000e+00,  1.00000000e+00)
(-6.36326016e-01,  1.00000000e+00,  1.00000000e+00)
(-0.00000000e+00,  1.00000000e+00,  1.00000000e+00)
( 6.36326016e-01,  1.00000000e+00,  1.00000000e+00)
( 1.00000000e+00,  1.00000000e+00,  1.00000000e+00)};\label{line:nodes:0}

\end{axis}
\end{tikzpicture}

%% file: _py/compare_tet.tikz
\begin{tikzpicture}
\begin{groupplot} [
group style={group size = 2 by 1, horizontal sep = 2.0cm}]
\nextgroupplot[xmin=1, xmax=9, ymode=log, xlabel={Polynomial degree}, ylabel={Lebesgue Constant}]
\addplot [black, thick, mark=*, mark size=2, mark options={solid}]
coordinates {
( 1.00000000e+00,  1.00000000e+00)
( 2.00000000e+00,  2.00000000e+00)
( 3.00000000e+00,  2.92941348e+00)
( 4.00000000e+00,  4.09065652e+00)
( 5.00000000e+00,  5.58774742e+00)
( 6.00000000e+00,  7.37566349e+00)
( 7.00000000e+00,  9.32999951e+00)
( 8.00000000e+00,  1.23286163e+01)
( 9.00000000e+00,  1.57216980e+01)};\label{line:lebcnst:0}

\addplot [black, thick, dashed, mark=*, mark size=2, mark options={solid}]
coordinates {
( 1.00000000e+00,  5.33333333e-01)
( 2.00000000e+00,  6.85714286e-01)
( 3.00000000e+00,  8.60260822e-01)
( 4.00000000e+00,  1.02450499e+00)
( 5.00000000e+00,  1.18675508e+00)
( 6.00000000e+00,  1.36485790e+00)
( 7.00000000e+00,  1.58466763e+00)
( 8.00000000e+00,  1.88594663e+00)
( 9.00000000e+00,  2.33356473e+00)};\label{line:lebobj:0}

\addplot [blue, thick, mark=+, mark size=2, mark options={solid}]
coordinates {
( 1.00000000e+00,  1.00000000e+00)
( 2.00000000e+00,  2.00000000e+00)
( 3.00000000e+00,  3.01978550e+00)
( 4.00000000e+00,  4.87842336e+00)
( 5.00000000e+00,  8.09171200e+00)
( 6.00000000e+00,  1.36566790e+01)
( 7.00000000e+00,  2.33618135e+01)
( 8.00000000e+00,  4.05449705e+01)
( 9.00000000e+00,  7.11124468e+01)};\label{line:lebcnst:1}

\addplot [blue, thick, dashed, mark=+, mark size=2, mark options={solid}]
coordinates {
( 1.00000000e+00,  5.33333333e-01)
( 2.00000000e+00,  6.85714286e-01)
( 3.00000000e+00,  9.19047619e-01)
( 4.00000000e+00,  1.24393138e+00)
( 5.00000000e+00,  1.72204945e+00)
( 6.00000000e+00,  2.49613986e+00)
( 7.00000000e+00,  3.88757545e+00)
( 8.00000000e+00,  6.64726051e+00)
( 9.00000000e+00,  1.26050666e+01)};\label{line:lebobj:1}

\addplot [red, thick, mark=triangle, mark size=2, mark options={solid}]
coordinates {
( 1.00000000e+00,  1.00000000e+00)
( 2.00000000e+00,  2.00000000e+00)
( 3.00000000e+00,  2.93253461e+00)
( 4.00000000e+00,  4.09289698e+00)
( 5.00000000e+00,  5.54450407e+00)
( 6.00000000e+00,  7.16890941e+00)
( 7.00000000e+00,  9.19199350e+00)
( 8.00000000e+00,  1.20512955e+01)
( 9.00000000e+00,  1.55835271e+01)};\label{line:lebcnst:2}

\addplot [red, thick, dashed, mark=triangle, mark size=2, mark options={solid}]
coordinates {
( 1.00000000e+00,  5.33333333e-01)
( 2.00000000e+00,  6.85714286e-01)
( 3.00000000e+00,  8.58201058e-01)
( 4.00000000e+00,  1.02347769e+00)
( 5.00000000e+00,  1.18842549e+00)
( 6.00000000e+00,  1.36886167e+00)
( 7.00000000e+00,  1.58868052e+00)
( 8.00000000e+00,  1.88619803e+00)
( 9.00000000e+00,  2.32882680e+00)};\label{line:lebobj:2}

\addplot [magenta, thick, mark=x, mark size=2, mark options={solid}]
coordinates {
( 1.00000000e+00,  1.00000000e+00)
( 2.00000000e+00,  2.00000000e+00)
( 3.00000000e+00,  2.93253461e+00)
( 4.00000000e+00,  4.07407306e+00)
( 5.00000000e+00,  5.31566031e+00)
( 6.00000000e+00,  7.00883352e+00)
( 7.00000000e+00,  9.20419258e+00)
( 8.00000000e+00,  1.25352735e+01)
( 9.00000000e+00,  1.70158309e+01)};\label{line:lebcnst:3}

\addplot [magenta, thick, dashed, mark=x, mark size=2, mark options={solid}]
coordinates {
( 1.00000000e+00,  5.33333333e-01)
( 2.00000000e+00,  6.85714286e-01)
( 3.00000000e+00,  8.58201058e-01)
( 4.00000000e+00,  1.02560667e+00)
( 5.00000000e+00,  1.19765146e+00)
( 6.00000000e+00,  1.39533966e+00)
( 7.00000000e+00,  1.67266675e+00)
( 8.00000000e+00,  2.08039666e+00)
( 9.00000000e+00,  2.74601451e+00)};\label{line:lebobj:3}

\nextgroupplot[xmin=1, xmax=9, ymode=log, xlabel={Polynomial degree}, ylabel={Mass Matrix Condition Number}]
\addplot [black, thick, mark=*, mark size=2, mark options={solid}]
coordinates {
( 1.00000000e+00,  5.00000000e+00)
( 2.00000000e+00,  3.59751792e+01)
( 3.00000000e+00,  1.09260378e+02)
( 4.00000000e+00,  2.45796017e+02)
( 5.00000000e+00,  3.40544589e+02)
( 6.00000000e+00,  6.42050609e+02)
( 7.00000000e+00,  1.33529528e+03)
( 8.00000000e+00,  3.03799029e+03)
( 9.00000000e+00,  7.15398210e+03)};\label{line:Mcond:0}

\addplot [blue, thick, mark=+, mark size=2, mark options={solid}]
coordinates {
( 1.00000000e+00,  5.00000000e+00)
( 2.00000000e+00,  3.59751792e+01)
( 3.00000000e+00,  1.05028447e+02)
( 4.00000000e+00,  2.36921262e+02)
( 5.00000000e+00,  4.09868811e+02)
( 6.00000000e+00,  1.08881444e+03)
( 7.00000000e+00,  2.73538625e+03)
( 8.00000000e+00,  8.45369188e+03)
( 9.00000000e+00,  2.63141068e+04)};\label{line:Mcond:1}

\addplot [red, thick, mark=triangle, mark size=2, mark options={solid}]
coordinates {
( 1.00000000e+00,  5.00000000e+00)
( 2.00000000e+00,  3.59751792e+01)
( 3.00000000e+00,  1.10353513e+02)
( 4.00000000e+00,  2.50164497e+02)
( 5.00000000e+00,  3.52168991e+02)
( 6.00000000e+00,  6.68291464e+02)
( 7.00000000e+00,  1.39418732e+03)
( 8.00000000e+00,  3.12532984e+03)
( 9.00000000e+00,  7.45030702e+03)};\label{line:Mcond:2}

\addplot [magenta, thick, mark=x, mark size=2, mark options={solid}]
coordinates {
( 1.00000000e+00,  5.00000000e+00)
( 2.00000000e+00,  3.59751792e+01)
( 3.00000000e+00,  1.10353513e+02)
( 4.00000000e+00,  2.50036908e+02)
( 5.00000000e+00,  3.62138789e+02)
( 6.00000000e+00,  6.92228918e+02)
( 7.00000000e+00,  1.38563216e+03)
( 8.00000000e+00,  3.17101356e+03)
( 9.00000000e+00,  7.80112146e+03)};\label{line:Mcond:3}

\end{groupplot}\end{tikzpicture}

%% file: _py/compare_tet_zoom.tikz
\begin{tikzpicture}
\begin{groupplot} [
group style={group size = 2 by 1, horizontal sep = 2.0cm}]
\nextgroupplot[xmin=1, xmax=9, ymode=linear, xlabel={Polynomial degree}, ylabel={Lebesgue Constant}, ymin=0, ymax=8]
\addplot [black, thick, mark=*, mark size=2, mark options={solid}]
coordinates {
( 1.00000000e+00,  1.00000000e+00)
( 2.00000000e+00,  2.00000000e+00)
( 3.00000000e+00,  2.92941348e+00)
( 4.00000000e+00,  4.09065652e+00)
( 5.00000000e+00,  5.58774742e+00)
( 6.00000000e+00,  7.37566349e+00)
( 7.00000000e+00,  9.32999951e+00)};\label{line:lebcnst:0}

\addplot [black, thick, dashed, mark=*, mark size=2, mark options={solid}]
coordinates {
( 1.00000000e+00,  5.33333333e-01)
( 2.00000000e+00,  6.85714286e-01)
( 3.00000000e+00,  8.60260822e-01)
( 4.00000000e+00,  1.02450499e+00)
( 5.00000000e+00,  1.18675508e+00)
( 6.00000000e+00,  1.36485790e+00)
( 7.00000000e+00,  1.58466763e+00)
( 8.00000000e+00,  1.88594663e+00)
( 9.00000000e+00,  2.33356473e+00)};\label{line:lebobj:0}

\addplot [blue, thick, mark=+, mark size=2, mark options={solid}]
coordinates {
( 1.00000000e+00,  1.00000000e+00)
( 2.00000000e+00,  2.00000000e+00)
( 3.00000000e+00,  3.01978550e+00)
( 4.00000000e+00,  4.87842336e+00)
( 5.00000000e+00,  8.09171200e+00)};\label{line:lebcnst:1}

\addplot [blue, thick, dashed, mark=+, mark size=2, mark options={solid}]
coordinates {
( 1.00000000e+00,  5.33333333e-01)
( 2.00000000e+00,  6.85714286e-01)
( 3.00000000e+00,  9.19047619e-01)
( 4.00000000e+00,  1.24393138e+00)
( 5.00000000e+00,  1.72204945e+00)
( 6.00000000e+00,  2.49613986e+00)
( 7.00000000e+00,  3.88757545e+00)
( 8.00000000e+00,  6.64726051e+00)
( 9.00000000e+00,  1.26050666e+01)};\label{line:lebobj:1}

\addplot [red, thick, mark=triangle, mark size=2, mark options={solid}]
coordinates {
( 1.00000000e+00,  1.00000000e+00)
( 2.00000000e+00,  2.00000000e+00)
( 3.00000000e+00,  2.93253461e+00)
( 4.00000000e+00,  4.09289698e+00)
( 5.00000000e+00,  5.54450407e+00)
( 6.00000000e+00,  7.16890941e+00)
( 7.00000000e+00,  9.19199350e+00)};\label{line:lebcnst:2}

\addplot [red, thick, dashed, mark=triangle, mark size=2, mark options={solid}]
coordinates {
( 1.00000000e+00,  5.33333333e-01)
( 2.00000000e+00,  6.85714286e-01)
( 3.00000000e+00,  8.58201058e-01)
( 4.00000000e+00,  1.02347769e+00)
( 5.00000000e+00,  1.18842549e+00)
( 6.00000000e+00,  1.36886167e+00)
( 7.00000000e+00,  1.58868052e+00)
( 8.00000000e+00,  1.88619803e+00)
( 9.00000000e+00,  2.32882680e+00)};\label{line:lebobj:2}

\addplot [magenta, thick, mark=x, mark size=2, mark options={solid}]
coordinates {
( 1.00000000e+00,  1.00000000e+00)
( 2.00000000e+00,  2.00000000e+00)
( 3.00000000e+00,  2.93253461e+00)
( 4.00000000e+00,  4.07407306e+00)
( 5.00000000e+00,  5.31566031e+00)
( 6.00000000e+00,  7.00883352e+00)
( 7.00000000e+00,  9.20419258e+00)};\label{line:lebcnst:3}

\addplot [magenta, thick, dashed, mark=x, mark size=2, mark options={solid}]
coordinates {
( 1.00000000e+00,  5.33333333e-01)
( 2.00000000e+00,  6.85714286e-01)
( 3.00000000e+00,  8.58201058e-01)
( 4.00000000e+00,  1.02560667e+00)
( 5.00000000e+00,  1.19765146e+00)
( 6.00000000e+00,  1.39533966e+00)
( 7.00000000e+00,  1.67266675e+00)
( 8.00000000e+00,  2.08039666e+00)
( 9.00000000e+00,  2.74601451e+00)};\label{line:lebobj:3}

\nextgroupplot[xmin=1, xmax=9, ymode=linear, xlabel={Polynomial degree}, ylabel={Mass Matrix Condition Number}, ymin=1, ymax=800]
\addplot [black, thick, mark=*, mark size=2, mark options={solid}]
coordinates {
( 1.00000000e+00,  5.00000000e+00)
( 2.00000000e+00,  3.59751792e+01)
( 3.00000000e+00,  1.09260378e+02)
( 4.00000000e+00,  2.45796017e+02)
( 5.00000000e+00,  3.40544589e+02)
( 6.00000000e+00,  6.42050609e+02)
( 7.00000000e+00,  1.33529528e+03)};\label{line:Mcond:0}

\addplot [blue, thick, mark=+, mark size=2, mark options={solid}]
coordinates {
( 1.00000000e+00,  5.00000000e+00)
( 2.00000000e+00,  3.59751792e+01)
( 3.00000000e+00,  1.05028447e+02)
( 4.00000000e+00,  2.36921262e+02)
( 5.00000000e+00,  4.09868811e+02)
( 6.00000000e+00,  1.08881444e+03)};\label{line:Mcond:1}

\addplot [red, thick, mark=triangle, mark size=2, mark options={solid}]
coordinates {
( 1.00000000e+00,  5.00000000e+00)
( 2.00000000e+00,  3.59751792e+01)
( 3.00000000e+00,  1.10353513e+02)
( 4.00000000e+00,  2.50164497e+02)
( 5.00000000e+00,  3.52168991e+02)
( 6.00000000e+00,  6.68291464e+02)
( 7.00000000e+00,  1.39418732e+03)};\label{line:Mcond:2}

\addplot [magenta, thick, mark=x, mark size=2, mark options={solid}]
coordinates {
( 1.00000000e+00,  5.00000000e+00)
( 2.00000000e+00,  3.59751792e+01)
( 3.00000000e+00,  1.10353513e+02)
( 4.00000000e+00,  2.50036908e+02)
( 5.00000000e+00,  3.62138789e+02)
( 6.00000000e+00,  6.92228918e+02)
( 7.00000000e+00,  1.38563216e+03)};\label{line:Mcond:3}

\end{groupplot}\end{tikzpicture}

%% file: _py/nodes_tet.tikz
\begin{tikzpicture}
\begin{axis}[
axis equal image,
xmin=-1,
xmax=1,
ymin=-1,
ymax=1,
zmin=-1,
zmax=1,
axis lines=none,
width=0.8\textwidth]
\addplot3 [black, thick, forget plot]
coordinates {
(-1.00000000e+00, -1.00000000e+00, -1.00000000e+00)
( 1.00000000e+00, -1.00000000e+00, -1.00000000e+00)};

\addplot3 [black, thick, forget plot]
coordinates {
(-1.00000000e+00, -1.00000000e+00, -1.00000000e+00)
(-1.00000000e+00,  1.00000000e+00, -1.00000000e+00)};

\addplot3 [black, thick, forget plot]
coordinates {
(-1.00000000e+00, -1.00000000e+00, -1.00000000e+00)
(-1.00000000e+00, -1.00000000e+00,  1.00000000e+00)};

\addplot3 [black, thick, forget plot]
coordinates {
( 1.00000000e+00, -1.00000000e+00, -1.00000000e+00)
(-1.00000000e+00,  1.00000000e+00, -1.00000000e+00)};

\addplot3 [black, thick, forget plot]
coordinates {
(-1.00000000e+00, -1.00000000e+00,  1.00000000e+00)
(-1.00000000e+00,  1.00000000e+00, -1.00000000e+00)};

\addplot3 [black, thick, forget plot]
coordinates {
(-1.00000000e+00, -1.00000000e+00,  1.00000000e+00)
( 1.00000000e+00, -1.00000000e+00, -1.00000000e+00)};

\addplot3 [black, thick, mark=*, mark size=2, mark options={solid}, only marks]
coordinates {
(-1.00000000e+00, -1.00000000e+00, -1.00000000e+00)
(-6.36326016e-01, -1.00000000e+00, -1.00000000e+00)
( 0.00000000e+00, -1.00000000e+00, -1.00000000e+00)
( 6.36326016e-01, -1.00000000e+00, -1.00000000e+00)
( 1.00000000e+00, -1.00000000e+00, -1.00000000e+00)
(-1.00000000e+00, -6.36326016e-01, -1.00000000e+00)
(-5.57815335e-01, -5.57815335e-01, -1.00000000e+00)
( 1.15630670e-01, -5.57815335e-01, -1.00000000e+00)
( 6.36326016e-01, -6.36326016e-01, -1.00000000e+00)
(-1.00000000e+00,  0.00000000e+00, -1.00000000e+00)
(-5.57815335e-01,  1.15630670e-01, -1.00000000e+00)
( 0.00000000e+00,  0.00000000e+00, -1.00000000e+00)
(-1.00000000e+00,  6.36326016e-01, -1.00000000e+00)
(-6.36326016e-01,  6.36326016e-01, -1.00000000e+00)
(-1.00000000e+00,  1.00000000e+00, -1.00000000e+00)
(-1.00000000e+00, -1.00000000e+00, -6.36326016e-01)
(-5.57815335e-01, -1.00000000e+00, -5.57815335e-01)
( 1.15630670e-01, -1.00000000e+00, -5.57815335e-01)
( 6.36326016e-01, -1.00000000e+00, -6.36326016e-01)
(-1.00000000e+00, -5.57815335e-01, -5.57815335e-01)
(-5.00000000e-01, -5.00000000e-01, -5.00000000e-01)
( 1.15630670e-01, -5.57815335e-01, -5.57815335e-01)
(-1.00000000e+00,  1.15630670e-01, -5.57815335e-01)
(-5.57815335e-01,  1.15630670e-01, -5.57815335e-01)
(-1.00000000e+00,  6.36326016e-01, -6.36326016e-01)
(-1.00000000e+00, -1.00000000e+00,  0.00000000e+00)
(-5.57815335e-01, -1.00000000e+00,  1.15630670e-01)
( 0.00000000e+00, -1.00000000e+00,  0.00000000e+00)
(-1.00000000e+00, -5.57815335e-01,  1.15630670e-01)
(-5.57815335e-01, -5.57815335e-01,  1.15630670e-01)
(-1.00000000e+00,  0.00000000e+00,  0.00000000e+00)
(-1.00000000e+00, -1.00000000e+00,  6.36326016e-01)
(-6.36326016e-01, -1.00000000e+00,  6.36326016e-01)
(-1.00000000e+00, -6.36326016e-01,  6.36326016e-01)
(-1.00000000e+00, -1.00000000e+00,  1.00000000e+00)};\label{line:nodes:0}

\end{axis}
\end{tikzpicture}

%% file: _py/compare_triprism.tikz
\begin{tikzpicture}
\begin{groupplot} [
group style={group size = 2 by 1, horizontal sep = 2.0cm}]
\nextgroupplot[xmin=1, xmax=9, ymode=log, xlabel={Polynomial degree}, ylabel={Lebesgue Constant}]
\addplot [black, thick, mark=*, mark size=2, mark options={solid}]
coordinates {
( 1.00000000e+00,  1.00000000e+00)
( 2.00000000e+00,  2.08325000e+00)
( 3.00000000e+00,  3.08620784e+00)
( 4.00000000e+00,  4.27186954e+00)
( 5.00000000e+00,  5.61851096e+00)
( 6.00000000e+00,  6.84297142e+00)
( 7.00000000e+00,  8.37102068e+00)
( 8.00000000e+00,  1.00889534e+01)
( 9.00000000e+00,  1.21829106e+01)};\label{line:lebcnst:0}

\addplot [black, thick, dashed, mark=*, mark size=2, mark options={solid}]
coordinates {
( 1.00000000e+00,  1.33333333e+00)
( 2.00000000e+00,  2.02666667e+00)
( 3.00000000e+00,  2.53900178e+00)
( 4.00000000e+00,  2.91409044e+00)
( 5.00000000e+00,  3.21272456e+00)
( 6.00000000e+00,  3.47212714e+00)
( 7.00000000e+00,  3.71757931e+00)
( 8.00000000e+00,  3.97009177e+00)
( 9.00000000e+00,  4.25131115e+00)};\label{line:lebobj:0}

\addplot [blue, thick, mark=+, mark size=2, mark options={solid}]
coordinates {
( 1.00000000e+00,  1.00000000e+00)
( 2.00000000e+00,  2.08325000e+00)
( 3.00000000e+00,  3.70185776e+00)
( 4.00000000e+00,  7.67153626e+00)
( 5.00000000e+00,  1.69277451e+01)
( 6.00000000e+00,  3.97593359e+01)
( 7.00000000e+00,  9.93337816e+01)
( 8.00000000e+00,  2.62374909e+02)
( 9.00000000e+00,  7.28122773e+02)};\label{line:lebcnst:1}

\addplot [blue, thick, dashed, mark=+, mark size=2, mark options={solid}]
coordinates {
( 1.00000000e+00,  1.33333333e+00)
( 2.00000000e+00,  2.02666667e+00)
( 3.00000000e+00,  2.97598639e+00)
( 4.00000000e+00,  4.35796684e+00)
( 5.00000000e+00,  6.70966444e+00)
( 6.00000000e+00,  1.15687190e+01)
( 7.00000000e+00,  2.40574519e+01)
( 8.00000000e+00,  6.45730247e+01)
( 9.00000000e+00,  2.30245959e+02)};\label{line:lebobj:1}

\addplot [red, thick, mark=triangle, mark size=2, mark options={solid}]
coordinates {
( 1.00000000e+00,  1.00000000e+00)
( 2.00000000e+00,  2.08325000e+00)
( 3.00000000e+00,  3.16880081e+00)
( 4.00000000e+00,  4.38137443e+00)
( 5.00000000e+00,  6.05876178e+00)
( 6.00000000e+00,  7.31305305e+00)
( 7.00000000e+00,  8.83403306e+00)
( 8.00000000e+00,  1.04348562e+01)
( 9.00000000e+00,  1.24555840e+01)};\label{line:lebcnst:2}

\addplot [red, thick, dashed, mark=triangle, mark size=2, mark options={solid}]
coordinates {
( 1.00000000e+00,  1.33333333e+00)
( 2.00000000e+00,  2.02666667e+00)
( 3.00000000e+00,  2.53877551e+00)
( 4.00000000e+00,  2.91509693e+00)
( 5.00000000e+00,  3.21939315e+00)
( 6.00000000e+00,  3.48809379e+00)
( 7.00000000e+00,  3.74528735e+00)
( 8.00000000e+00,  4.01108768e+00)
( 9.00000000e+00,  4.30657345e+00)};\label{line:lebobj:2}

\addplot [magenta, thick, mark=x, mark size=2, mark options={solid}]
coordinates {
( 1.00000000e+00,  1.00000000e+00)
( 2.00000000e+00,  2.08325000e+00)
( 3.00000000e+00,  3.16880081e+00)
( 4.00000000e+00,  4.35457597e+00)
( 5.00000000e+00,  5.55085448e+00)
( 6.00000000e+00,  6.93499682e+00)
( 7.00000000e+00,  8.42858895e+00)
( 8.00000000e+00,  1.01419696e+01)
( 9.00000000e+00,  1.21605523e+01)};\label{line:lebcnst:3}

\addplot [magenta, thick, dashed, mark=x, mark size=2, mark options={solid}]
coordinates {
( 1.00000000e+00,  1.33333333e+00)
( 2.00000000e+00,  2.02666667e+00)
( 3.00000000e+00,  2.53877551e+00)
( 4.00000000e+00,  2.91792453e+00)
( 5.00000000e+00,  3.22531438e+00)
( 6.00000000e+00,  3.48504184e+00)
( 7.00000000e+00,  3.73705234e+00)
( 8.00000000e+00,  3.99416767e+00)
( 9.00000000e+00,  4.28322651e+00)};\label{line:lebobj:3}

\addplot [green, thick, mark=diamond, mark size=2, mark options={solid}]
coordinates {
( 1.00000000e+00,  1.00000000e+00)
( 2.00000000e+00,  2.08325000e+00)
( 3.00000000e+00,  3.08620784e+00)
( 4.00000000e+00,  4.27186954e+00)
( 5.00000000e+00,  5.61851096e+00)
( 6.00000000e+00,  6.84297142e+00)
( 7.00000000e+00,  8.37102068e+00)
( 8.00000000e+00,  1.00889534e+01)
( 9.00000000e+00,  1.21829300e+01)};\label{line:lebcnst:4}

\addplot [green, thick, dashed, mark=diamond, mark size=2, mark options={solid}]
coordinates {
( 1.00000000e+00,  1.33333333e+00)
( 2.00000000e+00,  2.02666667e+00)
( 3.00000000e+00,  2.53900178e+00)
( 4.00000000e+00,  2.91409044e+00)
( 5.00000000e+00,  3.21272456e+00)
( 6.00000000e+00,  3.47212714e+00)
( 7.00000000e+00,  3.71757931e+00)
( 8.00000000e+00,  3.97009177e+00)
( 9.00000000e+00,  4.25131116e+00)};\label{line:lebobj:4}

\nextgroupplot[xmin=1, xmax=9, ymode=log, xlabel={Polynomial degree}, ylabel={Mass Matrix Condition Number}]
\addplot [black, thick, mark=*, mark size=2, mark options={solid}]
coordinates {
( 1.00000000e+00,  1.20000000e+01)
( 2.00000000e+00,  1.18388733e+02)
( 3.00000000e+00,  2.89863818e+02)
( 4.00000000e+00,  4.81298565e+02)
( 5.00000000e+00,  7.60266943e+02)
( 6.00000000e+00,  1.21271540e+03)
( 7.00000000e+00,  1.80373508e+03)
( 8.00000000e+00,  3.06762087e+03)
( 9.00000000e+00,  4.93192128e+03)};\label{line:Mcond:0}

\addplot [blue, thick, mark=+, mark size=2, mark options={solid}]
coordinates {
( 1.00000000e+00,  1.20000000e+01)
( 2.00000000e+00,  1.18388733e+02)
( 3.00000000e+00,  3.17985857e+02)
( 4.00000000e+00,  8.16447094e+02)
( 5.00000000e+00,  2.44443559e+03)
( 6.00000000e+00,  9.49986515e+03)
( 7.00000000e+00,  4.60124907e+04)
( 8.00000000e+00,  2.78814314e+05)
( 9.00000000e+00,  2.06087144e+06)};\label{line:Mcond:1}

\addplot [red, thick, mark=triangle, mark size=2, mark options={solid}]
coordinates {
( 1.00000000e+00,  1.20000000e+01)
( 2.00000000e+00,  1.18388733e+02)
( 3.00000000e+00,  3.01448155e+02)
( 4.00000000e+00,  4.89386281e+02)
( 5.00000000e+00,  8.28577495e+02)
( 6.00000000e+00,  1.25420213e+03)
( 7.00000000e+00,  1.93026005e+03)
( 8.00000000e+00,  3.26744732e+03)
( 9.00000000e+00,  5.30515609e+03)};\label{line:Mcond:2}

\addplot [magenta, thick, mark=x, mark size=2, mark options={solid}]
coordinates {
( 1.00000000e+00,  1.20000000e+01)
( 2.00000000e+00,  1.18388733e+02)
( 3.00000000e+00,  3.01448155e+02)
( 4.00000000e+00,  4.78189533e+02)
( 5.00000000e+00,  7.37902045e+02)
( 6.00000000e+00,  1.25170809e+03)
( 7.00000000e+00,  1.88851654e+03)
( 8.00000000e+00,  3.22938697e+03)
( 9.00000000e+00,  5.22542426e+03)};\label{line:Mcond:3}

\addplot [green, thick, mark=diamond, mark size=2, mark options={solid}]
coordinates {
( 1.00000000e+00,  1.20000000e+01)
( 2.00000000e+00,  1.18388733e+02)
( 3.00000000e+00,  2.89863818e+02)
( 4.00000000e+00,  4.81298565e+02)
( 5.00000000e+00,  7.60266943e+02)
( 6.00000000e+00,  1.21271540e+03)
( 7.00000000e+00,  1.80373508e+03)
( 8.00000000e+00,  3.06762088e+03)
( 9.00000000e+00,  4.93192866e+03)};\label{line:Mcond:4}

\end{groupplot}\end{tikzpicture}

%% file: _py/compare_triprism_zoom.tikz
\begin{tikzpicture}
\begin{groupplot} [
group style={group size = 2 by 1, horizontal sep = 2.0cm}]
\nextgroupplot[xmin=1, xmax=9, ymode=linear, xlabel={Polynomial degree}, ylabel={Lebesgue Constant}, ymin=1, ymax=8]
\addplot [black, thick, mark=*, mark size=2, mark options={solid}]
coordinates {
( 1.00000000e+00,  1.00000000e+00)
( 2.00000000e+00,  2.08325000e+00)
( 3.00000000e+00,  3.08620784e+00)
( 4.00000000e+00,  4.27186954e+00)
( 5.00000000e+00,  5.61851096e+00)
( 6.00000000e+00,  6.84297142e+00)
( 7.00000000e+00,  8.37102068e+00)};\label{line:lebcnst:0}

\addplot [black, thick, dashed, mark=*, mark size=2, mark options={solid}]
coordinates {
( 1.00000000e+00,  1.33333333e+00)
( 2.00000000e+00,  2.02666667e+00)
( 3.00000000e+00,  2.53900178e+00)
( 4.00000000e+00,  2.91409044e+00)
( 5.00000000e+00,  3.21272456e+00)
( 6.00000000e+00,  3.47212714e+00)
( 7.00000000e+00,  3.71757931e+00)
( 8.00000000e+00,  3.97009177e+00)
( 9.00000000e+00,  4.25131115e+00)};\label{line:lebobj:0}

\addplot [blue, thick, mark=+, mark size=2, mark options={solid}]
coordinates {
( 1.00000000e+00,  1.00000000e+00)
( 2.00000000e+00,  2.08325000e+00)
( 3.00000000e+00,  3.70185776e+00)
( 4.00000000e+00,  7.67153626e+00)
( 5.00000000e+00,  1.69277451e+01)};\label{line:lebcnst:1}

\addplot [blue, thick, dashed, mark=+, mark size=2, mark options={solid}]
coordinates {
( 1.00000000e+00,  1.33333333e+00)
( 2.00000000e+00,  2.02666667e+00)
( 3.00000000e+00,  2.97598639e+00)
( 4.00000000e+00,  4.35796684e+00)
( 5.00000000e+00,  6.70966444e+00)
( 6.00000000e+00,  1.15687190e+01)};\label{line:lebobj:1}

\addplot [red, thick, mark=triangle, mark size=2, mark options={solid}]
coordinates {
( 1.00000000e+00,  1.00000000e+00)
( 2.00000000e+00,  2.08325000e+00)
( 3.00000000e+00,  3.16880081e+00)
( 4.00000000e+00,  4.38137443e+00)
( 5.00000000e+00,  6.05876178e+00)
( 6.00000000e+00,  7.31305305e+00)
( 7.00000000e+00,  8.83403306e+00)};\label{line:lebcnst:2}

\addplot [red, thick, dashed, mark=triangle, mark size=2, mark options={solid}]
coordinates {
( 1.00000000e+00,  1.33333333e+00)
( 2.00000000e+00,  2.02666667e+00)
( 3.00000000e+00,  2.53877551e+00)
( 4.00000000e+00,  2.91509693e+00)
( 5.00000000e+00,  3.21939315e+00)
( 6.00000000e+00,  3.48809379e+00)
( 7.00000000e+00,  3.74528735e+00)
( 8.00000000e+00,  4.01108768e+00)
( 9.00000000e+00,  4.30657345e+00)};\label{line:lebobj:2}

\addplot [magenta, thick, mark=x, mark size=2, mark options={solid}]
coordinates {
( 1.00000000e+00,  1.00000000e+00)
( 2.00000000e+00,  2.08325000e+00)
( 3.00000000e+00,  3.16880081e+00)
( 4.00000000e+00,  4.35457597e+00)
( 5.00000000e+00,  5.55085448e+00)
( 6.00000000e+00,  6.93499682e+00)
( 7.00000000e+00,  8.42858895e+00)};\label{line:lebcnst:3}

\addplot [magenta, thick, dashed, mark=x, mark size=2, mark options={solid}]
coordinates {
( 1.00000000e+00,  1.33333333e+00)
( 2.00000000e+00,  2.02666667e+00)
( 3.00000000e+00,  2.53877551e+00)
( 4.00000000e+00,  2.91792453e+00)
( 5.00000000e+00,  3.22531438e+00)
( 6.00000000e+00,  3.48504184e+00)
( 7.00000000e+00,  3.73705234e+00)
( 8.00000000e+00,  3.99416767e+00)
( 9.00000000e+00,  4.28322651e+00)};\label{line:lebobj:3}

\addplot [green, thick, mark=diamond, mark size=2, mark options={solid}]
coordinates {
( 1.00000000e+00,  1.00000000e+00)
( 2.00000000e+00,  2.08325000e+00)
( 3.00000000e+00,  3.08620784e+00)
( 4.00000000e+00,  4.27186954e+00)
( 5.00000000e+00,  5.61851096e+00)
( 6.00000000e+00,  6.84297142e+00)
( 7.00000000e+00,  8.37102068e+00)};\label{line:lebcnst:4}

\addplot [green, thick, dashed, mark=diamond, mark size=2, mark options={solid}]
coordinates {
( 1.00000000e+00,  1.33333333e+00)
( 2.00000000e+00,  2.02666667e+00)
( 3.00000000e+00,  2.53900178e+00)
( 4.00000000e+00,  2.91409044e+00)
( 5.00000000e+00,  3.21272456e+00)
( 6.00000000e+00,  3.47212714e+00)
( 7.00000000e+00,  3.71757931e+00)
( 8.00000000e+00,  3.97009177e+00)
( 9.00000000e+00,  4.25131116e+00)};\label{line:lebobj:4}

\nextgroupplot[xmin=1, xmax=9, ymode=linear, xlabel={Polynomial degree}, ylabel={Mass Matrix Condition Number}, ymin=1, ymax=1500]
\addplot [black, thick, mark=*, mark size=2, mark options={solid}]
coordinates {
( 1.00000000e+00,  1.20000000e+01)
( 2.00000000e+00,  1.18388733e+02)
( 3.00000000e+00,  2.89863818e+02)
( 4.00000000e+00,  4.81298565e+02)
( 5.00000000e+00,  7.60266943e+02)
( 6.00000000e+00,  1.21271540e+03)
( 7.00000000e+00,  1.80373508e+03)};\label{line:Mcond:0}

\addplot [blue, thick, mark=+, mark size=2, mark options={solid}]
coordinates {
( 1.00000000e+00,  1.20000000e+01)
( 2.00000000e+00,  1.18388733e+02)
( 3.00000000e+00,  3.17985857e+02)
( 4.00000000e+00,  8.16447094e+02)
( 5.00000000e+00,  2.44443559e+03)};\label{line:Mcond:1}

\addplot [red, thick, mark=triangle, mark size=2, mark options={solid}]
coordinates {
( 1.00000000e+00,  1.20000000e+01)
( 2.00000000e+00,  1.18388733e+02)
( 3.00000000e+00,  3.01448155e+02)
( 4.00000000e+00,  4.89386281e+02)
( 5.00000000e+00,  8.28577495e+02)
( 6.00000000e+00,  1.25420213e+03)
( 7.00000000e+00,  1.93026005e+03)};\label{line:Mcond:2}

\addplot [magenta, thick, mark=x, mark size=2, mark options={solid}]
coordinates {
( 1.00000000e+00,  1.20000000e+01)
( 2.00000000e+00,  1.18388733e+02)
( 3.00000000e+00,  3.01448155e+02)
( 4.00000000e+00,  4.78189533e+02)
( 5.00000000e+00,  7.37902045e+02)
( 6.00000000e+00,  1.25170809e+03)
( 7.00000000e+00,  1.88851654e+03)};\label{line:Mcond:3}

\addplot [green, thick, mark=diamond, mark size=2, mark options={solid}]
coordinates {
( 1.00000000e+00,  1.20000000e+01)
( 2.00000000e+00,  1.18388733e+02)
( 3.00000000e+00,  2.89863818e+02)
( 4.00000000e+00,  4.81298565e+02)
( 5.00000000e+00,  7.60266943e+02)
( 6.00000000e+00,  1.21271540e+03)
( 7.00000000e+00,  1.80373508e+03)};\label{line:Mcond:4}

\end{groupplot}\end{tikzpicture}

%% file: _py/nodes_triprism.tikz
\begin{tikzpicture}
\begin{axis}[
axis equal image,
xmin=-1,
xmax=1,
ymin=-1,
ymax=1,
zmin=-1,
zmax=1,
axis lines=none,
width=0.6\textwidth]
\addplot3 [black, thick, forget plot]
coordinates {
(-1.00000000e+00, -1.00000000e+00, -1.00000000e+00)
( 1.00000000e+00, -1.00000000e+00, -1.00000000e+00)
(-1.00000000e+00,  1.00000000e+00, -1.00000000e+00)
(-1.00000000e+00, -1.00000000e+00, -1.00000000e+00)};

\addplot3 [black, thick, forget plot]
coordinates {
(-1.00000000e+00, -1.00000000e+00,  1.00000000e+00)
( 1.00000000e+00, -1.00000000e+00,  1.00000000e+00)
(-1.00000000e+00,  1.00000000e+00,  1.00000000e+00)
(-1.00000000e+00, -1.00000000e+00,  1.00000000e+00)};

\addplot3 [black, thick, forget plot]
coordinates {
(-1.00000000e+00, -1.00000000e+00, -1.00000000e+00)
(-1.00000000e+00, -1.00000000e+00,  1.00000000e+00)};

\addplot3 [black, thick, forget plot]
coordinates {
( 1.00000000e+00, -1.00000000e+00, -1.00000000e+00)
( 1.00000000e+00, -1.00000000e+00,  1.00000000e+00)};

\addplot3 [black, thick, forget plot]
coordinates {
(-1.00000000e+00,  1.00000000e+00, -1.00000000e+00)
(-1.00000000e+00,  1.00000000e+00,  1.00000000e+00)};

\addplot3 [black, thick, mark=*, mark size=2, mark options={solid}, only marks]
coordinates {
(-1.00000000e+00, -1.00000000e+00, -1.00000000e+00)
(-6.36326016e-01, -1.00000000e+00, -1.00000000e+00)
( 0.00000000e+00, -1.00000000e+00, -1.00000000e+00)
( 6.36326016e-01, -1.00000000e+00, -1.00000000e+00)
( 1.00000000e+00, -1.00000000e+00, -1.00000000e+00)
(-1.00000000e+00, -6.36326016e-01, -1.00000000e+00)
(-5.57815335e-01, -5.57815335e-01, -1.00000000e+00)
( 1.15630670e-01, -5.57815335e-01, -1.00000000e+00)
( 6.36326016e-01, -6.36326016e-01, -1.00000000e+00)
(-1.00000000e+00,  0.00000000e+00, -1.00000000e+00)
(-5.57815335e-01,  1.15630670e-01, -1.00000000e+00)
( 0.00000000e+00,  0.00000000e+00, -1.00000000e+00)
(-1.00000000e+00,  6.36326016e-01, -1.00000000e+00)
(-6.36326016e-01,  6.36326016e-01, -1.00000000e+00)
(-1.00000000e+00,  1.00000000e+00, -1.00000000e+00)
(-1.00000000e+00, -1.00000000e+00, -6.36326016e-01)
(-6.36326016e-01, -1.00000000e+00, -6.36326016e-01)
( 0.00000000e+00, -1.00000000e+00, -6.36326016e-01)
( 6.36326016e-01, -1.00000000e+00, -6.36326016e-01)
( 1.00000000e+00, -1.00000000e+00, -6.36326016e-01)
(-1.00000000e+00, -6.36326016e-01, -6.36326016e-01)
(-5.57815335e-01, -5.57815335e-01, -6.36326016e-01)
( 1.15630670e-01, -5.57815335e-01, -6.36326016e-01)
( 6.36326016e-01, -6.36326016e-01, -6.36326016e-01)
(-1.00000000e+00,  0.00000000e+00, -6.36326016e-01)
(-5.57815335e-01,  1.15630670e-01, -6.36326016e-01)
( 0.00000000e+00,  0.00000000e+00, -6.36326016e-01)
(-1.00000000e+00,  6.36326016e-01, -6.36326016e-01)
(-6.36326016e-01,  6.36326016e-01, -6.36326016e-01)
(-1.00000000e+00,  1.00000000e+00, -6.36326016e-01)
(-1.00000000e+00, -1.00000000e+00,  0.00000000e+00)
(-6.36326016e-01, -1.00000000e+00,  0.00000000e+00)
( 0.00000000e+00, -1.00000000e+00,  0.00000000e+00)
( 6.36326016e-01, -1.00000000e+00,  0.00000000e+00)
( 1.00000000e+00, -1.00000000e+00,  0.00000000e+00)
(-1.00000000e+00, -6.36326016e-01,  0.00000000e+00)
(-5.57815335e-01, -5.57815335e-01,  0.00000000e+00)
( 1.15630670e-01, -5.57815335e-01,  0.00000000e+00)
( 6.36326016e-01, -6.36326016e-01,  0.00000000e+00)
(-1.00000000e+00,  0.00000000e+00,  0.00000000e+00)
(-5.57815335e-01,  1.15630670e-01,  0.00000000e+00)
( 0.00000000e+00,  0.00000000e+00,  0.00000000e+00)
(-1.00000000e+00,  6.36326016e-01,  0.00000000e+00)
(-6.36326016e-01,  6.36326016e-01,  0.00000000e+00)
(-1.00000000e+00,  1.00000000e+00,  0.00000000e+00)
(-1.00000000e+00, -1.00000000e+00,  6.36326016e-01)
(-6.36326016e-01, -1.00000000e+00,  6.36326016e-01)
( 0.00000000e+00, -1.00000000e+00,  6.36326016e-01)
( 6.36326016e-01, -1.00000000e+00,  6.36326016e-01)
( 1.00000000e+00, -1.00000000e+00,  6.36326016e-01)
(-1.00000000e+00, -6.36326016e-01,  6.36326016e-01)
(-5.57815335e-01, -5.57815335e-01,  6.36326016e-01)
( 1.15630670e-01, -5.57815335e-01,  6.36326016e-01)
( 6.36326016e-01, -6.36326016e-01,  6.36326016e-01)
(-1.00000000e+00,  0.00000000e+00,  6.36326016e-01)
(-5.57815335e-01,  1.15630670e-01,  6.36326016e-01)
( 0.00000000e+00,  0.00000000e+00,  6.36326016e-01)
(-1.00000000e+00,  6.36326016e-01,  6.36326016e-01)
(-6.36326016e-01,  6.36326016e-01,  6.36326016e-01)
(-1.00000000e+00,  1.00000000e+00,  6.36326016e-01)
(-1.00000000e+00, -1.00000000e+00,  1.00000000e+00)
(-6.36326016e-01, -1.00000000e+00,  1.00000000e+00)
( 0.00000000e+00, -1.00000000e+00,  1.00000000e+00)
( 6.36326016e-01, -1.00000000e+00,  1.00000000e+00)
( 1.00000000e+00, -1.00000000e+00,  1.00000000e+00)
(-1.00000000e+00, -6.36326016e-01,  1.00000000e+00)
(-5.57815335e-01, -5.57815335e-01,  1.00000000e+00)
( 1.15630670e-01, -5.57815335e-01,  1.00000000e+00)
( 6.36326016e-01, -6.36326016e-01,  1.00000000e+00)
(-1.00000000e+00,  0.00000000e+00,  1.00000000e+00)
(-5.57815335e-01,  1.15630670e-01,  1.00000000e+00)
( 0.00000000e+00,  0.00000000e+00,  1.00000000e+00)
(-1.00000000e+00,  6.36326016e-01,  1.00000000e+00)
(-6.36326016e-01,  6.36326016e-01,  1.00000000e+00)
(-1.00000000e+00,  1.00000000e+00,  1.00000000e+00)};\label{line:nodes:0}

\end{axis}
\end{tikzpicture}

%% file: _py/compare_quadpyrmd.tikz
\begin{tikzpicture}
\begin{groupplot} [
group style={group size = 2 by 1, horizontal sep = 2.0cm}]
\nextgroupplot[xmin=1, xmax=9, ymode=log, xlabel={Polynomial degree}, ylabel={Lebesgue Constant}]
\addplot [black, thick, mark=*, mark size=2, mark options={solid}]
coordinates {
( 1.00000000e+00,  1.00000000e+00)
( 2.00000000e+00,  1.84205915e+00)
( 3.00000000e+00,  2.75550443e+00)
( 4.00000000e+00,  3.85985041e+00)
( 5.00000000e+00,  4.94776905e+00)
( 6.00000000e+00,  6.43828204e+00)
( 7.00000000e+00,  8.06288845e+00)
( 8.00000000e+00,  1.05474265e+01)
( 9.00000000e+00,  1.34543576e+01)};\label{line:lebcnst:0}

\addplot [black, thick, dashed, mark=*, mark size=2, mark options={solid}]
coordinates {
( 1.00000000e+00,  9.77777778e-01)
( 2.00000000e+00,  1.30285714e+00)
( 3.00000000e+00,  1.65658611e+00)
( 4.00000000e+00,  1.93608723e+00)
( 5.00000000e+00,  2.16701832e+00)
( 6.00000000e+00,  2.37474740e+00)
( 7.00000000e+00,  2.57798590e+00)
( 8.00000000e+00,  2.79301530e+00)
( 9.00000000e+00,  3.03730876e+00)};\label{line:lebobj:0}

\addplot [blue, thick, mark=+, mark size=2, mark options={solid}]
coordinates {
( 1.00000000e+00,  1.00000000e+00)
( 2.00000000e+00,  1.84205915e+00)
( 3.00000000e+00,  3.15377004e+00)
( 4.00000000e+00,  5.94407419e+00)
( 5.00000000e+00,  1.18634798e+01)
( 6.00000000e+00,  2.51017471e+01)
( 7.00000000e+00,  5.66299432e+01)
( 8.00000000e+00,  1.35956762e+02)
( 9.00000000e+00,  3.49739035e+02)};\label{line:lebcnst:1}

\addplot [blue, thick, dashed, mark=+, mark size=2, mark options={solid}]
coordinates {
( 1.00000000e+00,  9.77777778e-01)
( 2.00000000e+00,  1.30285714e+00)
( 3.00000000e+00,  1.83381519e+00)
( 4.00000000e+00,  2.56254232e+00)
( 5.00000000e+00,  3.67333545e+00)
( 6.00000000e+00,  5.62826608e+00)
( 7.00000000e+00,  9.67908743e+00)
( 8.00000000e+00,  1.97399289e+01)
( 9.00000000e+00,  5.04048669e+01)};\label{line:lebobj:1}

\addplot [red, thick, mark=triangle, mark size=2, mark options={solid}]
coordinates {
( 1.00000000e+00,  1.00000000e+00)
( 2.00000000e+00,  1.84205915e+00)
( 3.00000000e+00,  2.82585341e+00)
( 4.00000000e+00,  4.28988817e+00)
( 5.00000000e+00,  6.84022547e+00)
( 6.00000000e+00,  1.01028275e+01)
( 7.00000000e+00,  1.41923068e+01)
( 8.00000000e+00,  2.04189792e+01)
( 9.00000000e+00,  3.11271126e+01)};\label{line:lebcnst:2}

\addplot [red, thick, dashed, mark=triangle, mark size=2, mark options={solid}]
coordinates {
( 1.00000000e+00,  9.77777778e-01)
( 2.00000000e+00,  1.30285714e+00)
( 3.00000000e+00,  1.69075586e+00)
( 4.00000000e+00,  2.04964027e+00)
( 5.00000000e+00,  2.41420571e+00)
( 6.00000000e+00,  2.83488062e+00)
( 7.00000000e+00,  3.39434728e+00)
( 8.00000000e+00,  4.24709522e+00)
( 9.00000000e+00,  5.71124595e+00)};\label{line:lebobj:2}

\addplot [magenta, thick, mark=x, mark size=2, mark options={solid}]
coordinates {
( 1.00000000e+00,  1.00000000e+00)
( 2.00000000e+00,  2.08000000e+00)
( 3.00000000e+00,  2.73000000e+00)
( 4.00000000e+00,  4.13000000e+00)
( 5.00000000e+00,  5.53000000e+00)
( 6.00000000e+00,  7.35000000e+00)
( 7.00000000e+00,  9.71000000e+00)
( 8.00000000e+00,  1.28000000e+01)
( 9.00000000e+00,  1.72000000e+01)};\label{line:lebcnst:3}

\addplot [magenta, thick, dashed, mark=x, mark size=2, mark options={solid}]
coordinates {
( 1.00000000e+00,          nan)
( 2.00000000e+00,          nan)
( 3.00000000e+00,          nan)
( 4.00000000e+00,          nan)
( 5.00000000e+00,          nan)
( 6.00000000e+00,          nan)
( 7.00000000e+00,          nan)
( 8.00000000e+00,          nan)
( 9.00000000e+00,          nan)};\label{line:lebobj:3}

\nextgroupplot[xmin=1, xmax=9, ymode=log, xlabel={Polynomial degree}, ylabel={Mass Matrix Condition Number}]
\addplot [black, thick, mark=*, mark size=2, mark options={solid}]
coordinates {
( 1.00000000e+00,  1.22434165e+01)
( 2.00000000e+00,  7.70350298e+01)
( 3.00000000e+00,  2.58375834e+02)
( 4.00000000e+00,  3.83866308e+02)
( 5.00000000e+00,  6.12389895e+02)
( 6.00000000e+00,  1.08178880e+03)
( 7.00000000e+00,  1.98734696e+03)
( 8.00000000e+00,  3.75294018e+03)
( 9.00000000e+00,  7.18623438e+03)};\label{line:Mcond:0}

\addplot [blue, thick, mark=+, mark size=2, mark options={solid}]
coordinates {
( 1.00000000e+00,  1.22434165e+01)
( 2.00000000e+00,  7.70350298e+01)
( 3.00000000e+00,  2.51003001e+02)
( 4.00000000e+00,  4.90405701e+02)
( 5.00000000e+00,  1.21024544e+03)
( 6.00000000e+00,  3.70102282e+03)
( 7.00000000e+00,  1.52420665e+04)
( 8.00000000e+00,  9.09933108e+04)
( 9.00000000e+00,  6.56205206e+05)};\label{line:Mcond:1}

\addplot [red, thick, mark=triangle, mark size=2, mark options={solid}]
coordinates {
( 1.00000000e+00,  1.22434165e+01)
( 2.00000000e+00,  7.70350298e+01)
( 3.00000000e+00,  2.69873493e+02)
( 4.00000000e+00,  4.23150154e+02)
( 5.00000000e+00,  8.81235463e+02)
( 6.00000000e+00,  1.46405400e+03)
( 7.00000000e+00,  2.81112232e+03)
( 8.00000000e+00,  6.53735824e+03)
( 9.00000000e+00,  1.71944456e+04)};\label{line:Mcond:2}

\addplot [magenta, thick, mark=x, mark size=2, mark options={solid}]
coordinates {
( 1.00000000e+00,          nan)
( 2.00000000e+00,          nan)
( 3.00000000e+00,          nan)
( 4.00000000e+00,          nan)
( 5.00000000e+00,          nan)
( 6.00000000e+00,          nan)
( 7.00000000e+00,          nan)
( 8.00000000e+00,          nan)
( 9.00000000e+00,          nan)};\label{line:Mcond:3}

\end{groupplot}\end{tikzpicture}

%% file: _py/compare_quadpyrmd_zoom.tikz
\begin{tikzpicture}
\begin{groupplot} [
group style={group size = 2 by 1, horizontal sep = 2.0cm}]
\nextgroupplot[xmin=1, xmax=9, ymode=linear, xlabel={Polynomial degree}, ylabel={Lebesgue Constant}, ymin=1, ymax=8]
\addplot [black, thick, mark=*, mark size=2, mark options={solid}]
coordinates {
( 1.00000000e+00,  1.00000000e+00)
( 2.00000000e+00,  1.84205915e+00)
( 3.00000000e+00,  2.75550443e+00)
( 4.00000000e+00,  3.85985041e+00)
( 5.00000000e+00,  4.94776905e+00)
( 6.00000000e+00,  6.43828204e+00)
( 7.00000000e+00,  8.06288845e+00)};\label{line:lebcnst:0}

\addplot [black, thick, dashed, mark=*, mark size=2, mark options={solid}]
coordinates {
( 1.00000000e+00,  9.77777778e-01)
( 2.00000000e+00,  1.30285714e+00)
( 3.00000000e+00,  1.65658611e+00)
( 4.00000000e+00,  1.93608723e+00)
( 5.00000000e+00,  2.16701832e+00)
( 6.00000000e+00,  2.37474740e+00)
( 7.00000000e+00,  2.57798590e+00)
( 8.00000000e+00,  2.79301530e+00)
( 9.00000000e+00,  3.03730876e+00)};\label{line:lebobj:0}

\addplot [blue, thick, mark=+, mark size=2, mark options={solid}]
coordinates {
( 1.00000000e+00,  1.00000000e+00)
( 2.00000000e+00,  1.84205915e+00)
( 3.00000000e+00,  3.15377004e+00)
( 4.00000000e+00,  5.94407419e+00)
( 5.00000000e+00,  1.18634798e+01)};\label{line:lebcnst:1}

\addplot [blue, thick, dashed, mark=+, mark size=2, mark options={solid}]
coordinates {
( 1.00000000e+00,  9.77777778e-01)
( 2.00000000e+00,  1.30285714e+00)
( 3.00000000e+00,  1.83381519e+00)
( 4.00000000e+00,  2.56254232e+00)
( 5.00000000e+00,  3.67333545e+00)
( 6.00000000e+00,  5.62826608e+00)
( 7.00000000e+00,  9.67908743e+00)};\label{line:lebobj:1}

\addplot [red, thick, mark=triangle, mark size=2, mark options={solid}]
coordinates {
( 1.00000000e+00,  1.00000000e+00)
( 2.00000000e+00,  1.84205915e+00)
( 3.00000000e+00,  2.82585341e+00)
( 4.00000000e+00,  4.28988817e+00)
( 5.00000000e+00,  6.84022547e+00)
( 6.00000000e+00,  1.01028275e+01)};\label{line:lebcnst:2}

\addplot [red, thick, dashed, mark=triangle, mark size=2, mark options={solid}]
coordinates {
( 1.00000000e+00,  9.77777778e-01)
( 2.00000000e+00,  1.30285714e+00)
( 3.00000000e+00,  1.69075586e+00)
( 4.00000000e+00,  2.04964027e+00)
( 5.00000000e+00,  2.41420571e+00)
( 6.00000000e+00,  2.83488062e+00)
( 7.00000000e+00,  3.39434728e+00)
( 8.00000000e+00,  4.24709522e+00)
( 9.00000000e+00,  5.71124595e+00)};\label{line:lebobj:2}

\addplot [magenta, thick, mark=x, mark size=2, mark options={solid}]
coordinates {
( 1.00000000e+00,  1.00000000e+00)
( 2.00000000e+00,  2.08000000e+00)
( 3.00000000e+00,  2.73000000e+00)
( 4.00000000e+00,  4.13000000e+00)
( 5.00000000e+00,  5.53000000e+00)
( 6.00000000e+00,  7.35000000e+00)
( 7.00000000e+00,  9.71000000e+00)};\label{line:lebcnst:3}

\addplot [magenta, thick, dashed, mark=x, mark size=2, mark options={solid}]
coordinates {
( 1.00000000e+00,          nan)
( 2.00000000e+00,          nan)
( 3.00000000e+00,          nan)
( 4.00000000e+00,          nan)
( 5.00000000e+00,          nan)
( 6.00000000e+00,          nan)
( 7.00000000e+00,          nan)
( 8.00000000e+00,          nan)
( 9.00000000e+00,          nan)};\label{line:lebobj:3}

\nextgroupplot[xmin=1, xmax=9, ymode=linear, xlabel={Polynomial degree}, ylabel={Mass Matrix Condition Number}, ymin=1, ymax=1500]
\addplot [black, thick, mark=*, mark size=2, mark options={solid}]
coordinates {
( 1.00000000e+00,  1.22434165e+01)
( 2.00000000e+00,  7.70350298e+01)
( 3.00000000e+00,  2.58375834e+02)
( 4.00000000e+00,  3.83866308e+02)
( 5.00000000e+00,  6.12389895e+02)
( 6.00000000e+00,  1.08178880e+03)
( 7.00000000e+00,  1.98734696e+03)};\label{line:Mcond:0}

\addplot [blue, thick, mark=+, mark size=2, mark options={solid}]
coordinates {
( 1.00000000e+00,  1.22434165e+01)
( 2.00000000e+00,  7.70350298e+01)
( 3.00000000e+00,  2.51003001e+02)
( 4.00000000e+00,  4.90405701e+02)
( 5.00000000e+00,  1.21024544e+03)
( 6.00000000e+00,  3.70102282e+03)};\label{line:Mcond:1}

\addplot [red, thick, mark=triangle, mark size=2, mark options={solid}]
coordinates {
( 1.00000000e+00,  1.22434165e+01)
( 2.00000000e+00,  7.70350298e+01)
( 3.00000000e+00,  2.69873493e+02)
( 4.00000000e+00,  4.23150154e+02)
( 5.00000000e+00,  8.81235463e+02)
( 6.00000000e+00,  1.46405400e+03)
( 7.00000000e+00,  2.81112232e+03)};\label{line:Mcond:2}

\addplot [magenta, thick, mark=x, mark size=2, mark options={solid}]
coordinates {
( 1.00000000e+00,          nan)
( 2.00000000e+00,          nan)
( 3.00000000e+00,          nan)
( 4.00000000e+00,          nan)
( 5.00000000e+00,          nan)
( 6.00000000e+00,          nan)
( 7.00000000e+00,          nan)
( 8.00000000e+00,          nan)
( 9.00000000e+00,          nan)};\label{line:Mcond:3}

\end{groupplot}\end{tikzpicture}

%% file: _py/nodes_quadpyrmd.tikz
\begin{tikzpicture}
\begin{axis}[
axis equal image,
xmin=-1,
xmax=1,
ymin=-1,
ymax=1,
zmin=-1,
zmax=1,
axis lines=none,
width=0.6\textwidth]
\addplot3 [black, thick, forget plot]
coordinates {
(-1.00000000e+00, -1.00000000e+00, -1.00000000e+00)
( 1.00000000e+00, -1.00000000e+00, -1.00000000e+00)
( 1.00000000e+00,  1.00000000e+00, -1.00000000e+00)
(-1.00000000e+00,  1.00000000e+00, -1.00000000e+00)
(-1.00000000e+00, -1.00000000e+00, -1.00000000e+00)};

\addplot3 [black, thick, forget plot]
coordinates {
(-1.00000000e+00, -1.00000000e+00, -1.00000000e+00)
( 0.00000000e+00,  0.00000000e+00,  1.00000000e+00)};

\addplot3 [black, thick, forget plot]
coordinates {
( 1.00000000e+00, -1.00000000e+00, -1.00000000e+00)
( 0.00000000e+00,  0.00000000e+00,  1.00000000e+00)};

\addplot3 [black, thick, forget plot]
coordinates {
( 1.00000000e+00,  1.00000000e+00, -1.00000000e+00)
( 0.00000000e+00,  0.00000000e+00,  1.00000000e+00)};

\addplot3 [black, thick, forget plot]
coordinates {
(-1.00000000e+00,  1.00000000e+00, -1.00000000e+00)
( 0.00000000e+00,  0.00000000e+00,  1.00000000e+00)};

\addplot3 [black, thick, mark=*, mark size=2, mark options={solid}, only marks]
coordinates {
(-1.00000000e+00, -1.00000000e+00, -1.00000000e+00)
(-6.36326016e-01, -1.00000000e+00, -1.00000000e+00)
( 0.00000000e+00, -1.00000000e+00, -1.00000000e+00)
( 6.36326016e-01, -1.00000000e+00, -1.00000000e+00)
( 1.00000000e+00, -1.00000000e+00, -1.00000000e+00)
(-1.00000000e+00, -6.36326016e-01, -1.00000000e+00)
(-6.36326016e-01, -6.36326016e-01, -1.00000000e+00)
( 0.00000000e+00, -6.36326016e-01, -1.00000000e+00)
( 6.36326016e-01, -6.36326016e-01, -1.00000000e+00)
( 1.00000000e+00, -6.36326016e-01, -1.00000000e+00)
(-1.00000000e+00,  0.00000000e+00, -1.00000000e+00)
(-6.36326016e-01,  0.00000000e+00, -1.00000000e+00)
( 0.00000000e+00,  0.00000000e+00, -1.00000000e+00)
( 6.36326016e-01,  0.00000000e+00, -1.00000000e+00)
( 1.00000000e+00,  0.00000000e+00, -1.00000000e+00)
(-1.00000000e+00,  6.36326016e-01, -1.00000000e+00)
(-6.36326016e-01,  6.36326016e-01, -1.00000000e+00)
( 0.00000000e+00,  6.36326016e-01, -1.00000000e+00)
( 6.36326016e-01,  6.36326016e-01, -1.00000000e+00)
( 1.00000000e+00,  6.36326016e-01, -1.00000000e+00)
(-1.00000000e+00,  1.00000000e+00, -1.00000000e+00)
(-6.36326016e-01,  1.00000000e+00, -1.00000000e+00)
( 0.00000000e+00,  1.00000000e+00, -1.00000000e+00)
( 6.36326016e-01,  1.00000000e+00, -1.00000000e+00)
( 1.00000000e+00,  1.00000000e+00, -1.00000000e+00)
(-8.18163008e-01, -8.18163008e-01, -6.36326016e-01)
(-3.36723003e-01, -7.78907668e-01, -5.57815335e-01)
( 3.36723003e-01, -7.78907668e-01, -5.57815335e-01)
( 8.18163008e-01, -8.18163008e-01, -6.36326016e-01)
(-7.78907668e-01, -3.36723003e-01, -5.57815335e-01)
(-3.31518348e-01, -3.31518348e-01, -5.50921374e-01)
( 3.31518348e-01, -3.31518348e-01, -5.50921374e-01)
( 7.78907668e-01, -3.36723003e-01, -5.57815335e-01)
(-7.78907668e-01,  3.36723003e-01, -5.57815335e-01)
(-3.31518348e-01,  3.31518348e-01, -5.50921374e-01)
( 3.31518348e-01,  3.31518348e-01, -5.50921374e-01)
( 7.78907668e-01,  3.36723003e-01, -5.57815335e-01)
(-8.18163008e-01,  8.18163008e-01, -6.36326016e-01)
(-3.36723003e-01,  7.78907668e-01, -5.57815335e-01)
( 3.36723003e-01,  7.78907668e-01, -5.57815335e-01)
( 8.18163008e-01,  8.18163008e-01, -6.36326016e-01)
(-5.00000000e-01, -5.00000000e-01, -2.93873588e-38)
( 0.00000000e+00, -4.42184665e-01,  1.15630670e-01)
( 5.00000000e-01, -5.00000000e-01, -2.93873588e-38)
(-4.42184665e-01,  0.00000000e+00,  1.15630670e-01)
( 0.00000000e+00,  0.00000000e+00,  1.13554189e-01)
( 4.42184665e-01,  0.00000000e+00,  1.15630670e-01)
(-5.00000000e-01,  5.00000000e-01, -2.93873588e-38)
( 0.00000000e+00,  4.42184665e-01,  1.15630670e-01)
( 5.00000000e-01,  5.00000000e-01, -2.93873588e-38)
(-1.81836992e-01, -1.81836992e-01,  6.36326016e-01)
( 1.81836992e-01, -1.81836992e-01,  6.36326016e-01)
(-1.81836992e-01,  1.81836992e-01,  6.36326016e-01)
( 1.81836992e-01,  1.81836992e-01,  6.36326016e-01)
( 0.00000000e+00,  0.00000000e+00,  1.00000000e+00)};\label{line:nodes:0}

\end{axis}
\end{tikzpicture}

%% file: _conclude.tex
In this paper, we introduce a new numerical framework to generate optimal, symmetric distributions for general elements that satisfy general linear constraints. Symmetry is ensured by constructing the nodal distributions from the symmetry orbits of the element geometry. In the context of finite element discretizations, symmetric nodal distributions ensure the approximation properties are invariant with respect to rotations in the symmetry group. Additional linear constraints are imposed on the symmetry orbits to endow the distribution with additional properties. In this work, we leverage this flexibility to ensure cross-element compatibility, i.e., nodes of adjacent elements will be co-located, in a conforming mesh (homogeneous or mixed elements). We used the framework to generate Lebesgue constant-minimizing nodal distributions with cross-element compatibility. These distributions were shown to generate finite element mass matrices with small condition number even for large polynomial degrees.

Future work could leverage the proposed abstract framework to generate optimal nodal distributions for other element geometries and higher dimensional elements. Other interesting avenues of research include the investigation of alternate objective functions and constraints to endow the distributions with other desirable properties, e.g., summation-by-parts on non-tensor-product elements. Additionally, some constraints on the distributions could be relaxed to create more optimal distributions. For example, relaxation of the cross-compatibility constraint or the requirement that nodes lie on each face could lead to nodal distributions with much lower Lebesgue constants; however, enforcing $C^0$ continuity in a finite element discretization would be much more difficult with such distributions.

%% file: _appendix.tex
\begin{table}[H]
\caption{Lebesgue constant and mass matrix condition number, formatted as $(\mathrm{leb\_cnst}, \mathrm{mass\_cond})$, for each distribution for the line element.}
\label{tab:line}
\centering
\begin{tabular}{c | c | c | c }
\toprule
 \bf{Poly.} & \textbf{Optimized Dist.}  &  \textbf{Unif. Dist.} & \textbf{GLL Dist.} \\
\hline
\input{_tab/line.tab}
\end{tabular}
\end{table}

\begin{table}[H]
\caption{Lebesgue constant and mass matrix condition number, formatted as $(\mathrm{leb\_cnst}, \mathrm{mass\_cond})$, for each distribution for the quadrilateral element.}
\label{tab:quad}
\centering
\begin{tabular}{c | c | c | c }
\toprule
 \bf{Poly.} & \textbf{Optimized Dist.} &  \textbf{Unif. Dist.} & \textbf{GLL Dist.} \\
\hline
\input{_tab/quad.tab}
\end{tabular}
\end{table}

\begin{table}[H]
\caption{Lebesgue constant and mass matrix condition number, formatted as $(\mathrm{leb\_cnst}, \mathrm{mass\_cond})$, for each distribution for the triangle element.}
\label{tab:tri}
\centering
\begin{tabular}{c | c | c | c | c }
\toprule
 \bf{Poly.} & \textbf{Optimized Dist.} &  \textbf{Unif. Dist.} & \textbf{Isaac Dist.} & \textbf{Warburton Dist.}  \\
\hline
\input{_tab/tri.tab}
\end{tabular}
\end{table}

\begin{table}[H]
\caption{Lebesgue constant and mass matrix condition number, formatted as $(\mathrm{leb\_cnst}, \mathrm{mass\_cond})$, for each distribution for the hexahedron element.}
\label{tab:hex}
\centering
\begin{tabular}{c | c | c | c }
\toprule
 \bf{Poly.} & \textbf{Optimized Dist.} &  \textbf{Unif. Dist.} & \textbf{GLL Dist.} \\
\hline
\input{_tab/hex.tab}
\end{tabular}
\end{table}

\begin{table}[H]
\caption{Lebesgue constant and mass matrix condition number, formatted as $(\mathrm{leb\_cnst}, \mathrm{mass\_cond})$, for each distribution for the tetrahedron element.}
\label{tab:tet}
\centering
\begin{tabular}{c | c | c | c | c }
\toprule
 \bf{Poly.} & \textbf{Optimized Dist.} &  \textbf{Unif. Dist.} & \textbf{Isaac Dist.} & \textbf{Warburton Dist.}  \\
\hline
\input{_tab/tet.tab}
\end{tabular}
\end{table}

\begin{sidewaystable}
\caption{Lebesgue constant and mass matrix condition number, formatted as $(\mathrm{leb\_cnst}, \mathrm{mass\_cond})$, for each distribution for the triangular prism element.}
\label{tab:triprism}
\centering
\begin{tabular}{c | c | c | c | c | c}
\toprule
 \bf{Poly.} & \textbf{Optimized Dist.} &  \textbf{Unif. Dist.} & \textbf{Isaac Dist.} & \textbf{Warburton Dist.} &  \textbf{Opt. Tri. $\times$ Opt. Line} \\
\hline
\input{_tab/triprism.tab}
\end{tabular}
\end{sidewaystable}

\begin{table}[H]
\caption{Lebesgue constant and mass matrix condition number, formatted as $(\mathrm{leb\_cnst}, \mathrm{mass\_cond})$, for each distribution for the pyramid element.}
\label{tab:quadpyrmd}
\centering
\begin{tabular}{c | c | c | c | c}
\toprule
 \bf{Poly.} & \textbf{Optimized Dist.} & \textbf{Unif. Dist.} & \textbf{Stroud Dist.} & \textbf{Chan Dist.} \\
\hline
\input{_tab/quadpyrmd.tab}
\end{tabular}
\end{table}

%% file: _tab/line.tab
1 & (1.00e+00, 3.00e+00) & (1.00e+00, 3.00e+00) & (1.00e+00, 3.00e+00)\\
2 & (1.25e+00, 6.88e+00) & (1.25e+00, 6.88e+00) & (1.25e+00, 6.88e+00)\\
3 & (1.46e+00, 8.36e+00) & (1.63e+00, 9.36e+00) & (1.50e+00, 8.65e+00)\\
4 & (1.59e+00, 1.02e+01) & (2.21e+00, 1.41e+01) & (1.64e+00, 1.04e+01)\\
5 & (1.72e+00, 1.17e+01) & (3.11e+00, 2.36e+01) & (1.78e+00, 1.20e+01)\\
6 & (1.81e+00, 1.33e+01) & (4.55e+00, 4.42e+01) & (1.87e+00, 1.36e+01)\\
7 & (1.91e+00, 1.49e+01) & (6.93e+00, 9.30e+01) & (1.97e+00, 1.52e+01)\\
8 & (1.98e+00, 1.65e+01) & (1.09e+01, 2.20e+02) & (2.05e+00, 1.67e+01)\\
9 & (2.05e+00, 1.81e+01) & (1.78e+01, 5.73e+02) & (2.12e+00, 1.83e+01)\\
10 & (2.11e+00, 1.97e+01) & (2.99e+01, 1.62e+03) & (2.18e+00, 1.99e+01)\\
11 & (2.17e+00, 2.12e+01) & (5.12e+01, 4.85e+03) & (2.24e+00, 2.14e+01)\\
12 & (2.22e+00, 2.28e+01) & (8.93e+01, 1.51e+04) & (2.29e+00, 2.30e+01)\\
13 & (2.27e+00, 2.44e+01) & (1.58e+02, 4.84e+04) & (2.34e+00, 2.46e+01)\\
14 & (2.31e+00, 2.60e+01) & (2.83e+02, 1.59e+05) & (2.39e+00, 2.61e+01)\\
15 & (2.35e+00, 2.75e+01) & (5.12e+02, 5.29e+05) & (2.43e+00, 2.77e+01)\\
16 & (2.39e+00, 2.91e+01) & (9.34e+02, 1.79e+06) & (2.47e+00, 2.93e+01)\\
17 & (2.43e+00, 3.07e+01) & (1.72e+03, 6.13e+06) & (2.51e+00, 3.08e+01)\\
18 & (2.46e+00, 3.23e+01) & (3.17e+03, 2.12e+07) & (2.54e+00, 3.24e+01)\\
19 & (2.49e+00, 3.38e+01) & (5.89e+03, 7.42e+07) & (2.58e+00, 3.40e+01)\\
20 & (2.53e+00, 3.54e+01) & (1.10e+04, 2.62e+08) & (2.61e+00, 3.55e+01)\\
21 & (2.56e+00, 3.70e+01) & (2.06e+04, 9.31e+08) & (2.64e+00, 3.71e+01)\\
22 & (2.58e+00, 3.86e+01) & (3.87e+04, 3.33e+09) & (2.67e+00, 3.87e+01)\\
23 & (2.61e+00, 4.02e+01) & (7.28e+04, 1.20e+10) & (2.69e+00, 4.02e+01)\\
24 & (2.64e+00, 4.17e+01) & (1.38e+05, 4.34e+10) & (2.72e+00, 4.18e+01)\\
25 & (2.66e+00, 4.33e+01) & (2.61e+05, 1.58e+11) & (2.75e+00, 4.34e+01)\\
26 & (2.69e+00, 4.49e+01) & (4.96e+05, 5.76e+11) & (2.77e+00, 4.49e+01)\\
27 & (2.71e+00, 4.65e+01) & (9.45e+05, 2.11e+12) & (2.79e+00, 4.65e+01)\\
28 & (2.73e+00, 4.80e+01) & (1.80e+06, 7.76e+12) & (2.82e+00, 4.81e+01)\\
29 & (2.75e+00, 4.96e+01) & (3.45e+06, 2.87e+13) & (2.84e+00, 4.96e+01)\\
30 & (2.77e+00, 5.12e+01) & (6.59e+06, 1.06e+14) & (2.86e+00, 5.12e+01)

%% file: _tab/quad.tab
1 & (1.00e+00, 9.00e+00) & (1.00e+00, 9.00e+00) & (1.00e+00, 9.00e+00)\\
2 & (1.56e+00, 4.73e+01) & (1.56e+00, 4.73e+01) & (1.56e+00, 4.73e+01)\\
3 & (2.12e+00, 6.99e+01) & (2.66e+00, 8.76e+01) & (2.25e+00, 7.48e+01)\\
4 & (2.54e+00, 1.03e+02) & (4.87e+00, 2.00e+02) & (2.68e+00, 1.08e+02)\\
5 & (2.95e+00, 1.36e+02) & (9.65e+00, 5.59e+02) & (3.16e+00, 1.44e+02)\\
6 & (3.29e+00, 1.78e+02) & (2.07e+01, 1.95e+03) & (3.51e+00, 1.85e+02)\\
7 & (3.63e+00, 2.21e+02) & (4.80e+01, 8.66e+03) & (3.89e+00, 2.30e+02)\\
8 & (3.91e+00, 2.72e+02) & (1.20e+02, 4.83e+04) & (4.18e+00, 2.80e+02)\\
9 & (4.20e+00, 3.26e+02) & (3.18e+02, 3.29e+05) & (4.50e+00, 3.35e+02)\\
10 & (4.45e+00, 3.86e+02) & (8.93e+02, 2.62e+06) & (4.75e+00, 3.95e+02)\\
11 & (4.69e+00, 4.50e+02) & (2.62e+03, 2.35e+07) & (5.02e+00, 4.59e+02)\\
12 & (4.91e+00, 5.20e+02) & (7.98e+03, 2.28e+08) & (5.25e+00, 5.29e+02)\\
13 & (5.13e+00, 5.94e+02) & (2.50e+04, 2.34e+09) & (5.49e+00, 6.03e+02)\\
14 & (5.33e+00, 6.74e+02) & (8.02e+04, 2.52e+10) & (5.69e+00, 6.83e+02)\\
15 & (5.53e+00, 7.58e+02) & (2.63e+05, 2.80e+11) & (5.91e+00, 7.67e+02)\\
16 & (5.71e+00, 8.48e+02) & (8.72e+05, 3.20e+12) & (6.09e+00, 8.56e+02)\\
17 & (5.89e+00, 9.42e+02) & (2.94e+06, 3.76e+13) & (6.29e+00, 9.50e+02)\\
18 & (6.06e+00, 1.04e+03) & (1.01e+07, 4.51e+14) & (6.46e+00, 1.05e+03)\\
19 & (6.23e+00, 1.15e+03) & (3.46e+07, 6.58e+15) & (6.63e+00, 1.15e+03)\\
20 & (6.38e+00, 1.25e+03) & (1.20e+08, 7.29e+16) & (6.79e+00, 1.26e+03)\\
21 & (6.53e+00, 1.37e+03) & (4.23e+08, 9.36e+18) & (6.96e+00, 1.38e+03)\\
22 & (6.67e+00, 1.49e+03) & (1.49e+09, 5.20e+19) & (7.10e+00, 1.49e+03)\\
23 & (6.76e+00, 1.61e+03) & (5.29e+09, 1.55e+20) & (7.26e+00, 1.62e+03)

%% file: _tab/tri.tab
1 & (1.00e+00, 4.00e+00) & (1.00e+00, 4.00e+00) & (1.00e+00, 4.00e+00) & (1.00e+00, 4.00e+00)\\
2 & (1.67e+00, 1.72e+01) & (1.67e+00, 1.72e+01) & (1.67e+00, 1.72e+01) & (1.67e+00, 1.72e+01)\\
3 & (2.12e+00, 3.47e+01) & (2.27e+00, 3.40e+01) & (2.11e+00, 3.48e+01) & (2.11e+00, 3.48e+01)\\
4 & (2.68e+00, 4.73e+01) & (3.47e+00, 5.77e+01) & (2.68e+00, 4.70e+01) & (2.66e+00, 4.59e+01)\\
5 & (3.27e+00, 6.51e+01) & (5.45e+00, 1.03e+02) & (3.41e+00, 6.91e+01) & (3.12e+00, 6.15e+01)\\
6 & (3.77e+00, 9.09e+01) & (8.75e+00, 2.15e+02) & (3.90e+00, 9.22e+01) & (3.70e+00, 9.20e+01)\\
7 & (4.39e+00, 1.21e+02) & (1.43e+01, 4.95e+02) & (4.48e+00, 1.27e+02) & (4.27e+00, 1.25e+02)\\
8 & (5.11e+00, 1.86e+02) & (2.40e+01, 1.27e+03) & (5.10e+00, 1.95e+02) & (4.96e+00, 1.93e+02)\\
9 & (5.95e+00, 2.73e+02) & (4.09e+01, 3.59e+03) & (5.87e+00, 2.90e+02) & (5.74e+00, 2.85e+02)\\
10 & (7.11e+00, 4.39e+02) & (7.09e+01, 1.09e+04) & (6.77e+00, 4.70e+02) & (6.67e+00, 4.70e+02)\\
11 & (8.39e+00, 7.02e+02) & (1.24e+02, 3.48e+04) & (8.04e+00, 7.49e+02) & (7.90e+00, 7.51e+02)\\
12 & (1.01e+01, 1.18e+03) & (2.21e+02, 1.19e+05) & (9.49e+00, 1.28e+03) & (9.36e+00, 1.31e+03)\\
13 & (1.21e+01, 1.98e+03) & (3.98e+02, 3.93e+05) & (1.17e+01, 2.17e+03) & (1.15e+01, 2.23e+03)\\
14 & (1.48e+01, 3.39e+03) & (7.20e+02, 1.44e+06) & (1.43e+01, 3.89e+03) & (1.40e+01, 4.06e+03)\\
15 & (1.79e+01, 5.80e+03) & (1.32e+03, 4.82e+06) & (1.80e+01, 6.98e+03) & (1.76e+01, 7.34e+03)\\
16 & (2.28e+01, 1.03e+04) & (2.42e+03, 1.84e+07) & (2.28e+01, 1.30e+04) & (6.08e+01, 2.90e+04)\\
17 & (2.90e+01, 1.99e+04) & (4.46e+03, 6.21e+07) & (2.97e+01, 2.44e+04) & (8.77e+01, 6.42e+04)\\
18 & (3.86e+01, 4.54e+04) & (8.34e+03, 2.43e+08) & (3.87e+01, 4.72e+04) & (1.29e+02, 1.48e+05)\\
19 & (5.14e+01, 8.77e+04) & (1.61e+04, 8.53e+08) & (5.16e+01, 9.13e+04) & (1.92e+02, 3.44e+05)\\
20 & (7.04e+01, 1.73e+05) & (5.43e+04, 8.35e+09) & (7.87e+01, 2.09e+05) & (3.15e+02, 9.47e+05)\\
21 & (1.27e+02, 3.96e+05) & (1.64e+05, 3.56e+11) & (1.92e+02, 1.47e+06) & (8.74e+02, 8.69e+06)\\
22 & (2.28e+02, 2.13e+06) & (7.22e+07, 9.66e+16) & (3.13e+02, 3.28e+06) & (9.27e+03, 1.15e+09)\\
23 & (6.80e+02, 1.59e+07) & (1.30e+06, 1.57e+13) & (7.92e+03, 2.61e+09) & (4.52e+04, 2.09e+10)

%% file: _tab/hex.tab
1 & (1.00e+00, 2.70e+01) & (1.00e+00, 2.70e+01) & (1.00e+00, 2.70e+01)\\
2 & (1.95e+00, 3.26e+02) & (1.95e+00, 3.26e+02) & (1.95e+00, 3.26e+02)\\
3 & (3.08e+00, 5.85e+02) & (4.34e+00, 8.20e+02) & (3.38e+00, 6.48e+02)\\
4 & (4.04e+00, 1.05e+03) & (1.08e+01, 2.83e+03) & (4.38e+00, 1.13e+03)\\
5 & (5.08e+00, 1.59e+03) & (2.99e+01, 1.32e+04) & (5.63e+00, 1.72e+03)\\
6 & (5.97e+00, 2.37e+03) & (9.41e+01, 8.64e+04) & (6.58e+00, 2.52e+03)\\
7 & (6.91e+00, 3.30e+03) & (3.33e+02, 8.05e+05) & (7.67e+00, 3.49e+03)\\
8 & (7.75e+00, 4.49e+03) & (1.31e+03, 1.06e+07) & (8.55e+00, 4.70e+03)\\
9 & (8.60e+00, 5.88e+03) & (5.69e+03, 1.89e+08) & (9.54e+00, 6.13e+03)

%% file: _tab/tet.tab
1 & (1.00e+00, 5.00e+00) & (1.00e+00, 5.00e+00) & (1.00e+00, 5.00e+00) & (1.00e+00, 5.00e+00)\\
2 & (2.00e+00, 3.60e+01) & (2.00e+00, 3.60e+01) & (2.00e+00, 3.60e+01) & (2.00e+00, 3.60e+01)\\
3 & (2.93e+00, 1.09e+02) & (3.02e+00, 1.05e+02) & (2.93e+00, 1.10e+02) & (2.93e+00, 1.10e+02)\\
4 & (4.09e+00, 2.46e+02) & (4.88e+00, 2.37e+02) & (4.09e+00, 2.50e+02) & (4.07e+00, 2.50e+02)\\
5 & (5.59e+00, 3.41e+02) & (8.09e+00, 4.10e+02) & (5.54e+00, 3.52e+02) & (5.32e+00, 3.62e+02)\\
6 & (7.38e+00, 6.42e+02) & (1.37e+01, 1.09e+03) & (7.17e+00, 6.68e+02) & (7.01e+00, 6.92e+02)\\
7 & (9.33e+00, 1.34e+03) & (2.34e+01, 2.74e+03) & (9.19e+00, 1.39e+03) & (9.20e+00, 1.39e+03)\\
8 & (1.23e+01, 3.04e+03) & (4.05e+01, 8.45e+03) & (1.21e+01, 3.13e+03) & (1.25e+01, 3.17e+03)\\
9 & (1.57e+01, 7.15e+03) & (7.11e+01, 2.63e+04) & (1.56e+01, 7.45e+03) & (1.70e+01, 7.80e+03)

%% file: _tab/triprism.tab
1 & (1.00e+00, 1.20e+01) & (1.00e+00, 1.20e+01) & (1.00e+00, 1.20e+01) & (1.00e+00, 1.20e+01) & (1.00e+00, 1.20e+01)\\
2 & (2.08e+00, 1.18e+02) & (2.08e+00, 1.18e+02) & (2.08e+00, 1.18e+02) & (2.08e+00, 1.18e+02) & (2.08e+00, 1.18e+02)\\
3 & (3.09e+00, 2.90e+02) & (3.70e+00, 3.18e+02) & (3.17e+00, 3.01e+02) & (3.17e+00, 3.01e+02) & (3.09e+00, 2.90e+02)\\
4 & (4.27e+00, 4.81e+02) & (7.67e+00, 8.16e+02) & (4.38e+00, 4.89e+02) & (4.35e+00, 4.78e+02) & (4.27e+00, 4.81e+02)\\
5 & (5.62e+00, 7.60e+02) & (1.69e+01, 2.44e+03) & (6.06e+00, 8.29e+02) & (5.55e+00, 7.38e+02) & (5.62e+00, 7.60e+02)\\
6 & (6.84e+00, 1.21e+03) & (3.98e+01, 9.50e+03) & (7.31e+00, 1.25e+03) & (6.93e+00, 1.25e+03) & (6.84e+00, 1.21e+03)\\
7 & (8.37e+00, 1.80e+03) & (9.93e+01, 4.60e+04) & (8.83e+00, 1.93e+03) & (8.43e+00, 1.89e+03) & (8.37e+00, 1.80e+03)\\
8 & (1.01e+01, 3.07e+03) & (2.62e+02, 2.79e+05) & (1.04e+01, 3.27e+03) & (1.01e+01, 3.23e+03) & (1.01e+01, 3.07e+03)\\
9 & (1.22e+01, 4.93e+03) & (7.28e+02, 2.06e+06) & (1.25e+01, 5.31e+03) & (1.22e+01, 5.23e+03) & (1.22e+01, 4.93e+03)

%% file: _tab/quadpyrmd.tab
1 & (1.00e+00, 1.22e+01) & (1.00e+00, 1.22e+01) & (1.00e+00, 1.22e+01) & (1.00e+00, -)\\
2 & (1.84e+00, 7.70e+01) & (1.84e+00, 7.70e+01) & (1.84e+00, 7.70e+01) & (2.08e+00, -)\\
3 & (2.76e+00, 2.58e+02) & (3.15e+00, 2.51e+02) & (2.83e+00, 2.70e+02) & (2.73e+00, -)\\
4 & (3.86e+00, 3.84e+02) & (5.94e+00, 4.90e+02) & (4.29e+00, 4.23e+02) & (4.13e+00, -)\\
5 & (4.95e+00, 6.12e+02) & (1.19e+01, 1.21e+03) & (6.84e+00, 8.81e+02) & (5.53e+00, -)\\
6 & (6.44e+00, 1.08e+03) & (2.51e+01, 3.70e+03) & (1.01e+01, 1.46e+03) & (7.35e+00, -)\\
7 & (8.06e+00, 1.99e+03) & (5.66e+01, 1.52e+04) & (1.42e+01, 2.81e+03) & (9.71e+00, -)\\
8 & (1.05e+01, 3.75e+03) & (1.36e+02, 9.10e+04) & (2.04e+01, 6.54e+03) & (1.28e+01, -)\\
9 & (1.35e+01, 7.19e+03) & (3.50e+02, 6.56e+05) & (3.11e+01, 1.72e+04) & (1.72e+01, -)